\documentclass[12pt, a4paper]{article}

\usepackage[latin1]{inputenc}
\usepackage[english]{babel}
\usepackage{t1enc}

 \usepackage{a4wide}
\usepackage{amsmath, amstext,amscd}
 \usepackage{graphics}
\usepackage{epic}
\usepackage{amssymb}

\def\wt{\widetilde{\tau}}
\def\Xr{\mathcal{X}}

\def\Z{\mathbb{Z}}

\def\C{\mathbb{C}}
\def\P{ \mbox{P\hspace{-.3em}I}}

\def\O{\mathcal{O} }
\def\P{\mathbb{P}}

\makeatletter

\@addtoreset{equation}{section} \makeatother

\newtheorem{prop}{Proposition}[section]

\newtheorem{lem}[prop]{Lemma}
\newtheorem{theo}[prop]{Theorem}
\newtheorem{cor}[prop]{Corollary}

\def\Rq{\stepcounter{prop} \noindent \theprop{.} \emph{Remark.}  }

\def\dem{\underline{Proof} : }

\def\tens{\otimes }

\def\Frac#1#2{{\displaystyle{{#1} \overwithdelims.. {#2}}}}

\begin{document}

\title{The action of the Frobenius map on rank $2$ vector bundles
over genus $2$ curves in small characteristics}
\author{Laurent Ducrohet}
\maketitle

\begin{center} {\bf \underline{Abstract}} \end{center}

Let $X$ be genus 2 curve  defined over an algebraically closed
field of characteristic $p$ and let $X_1$ be its $p$-twist. Let
$M_X$ (resp. $M_{X_1}$) be the (coarse) moduli space of
semi-stable rank 2 vector bundles with trivial determinant over
$X$ (resp. $X_1$). The moduli space $M_X$ is isomorphic to the 3
dimensional projective space and is endowed with an action of the
group $J[2]$ of order 2 line bundles over $X$. When $3\leq p \leq
7$, we show that the Verschiebung (i.e., the separable part of the
action of Frobenius by pull-back) $V : M_{X_1} \dashrightarrow
M_X$ is completely determined by its restrictions to the lines
that are invariant under the action of a non zero element of
$J[2]$. As those lines correspond to elliptic curves that appear
as Prym varieties, the Verschiebung restricts to the morphism
induced by multiplication by $p$ and we are able to compute the
explicit equations for the Verschiebung.

\section{Introduction}

Let $k$ be a algebraically closed field of positive characteristic
$p$ and let $X$ be  a proper and smooth (connected) curve of genus
2 over $k$. Let $X_s$ ($s\in \Z$ since $k$ is perfect) be the
$p^s$-twist of $X$ and let $J$ (resp. $J_s$) denote its Jacobian
variety (resp. the $p^s$-twist of its Jacobian variety). Denote by
$M_X(r)$ the (coarse) moduli space of semi-stable rank $r$ vector
bundles with trivial determinant over $X$. The map $E\mapsto
F_{\rm abs}^*\, E$ defines a rational map the separable part of
which, the generalized Verschiebung,  will be denoted by $V_r :
M_{X_1}(r) \dashrightarrow M_{X}(r)$.

For $r=2$ and $k=\C$, we let $M_X$ (resp. $V$) be $M_X$ (resp.
$V_2$). [NR] have constructed an isomorphism $D : M_X \to
|2\Theta|\cong \P^3$ that remains valid for an algebraically closed
field of positive characteristic (it is straightforward for $p\neq
2$ and [LP1] (section 5) give a sketch of proof for $p=2$).
Furthermore, the semistable boundary of $M_X$ identifies with the
Kummer surface ${\rm Kum}_X$, which is canonically contained in the
linear system $|2\Theta|$.  We have the commutative diagram
\begin{eqnarray} \end{eqnarray}
\begin{center} \unitlength=0.6cm
\begin{picture}(10,2)
\put(1.3,3){$M_{X_1}$}  \put(8,3){$M_{X}$}
\put(1.2,0.5){$|2\Theta_1|$}  \put(8.2,0.5){$|2\Theta|$}
\put(7.6,3.2){\vector(1,0){0.1}} \dashline{0.2}(3,3.2)(7.6,3.2)
\put(7.6,0.7){\vector(1,0){0.1}} \dashline{0.2}(3,0.7)(7.6,0.7)
\put(1.9,2.8){\vector(0,-1){1.5}} \put(8.7,2.8){\vector(0,-1){1.5}}
\footnotesize{\put(1.3,1.8){$D$} \put(8.9,1.8){$D$}
\put(5.2,3.4){$V$} \put(5.2,0.9){$\widetilde{V}$}}
\end{picture} \end{center}
and the induced rational map $\widetilde{V}$ is given by  degree
$p$ polynomials ([LP2]). Note that all the map in that diagram are equivariant under the action of $J[2]$.\\

Our interest in the situation described by the diagram (1.1) comes
from the fact ([LS]) that, given a proper and smooth curve of genus
$g$ over a field $k$, a semistable rank $r$ vector bundle $E$ over
$X$ corresponds to an (irreducible) continuous representation of the
algebraic fundamental group $\pi_1(X)$ in ${\rm GL}_r(\bar{k})$
(endowed with the discrete topology) if and only if one can find an
integer $n>0$ such that ${F_{\rm abs}^{(n)}}^*E\cong E$. Thus,
natural questions about the generalized Verschiebung $V_r :
M_{X_1}(r) \dashrightarrow M_{X}(r)$ arise like, e.g., its
surjectivity, its degree, the density of Frobenius-stable bundles,
and the loci of Frobenius-destabilized bundles.

For general $(g,\,r,\,p)$, not much seems to be known (see the
introductions of [LP1] and [LP2] for an overview of this subject).\\

When $g=2$, $r=2$, and $p=2$ and $X$ is an ordinary curve, [LP1]
determined the quadric equations of $\widetilde{V}$ in terms of the
generalized theta constants of the curve $X$. In [LP2] (resp. in
[Du]), one could give the equations of $\widetilde{V}$ in case of a
nonordinary curve $X$ with Hasse-Witt invariant equal to $1$ (resp.
a supersingular curve $X$) by specializing a family $\Xr$ of genus
$2$ curves parameterized by a discrete valuation ring with ordinary
generic fiber and special fiber isomorphic to $X$.

When $g=2$, $r=2$, and $p=3$, [LP2] determined the cubic equations
of $\widetilde{V}$ in showing that this rational map coincided with
the polar map of a Kummer surface, isomorphic to ${\rm Kum}_X$, in
$\P^3$. To reach this result, they used a striking relationship (see
[vG]) between cubics and quartics on $|2\Theta_1|$.

In this paper, we shall suppose that $p\geq 3$. Given $\tau$ a line
bundle of order 2 over $X$ and a $\tau$-invariant semi-stable vector
bundle $E$, i.e., satisfying $E\tens \tau \xrightarrow \sim E$, of
degree 0, one can give $E$ a structure of invertible $\O_X\oplus
\tau$-module. In other words, if $\pi : \widetilde{X} \to X$ is the
degree 2 étale cover corresponding to $\tau$, there is a degree 0
line bundle $q$ over $\widetilde{X}$ such that $E\cong \pi_*(q)$. On
the one hand, the $\tau$-invariant locus of $M_X$ is the union of
two projective lines. On the other hand, $\pi$ being étale, one has
$F_{\rm abs}^*(\pi_*(q))\cong \pi_*(F_{\rm abs}^*(q))$. Requiring
that $E$ has trivial determinant forces $q$ to be in some translate
of the Prym variety $P$ associated to $\pi$ (which is an elliptic
curve) and, as multiplication by $p$ over an elliptic curve commutes
with the inversion, it induces a map $P/\{\pm\} \cong \P^1 \to
P/\{\pm\} \cong \P^1$. Let $V_{\P^1}$ be the separable part of the
latter map  and choose a $\tau$-invariant line $\Delta(\tau)$ in
$M_X$. There is a natural isomorphism $P/\{\pm\} \xrightarrow \sim
\Delta(\tau)$ and $V_{\P^1}$ coincides with the restriction of
$V : M_{X_1} \dashrightarrow M_X$ to $\Delta(\tau)$. The main result of this paper is the following theorem (4.16) :\\

\noindent {\bf Theorem } {\it Let $X$ be a smooth and proper curve
of genus 2, sufficiently general, over an algebraically closed field
of characteristic $p=3,\,5$ or 7. The generalized Verschiebung $V :
M_{X_1} \dashrightarrow M_X$ is completely determined by its
restriction to the projective lines that are invariant under the
action of a non zero element of $J[2]$. }\\

In particular, $V$ can be computed explicitly in these cases and
we recover the result of [LP2] in the characteristic 3 case. As an
application, we recover if $p=3$, and we show if $p=5$, that
 there is a degree $2p-2$
surface $S$ is $|2\Theta_1|$ such that the equality of divisors in
$|2\Theta_1|$
$$\widetilde{V}^{-1}({\rm Kum}_{X})= {\rm Kum}_{X_1}+2 S $$ holds scheme-theoretically.  The computations have
been carried out using Maple 9 and Magma.

I would like to thank Y. Laszlo for having introduced me to this
question, for his help and encouragements all along my thesis. I
would like to thank D. Bernardi, P.-V. Koseleff, M. Chardin and more
specially G. Lecerf, for their help and explanations in the use of
the computation softwares.\\

\section{The (coarse) moduli space $M_X$}

\subsection{Preliminaries}

Let $k$ be an algebraically closed field of odd characteristic. Let
$X$ be a proper and smooth curve of genus 2 over $k$ and let $J$
(resp. $J^1$) be its Jacobian variety (resp. the moduli space of
degree one line bundles over $X$, which is a principal homogeneous
space under the action $J$). The embedding $X \hookrightarrow J^1$
defined by $\O(\Delta)$, where $\Delta$ is the diagonal $\Delta
\subseteq X\times X$, defines a theta divisor $\Theta$.

Let $M_X$ be the (coarse) moduli space of semi-stable rank 2 vector
bundles with trivial determinant over $X$. The set of its closed
points is the set of S-equivalence classes $[E]$ of semi-stable rank
2 vector bundles $E$ with trivial determinant over $X$. Note that
$M_X$ is endowed with a natural action of $J[2]$ defined
set-theoretically by $(\tau,\,[E]) \mapsto [E \tens \tau]$. If $E$
is strictly semi-stable, i.e., semi-stable but non-stable, there is
a degree 0 line bundle $j$ such that $E$ is an extension $0 \to
j^{-1} \to E \to j \to 0$, i.e., $E$ is S-equivalent to $j\oplus
j^{-1}$. If $\mathcal{L}$ is the Poincaré's bundle over $X\times J$,
one can consider the family $\mathcal{E}=\mathcal{L}\oplus
\mathcal{L}^{-1}$, parameterized by $J$, of semi-stable rank 2
bundles with trivial determinant over $X$. Because of the coarse
moduli property, there is a unique $b : J \to M_X$ such that, for
any $j$ in $J$, $b(j)=[\mathcal{E}_{|X\times \{j\}}]=[j\oplus
j^{-1}]$. Note that the semi-stable boundary in $M_X$, i.e., the
image of $b$, is (globally) $J[2]$-invariant and
that $b$ is $J[2]$-equivariant.\\

When $k=\C$ (but this result extends straightforwardly to any
algebraically closed
 field with characteristic different from 2, and, with a little more work (see [LP1] for a sketch of proof)
 to an algebraically closed field of characteristic 2), it has been shown ([NR]) that there is a canonical
 isomorphism  $D : M_X \xrightarrow \sim \P H^0(J^1,\, \O(2\Theta))$ that
 maps a semi-stable bundle $E$ over $X$ to the reduced (hence linearly
equivalent to $2\Theta$) divisor over $J^1$  with support the set
$\{\xi\in J^1 / h^0(X,\,E\tens \xi) \geq 1\}$  and recall by the way
that any divisor linearly equivalent to $2\Theta$ is completely
determined by its support ([NR], Proposition 6.4). In particular,
this divisor is reducible if and only if it is of the form
$D[j\oplus j^{-1}]$. In other words, $D$ maps the semi-stable
boundary of $M_X$ onto the set of reducible divisors of $|2\Theta|$,
namely $D[j\oplus j^{-1}]=T_j^*\Theta+T_{j^{-1}}^*\Theta$, where
$T_j$ (resp. $T_{j^{-1}}$) is the isomorphism $J^1 \xrightarrow \sim
J^1$ that corresponds to $j$ (resp. $j^{-1}$) via the action of $J$
on $J^1$ .

 Choose once for
all an effective theta characteristic $\kappa_0$ of $X$ (that is to
say one of the 6 Weierstrass points of $X$) and consider the
corresponding isomorphism $J \xrightarrow \sim J^1$ given by
$j\mapsto j\tens \kappa_0$. We still denote by $\Theta$ the divisor
of $J$ obtained as the inverse image of $\Theta$ by means of this
isomorphism. Its support is the set $\{\zeta \in J/ H^0(\zeta\tens
\kappa_0)\geq 1)$ and using the Riemann-Roch theorem, $\Theta$ is
symmetric. It defines the canonical principal polarization on $J$.
Using this isomorphism, we obtain an isomorphism
$$D : M_X
\xrightarrow \sim \P H^0(J,\, \O(2\Theta))=|2\Theta|$$ that maps
$[E]$ to the (unique) divisor in $|2\Theta|$ with support $\{\zeta
\in J / h^0(E \tens \kappa_0 \tens \zeta)\geq 1\}$. Of course, the
semi-stable boundary is $M_X$ still coincides with
the set of reducible divisors in $|2\Theta|$.\\

Now, we recall the basic facts of the theory of theta group schemes
associated to an ample line bundle $L$ over an abelian variety $A$
over $k$, as presented in [Mu2], Section 1. Let $\mathcal{G}(L)$
(resp. $K(L)$) be the group scheme (resp. the finite group scheme)
such that, for any $k$-scheme $S$,
$$\mathcal{G}(L)(S)=\{(x,\,\gamma)|\,x\in A(S),\ \gamma : L
\xrightarrow \sim T_x^*L \}\ \ ({\rm resp.}\  K(L)(S)=\{x\in
A(S)|\,T_x^*L \cong L\})$$ We suppose that $K(L)$ is
reduced-reduced, in which case $L$ is said to be of separable type.
The theta group $\mathcal{G}(L)$ is a central extension
$$1 \to \mathbb{G}_m \to \mathcal{G}(L) \to K(L) \to 0$$
Note that $\mathcal{G}(L)$ has a natural action on $H^0(A,\,L)$ and
that this action has weight 1 in the sense that, being given
$\lambda$ in $\mathbb{G}_m$ and $s$ in $H^0(A,\,L)$,
$\lambda.s=\lambda {\rm Id} (s)$. More generally, if
$\mathcal{G}(L)$ acts on a vector space $W$, one says that this
action has weight $r$ if, being given $\lambda$ in $\mathbb{G}_m$
and $s$ in $H^0(A,\,L)$, one has $\lambda.s=\lambda^r {\rm Id} (s)$.
As an example, one can consider the space of $r$-symmetric powers
${\rm Sym}^r\,H^0(A,\,L)$ on which  $\mathcal{G}(L)$ has a natural
action of weight $r$.

\begin{theo}
The vector space $H^0(A,\,L)$ is the unique (up to isomorphism)
irreducible representation of weight 1 of $\mathcal{G}(L)$.
\end{theo}
\dem Combine the Proposition 3 and the Theorem 2 of [Mu2],
Section 1.\\

Being given $\alpha,\, \beta$ in $K(L)$, $\widetilde{\alpha},\,
\widetilde{\beta}$ in $\mathcal{G}(L)$ above $\alpha$ and $\beta$
respectively, the invertible scalar
$\widetilde{\alpha}\widetilde{\beta}\widetilde{\alpha}^{-1}\widetilde{\beta}^{-1}$
depends only on the choice of $\alpha$ and $\beta$ and it induces
a skew-symmetric bilinear form $e^L : K(L) \times K(L) \to
\mathbb{G}_m$ that is non degenerate because $\mathbb{G}_m$ is the
center of $\mathcal{G}(L)$. Now, let $\pi : A \to B$ be a
separable isogeny between abelian varieties and suppose that $L
\cong \pi^* M$ where $M$ is an ample line bundle over $B$. Because
$L$ descends on $B$, there is an action of $\ker \pi$ on $L$
covering the action of $\ker \pi$ on $A$, i.e., a lifting $\ker
\pi \xrightarrow \sim \widetilde{K} \subseteq \mathcal{G}(L)$ and
$M \cong (\pi_*L)^{\widetilde{K}}$. In particular, $\ker \pi$ must
be an isotropic subgroup of $K(L)$ with respect to $e^L$.
Furthermore, $\gamma$ in $\mathcal{G}(L)$ induces an isomorphism
$\pi_*(\gamma)$ of $\mathcal{G}$ if and only if $\gamma$ lies in
the centralizer $Z(\widetilde{K})$ of $\widetilde{K}$ in
$\mathcal{G}(L)$ and, when it is the case, $\pi_*(\gamma)$ depends
only on the class of $\gamma$ in $Z(\widetilde{K})/\widetilde{K}$.
More precisely,

\begin{theo} (1) Suppose that  $\pi : A \to B$ is a separable isogeny
between abelian varieties and let $L$ be an ample line bundle over
$A$. Then, there is a one-to-one correspondence between :
\begin{itemize}
\item[-] The set of isomorphism classes of line bundles $M$ over $B$
such that $L \cong \pi^*M$. Such an $M$ is necessarily ample.
\item[-] The set of homomorphisms $\ker \pi \to \mathcal{G}(L)$ lifting
the inclusion $\ker \pi \hookrightarrow A$. \end{itemize} (2) Being
given a homomorphism $\ker \pi \to \widetilde{K} \subseteq
\mathcal{G}(L)$ as in (1) and letting $M \cong
(\pi_*L)^{\widetilde{K}}$ be the corresponding ample line bundle
over $B$, there are canonical isomorphisms $$\mathcal{G}(M) \cong
Z(\widetilde{K})/\widetilde{K} \text{ and } K(M) \cong (\ker
\pi)^\perp/\ker \pi$$ where $(\ker \pi)^\perp =\{\alpha \in K(L)|\,
e^L(\alpha,\,\beta)=1 \  for\ any \ \beta \in \ker \pi\}$.\\
In particular, if $M$ is a principal polarization, $\ker \pi$ is a
maximal isotropic subgroup of $K(L)$.
\end{theo}
\dem (1) [Mu1], Section 23, Theorem 2. The fact that $M$ is
necessarily ample comes from [Mu1], Section 6, Proposition 6 and the
finiteness of $K(M)$.\\
(2) [Mu2], Section 1, Theorem 4 and [Mu1], Section 23, Theorem 4 for the last assertion.\\

Returning to our situation, consider the ample line bundle
$\O(2\Theta)$ over $J$ and the associated theta group scheme
$\mathcal{G}(\O(2\Theta))$. Its natural action on
$H^0(J,\,\O(2\Theta))$ induces an action of $J[2]$ onto $|2\Theta|$
and the morphism $${\rm K}_X : J \to |2\Theta|,\ j \mapsto
T_j^*\Theta + T_{-j}^*\Theta$$ which obviously factors through the
quotient of $J$ under the involution $j\mapsto j^{-1}$, is
$J[2]$-equivariant. The following lemma is an obvious consequence of
the definitions :

\begin{lem} The diagram
\begin{center}
\unitlength=0.6cm
\begin{picture}(8,4)
\put(0,2){$J$} \put(6,3.7){$M_X$} \put(6,0.3){$|2\Theta|$}
\put(0.7,2.4){\vector(4,1){5}} \put(0.7,2){\vector(4,-1){5}}
\put(6.6,3.2){\vector(0,-1){2}} {\footnotesize \put(3.1,3.2){$b$}
\put(2.7,0.7){${\rm K}_X$}  \put(7,1.9){$D$}}
\end{picture}
\end{center}
of $J[2]$-equivariant morphisms is commutative.\\
\end{lem}

Let $H$ be the hyperplan in $|2\Theta|$ consisting in those
effective divisors passing through the origin of $J$. The
corresponding line bundle is $\O_{\P^3}(1)$. We let $\Delta$ be the
reduced divisor over $M_X$ with support $D^{-1}(H)$. The associated
invertible sheaf $\O(\Delta)$ is isomorphic to $D^*(\O_{\P^3}(1))$
and one has $$H^0(M_X,\,\O(\Delta)) \cong
H^0(|2\Theta|,\,\O(1))\cong H^0(J,\,\O(2\Theta))$$ Furthermore,
$b^*(\O(\Delta)) \cong \O(2\Theta)$. It can be shown that
$\O(\Delta)$ is a determinant line bundle is the following sense :
Let $S$ be a $k$-scheme $S$ and let $\mathcal{E}$ be a rank 2 vector
bundle with trivial determinant over $X \times S$ such that
$\mathcal{E}(s):=\mathcal{E}_{|X\times \{s\}}$ is semi-stable for
any closed point $s$ in $S$. Because of the coarse moduli property,
there is a map $\varphi : S \to M_X$ such that, for any closed point
$s \in S$, one has $\varphi(s)=[\mathcal{E}(s)]$. Denote by $q$ the
second projection $X \times S \to S$ and suppose moreover that
$h^0(X,\,\mathcal{E}(s) \tens \kappa_0))=0$ (which implies that $h^1
(X,\,\mathcal{E}(s)\tens \kappa_0)=0$) for $s$ generic in $S$.
Therefore, $q_*(\mathcal{E} \boxtimes \kappa_0)=0$ and one can find
a locally free resolution
$$0 \to E_0 \xrightarrow u E_1 \to R^1q_*(\mathcal{E} \boxtimes
\kappa_0) \to 0$$ The map $u$ is generically bijective and the
section ${\rm det}(u)$ defines an effective divisor on $S$ that does
not depend on the choice of the projective resolution, and that
coincides with $\varphi^{-1}(\Delta)$ (notice that ${\rm
Supp}\,\varphi^{-1}(\Delta)=\{s\in S|\,h^1(X,\,\mathcal{E}(s)\tens
\kappa_0)\geq 1\}$). The associated line bundle $({\rm det} \,
R^1q_*(\mathcal{E}\boxtimes \kappa_0))^{-1}$
therefore coincides with $\varphi^*(\O(\Delta))$.\\

Let $\varphi_{2\Theta} : J \to \P
H^0(J,\,\O(2\Theta))^*=|2\Theta|^*$ be the canonical morphism
associated to $\O(2\Theta)$. It is of course $J[2]$-equivariant.
One has the following lemma due to Wirtinger.

\begin{lem}[Wirtinger] There is a non degenerate bilinear pairing
on the vector space $W^*:=H^0(J,\,\O(2\Theta))^*$ such that the
induced isomorphism $B_W : W^* \xrightarrow{\sim}  W$ makes the
following diagram commutative
\begin{center}
\unitlength=0.6cm
\begin{picture}(8,4)
\put(0,2){$J$} \put(6.2,3.7){$|2\Theta|^*$}
\put(6.3,0.3){$|2\Theta|$} \put(0.7,2.4){\vector(4,1){5}}
\put(0.7,2){\vector(4,-1){5}} \put(6.9,3.2){\vector(0,-1){2}}
{\footnotesize \put(2.9,3.4){$\varphi_{2\Theta}$}
\put(2.7,0.7){${\rm K}_X$} \put(7.2,1.9){$\P B_W$}}
\end{picture}
\end{center}
This pairing is well-defined up to a non zero constant.
\end{lem}
\dem : We only sketch the proof that can be found in ([Mu3]). The
key idea is to consider the isogeny $\xi : J \times J\to J\times J$
given by $(x,\,y) \mapsto (x+y,\,x-y)$, the kernel of which is the
finite group $J[2]$, embedded diagonally in $J\times J$. One proves
that the pullback $\xi ^*(\O(\Theta)\boxtimes \O(\Theta))$ of the
principal polarization for $J \times J$ is isomorphic to
$\O(2\Theta)\boxtimes \O(2\Theta)$. Hence, there is a unique maximal
level subgroup $\widetilde{H}$ of the Heisenberg group
$\mathcal{G}(\O(2\Theta)\boxtimes \O(2\Theta))$ such that the
$\widetilde{H}$-invariant part of $\xi_*(\O(2\Theta)\boxtimes
\O(2\Theta))$ is isomorphic to $\O(\Theta)\boxtimes \O(\Theta)$. If
$\theta$ is "the" non-zero element of $H^0(J,\,\O(\Theta))$,
$\xi^*(\theta \boxtimes \theta)$ is "the" non-zero
$\widetilde{H}$-invariant element of $H^0(J\times
J,\,\O(2\Theta)\boxtimes \O(2\Theta)) \cong W\tens W$ and the
corresponding bilinear form on $W^*$ is non degenerate because the
action of $\mathcal{G}(\O(2\Theta))$ over $W$ is
irreducible.$\square$\\

\Rq Note that Wirtinger's result remains true whenever $J$ is a
principally polarized abelian variety and $\Theta$ is a symmetric
representative for this polarization. In the sequel, we will
always identify $|2\Theta|^*=\P H^0(J,\,\O(2\Theta))^*$ and
$|2\Theta| =\P H^0(J,\,\O(2\Theta))$)
using the isomorphism $\P B_W$.\\

The diagram in the lemma above is $J[2]$-equivariant as well. Thus,
gathering these preliminary results, we obtain the following
proposition (also found in [B]):

\begin{prop} The diagram
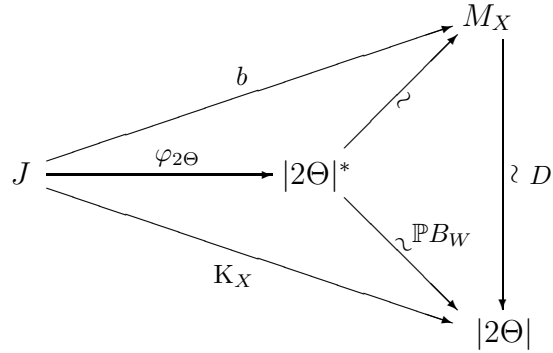
\begin{center}
\unitlength=0.6cm
\begin{picture}(15,8)
\put(0,4){$J$} \put(6,4){$|2\Theta|^*$}  \put(10,7.5){$M_X$}
\put(10.3,0.5){$|2\Theta|$} \put(10.9,7.2){\vector(0,-1){6}}
\put(7.4,3.7){\vector(1,-1){2.5}} \put(7.4,4.8){\vector(1,1){2.5}}
\put(0.8,4.2){\vector(1,0){5}} \put(0.8,4.5){\vector(3,1){9}}
\put(0.8,3.9){\vector(3,-1){9}} {\footnotesize \put(11.5,4){$D$}
\put(11.1,4){\rotatebox{90}{$\sim$}}
\put(3.2,4.5){$\varphi_{2\Theta}$}  \put(5,6.2){$b$}
\put(4.5,1.8){${\rm K}_X$} \put(8.5,2.6){\rotatebox{135}{$\sim$}}
\put(8.9,2.7){$\P B_W$} \put(8.5,5.6){\rotatebox{45}{$\sim$}}}
\end{picture}\end{center}
of $J[2]$-equivariant morphisms is commutative.
\end{prop}

\subsection{Choosing  a Theta structure}

Most of the material in that section comes from [Mu2], Section 1.\\

 Let us define the two groups
$$\left\{\begin{array}{l}
H=(\Z/2\Z)^2\\
\hat{H}=\text{Hom}((\Z/2\Z)^2,\,k^*)\end{array} \right.$$ We
identify $H$ and $\hat{H}$ by means of the bilinear form $(x,\,y)
\mapsto (-1)^{^tx.y}$ and the canonical evaluation map $H \times
\hat{H} \to k^*$ maps $(x,\,y^*)$ to $y^*(x)=(-1)^{^ty.x}$, where
$y$ corresponds to $y^*$ via the above mentioned identification
$\hat{H} \xrightarrow \sim H$. Define the Heisenberg group
$\mathcal{H}$ as the set $k^* \times H \times \hat{H}$, endowed with
the group structure given by
$$(t,\,x,\,x^*)(s,\,y,\,y^*)=(sty^*(x),\,x+y,\,x^*+y^*)$$
Denote by $E : (H\times \hat{H})\times (H\times \hat{H}) \to k^*$
the non degenerate bilinear form defined by the commutator in
$\mathcal{H}$, namely
$$E((x,\,x^*),\,(y,\,y^*))=[(1,\,x,\,x^*),\,(1,\,y,\,y^*)]=x^*(y)y^*(x)$$ It makes $H\times
\hat{H}$ self-dual, the sub-groups $H$ and $\hat{H}$ are isotropic
with respect to $E$ and the isomorphism of $H$ induced by the
restriction of $E$ to $(H \times \{0\})\times (\{0\}\times \hat{H})$
is the identity.

Let $e_2$ denote the bilinear pairing defined onto $J[2]$ by the
commutator in the theta group $\mathcal{G}(\O(2\Theta))$. As it is
non degenerate, it induces an isomorphism $$\phi : H\times \hat{H}
\xrightarrow \sim J[2]$$ symplectic with respect to $E$ and $e_2$,
namely a Göpel system for $J[2]$.

Because $H$ (resp. $\hat{H}$) is isotropic, one can choose a lifting
$H \xrightarrow \sim \widetilde{H}$ (resp. $\hat{H} \xrightarrow
\sim \widetilde{\hat{H}}$) in $\mathcal{G}(\O(2\Theta))$. There is a
unique isomorphism $$\widetilde{\phi} : \mathcal{H} \xrightarrow
\sim \mathcal{G}(\O(2 \Theta))$$ that induces identity on the
centers, that maps $H$ (resp. $\hat{H}$) onto $\widetilde{H}$ (resp.
$\widetilde{\hat{H}}$). Such an isomorphism is a theta structure on
$\mathcal{G}(\O(2\Theta))$ and we choose one once for all. In
particular, $W: = H^0(J,\,\O(2\Theta))$ is an irreducible
representation $U : \mathcal{H} \to {\rm GL}(W)$ of weight 1.

Given such a representation, one can construct a basis $\{X_x|\,
x\in H\}$, defined up to a multiplicative scalar, satisfying the
following properties :
$$y.X_x=X_{x+y} \text{ for any } x,\,y \in H,\ x^*.X_x=x^*(x)X_x \text{ for any } x,\,x^* \in
H\times \hat{H}.$$ Namely, the subgroup $\hat{H}\subset \mathcal{H}$
being abelian, the induced representation on $W$ is completely
reducible and splits in a direct sum of dimensional 1 subspaces
indexed by the group of character of $\hat{H}$, i.e., by $H$. Let
$X_{00}$ be a non-zero element of $W^{\hat{H}}$ and set
$X_x=(1,\,x,\,0).X_{00}$ for any $x\in H$. Then, for any $x^*\in
\hat{H}$, the equalities
$$(1,0,\,x^*)X_x=(1,0,\,x^*)(1,x,0)X_{00}=x^*(x)(1,\,x,0)(1,0,\,x^*)X_{00}=x^*(x)X_x$$
show that $X_x$ spans the subspace $W_x$ on which $(1,0,\,x^*)$ acts
like $x^*(x) Id$. Because the canonical evaluation map $(x,\,x^*)
\mapsto x^*(x)$ is non degenerate, the family $\{X_x|\, x\in H\}$ is
free and, since it contains 4 elements, we can conclude.

We call $\{X_{00},\,X_{01},\,X_{10},\,X_{11}\}$ the theta basis of
$W$ (associated to the theta structure $\phi : \mathcal{H}
\xrightarrow \sim \mathcal{G}(\O(2\Theta))$). In terms of that
basis, one can give an explicit description of $U : \mathcal{H}
\to {\rm GL}(W)$ : Set
$$\alpha_{01}=\left(\begin{array}{cccc}0 & 1 &0 & 0 \\
1& 0& 0& 0 \\
0& 0& 0& 1\\
0& 0& 1 &0
\end{array}\right); \hspace{1cm} \alpha_{10}=\left(\begin{array}{cccc}0 & 0 &1 & 0 \\
0& 0& 0& 1 \\
1& 0& 0& 0\\
0& 1& 0 &0
\end{array}\right)\hspace{1cm} \alpha_{11}=\left(\begin{array}{cccc}0 & 0 &0 & 1 \\
0& 0& 1& 0 \\
0& 1& 0& 0\\
1& 0& 0 &0
\end{array}\right)$$
and
$$ \beta_{01}=\left(\begin{array}{cccc}1 & 0 &0 &0 \\
0& -1& 0& 0 \\
0& 0& 1& 0\\
0& 0& 0 &-1
\end{array}\right); \hspace{1cm} \beta_{10}=\left(\begin{array}{cccc}1 & 0 &0 & 0 \\
0& 1& 0& 0 \\
0& 0& -1& 0\\
0& 0 & 0 &-1
\end{array}\right); \hspace{1cm} \beta_{11}=\left(\begin{array}{cccc}1 & 0 &0 & 0 \\
0& -1& 0& 0 \\
0& 0& -1& 0\\
0& 0 & 0 &1
\end{array}\right) $$
For convenience, set also $\alpha_{00}=\beta_{00}=Id$ and notice
that  $\alpha_{11}=\alpha_{01}\alpha_{10}=\alpha_{10}\alpha_{01}$
and $\beta_{11}=\beta_{01}\beta_{10}=\beta_{10}\beta_{01}$. Then
$U_{(t,\,x,\,x^*)}=t\beta_{x^*}\alpha_x$ for any $(t,\,x,\,x^*)\in
\mathcal{H}$. These six matrices are order 2 elements of ${\rm
GL}(W)$ with determinant 1 and the same statements hold for any
product $\beta_{x^*}\alpha_x$ with $(x,\,x^*)$ non-zero in
$H\times \hat{H}$.\\

\Rq Note that the isomorphism $B_W$ (Lemma 2.2) allows us to
identify $W$ and $W^*$. We let $\{x_\bullet\}$ be the
corresponding dual basis of $W^*$. It gives a set of homogeneous
coordinates for $\P W=|2\Theta|$, canonical in the sense that it
depends only on the
choice of the Göpel system $H\times \hat{H} \xrightarrow \sim J[2]$.\\

\subsection{The Kummer quartic surface}

Let ${\rm Kum}_X$ denote the image  of the map ${\rm K}_X: J \to
|2\Theta|$ appearing in Lemma 2.3. We have the following lemma

\begin{lem} Let $X$ be a genus 2 curve and let $J$ be its Jacobian. The map ${\rm K}_X : J \to |2\Theta|$ identifies with the quotient of $J$ under the
action of $\{\pm \}$. Its image is a reduced, irreducible,
$J[2]$-invariant quartic in $|2\Theta|$ with 16 nodes and no other
singularities, i.e., a \emph{Kummer surface}.
\end{lem}
\dem Because ${\rm K}_X$ coincides with $\varphi_{2\Theta}$ (Lemma
2.4) and because of the Riemann-Roch Theorem for abelian varieties
(see, e.g., [Mu1]), one has $$\deg (\varphi_{2\Theta}) . \deg ({\rm
Kum}_X) = (2\Theta)^2=8$$ Therefore, $\deg ({\rm Kum}_X)=4$. Because
$J$ is an irreducible abelian surface, ${\rm Kum}_X$ is reduced.
Thus, it is a quartic and ${\rm K}_X$ coincides with the quotient $J
\to J/\{\pm\}$ (see [GD], Proposition 4.23 for details). It is
therefore finite, surjective and separable, generically 2-1 hence of
degree 2 (see, e.g., [Mu1], Section 7) and ${\rm Kum}_X$ is
irreducible. The map ${\rm K_X}$ is furthermore generically \'etale
and it ramifies only at the $2$-torsion points of $J$ thus the
singular locus of ${\rm Kum}_X$ is contained in the image of $J[2]$.
An easy calculation in the formal completion of the local ring at a
point of $J[2]$ shows that it is a node (see [GD], Note 4.16). The
fact that it is
$J[2]$-invariant is clear. $\square$ \\

\Rq Let us consider a principally polarized abelian surface $A$ and
let $\Theta$ be a symmetric representative for that polarization.
Either $A$ is an irreducible abelian variety or $A$ is the product
of two elliptic curves. In the former case, one can show that
$\Theta$ is a non singular genus 2 curve the Jacobian of which is
isomorphic to $A$ and this is the situation of the lemma above. In
the latter case, $\Theta$ is the union of two elliptic curves
meeting transversally in one point and the morphism
$\varphi_{2\Theta} : A \to |2\Theta|^*$ associated to
 $\O(2\Theta)$ makes the following diagram commutative
\begin{center}
\unitlength=0.6cm
 \begin{picture}(10,4)
\put(0,1){\put(0,2.5){$E_1 \times E_2$}  \put(0.1,0){$\P^1 \times \P
^1$} \put(9,0){$|2\Theta|^* \cong \P^3$}
\put(1.1,2.3){\vector(0,-1){1.5}} \put(2.5,0.2){\vector(1,0){6}}
\put(2.5,2.4){\vector(3,-1){6}} \small
\put(4.5,2){$\varphi_{2\Theta}$}}\small \put(3,0.5){Segre embedding}
 \end{picture}
 \end{center}
where the map $E_1 \times E_2 \to \P^1\times \P ^1$ is the product
of the canonical maps $E_i \to \P^1$ ($i=1,\,2$), which identifies
with the quotient $E_i \to E_i/\{\pm\}$. The image of
the Segre embedding is a non singular quadric in $\P^3$.\\

\Rq Conversely, it corresponds to any quartic surface $S$ in $\P^3$
with the properties described in the lemma (i.e., a Kummer surface)
a principally polarized abelian surface $A$ whose principal
polarization is a non-singular curve ([GD], Proposition 4.22), hence
a genus 2 curve the Jacobian of which is isomorphic to $A$ (see the
previous Remark). In other words, one recovers the moduli space of
proper and smooth curves of genus 2 over $k$ (see [H], Chapter IV,
Ex. 2.2 for another (much more basic) description).\\

\begin{lem} (1) In the coordinate system $\{x_\bullet\}$ defined above,
there are scalars $k_{00},\,k_{01},\,k_{10}$ and $k_{11}$ such that
the equation defining the Kummer quartic surface ${\rm Kum}_X$ is
\begin{eqnarray}S+ 2k_{00}
P+k_{01}Q_{01}+k_{10}Q_{10}+k_{11}Q_{11}\end{eqnarray} where
$$\begin{array}{c}S= x_{00}^4+x_{01}^4+x_{10}^4+x_{11}^4,\hspace{1cm}
P=x_{00}x_{01}x_{10}x_{11},\\
Q_{01}=x_{00}^2x_{01}^2+x_{10}^2x_{11}^2,\hspace{1cm}Q_{10}=x_{00}^2x_{10}^2+x_{01}^2x_{11}^2,\hspace{1cm}
Q_{11}=x_{00}^2x_{11}^2+x_{01}^2x_{10}^2.\end{array}$$ (2) These
scalars $k_{00},\,k_{01},\,k_{10}$ and $k_{11}$ satisfy the cubic
relationship
\begin{eqnarray}4+k_{01}k_{10}k_{11}-k_{01}^2-k_{10}^2-k_{11}^2+k_{00}^2=0\end{eqnarray}
and one has \begin{eqnarray} \left\{\begin{array}{l} k_{01}\neq \pm
2, \ k_{10}\neq \pm
2,\  k_{11}\neq \pm 2,\\
k_{01}+k_{10}+k_{11}+2\pm k_{00}\neq 0,\\ k_{01}+k_{10}-k_{11}-2\pm
k_{00} \neq 0,\\ k_{01}-k_{10}+k_{11}-2\pm k_{00} \neq 0,\\
-k_{01}+k_{10}+k_{11}-2\pm k_{00} \neq 0 \end{array}\right.
\end{eqnarray}\end{lem}
\dem (1) Using the fact that the Kummer surface in $|2\Theta|$ is
$J[2]$-invariant, the element of ${\rm Sym}^4\,W^*$ that generates
its ideal is invariant (up to scalar) under the natural action (of
weight $-4$) of $\mathcal{H}$ on $ {\rm Sym}^4\,W^*$. Thus ([GD],
Theorem 2.20), one can look for the equation of ${\rm Kum}_X$
under the form (2.1).\\
(2) We let
$(\vartheta_{00}:\,\vartheta_{01}:\,\vartheta_{10}:\,\vartheta_{11})$
be the homogeneous coordinates of the image $\varphi_{2\Theta}(0)$
of the origin of $J$ in $|2\Theta|^*$, i.e., the generalized theta
constants in classical terminology of theta functions. Note that
the orbit of that point under the action of $J[2]$ consists in 16
distinct points. Because these points are the points of a 16-6
configuration (see [GD], Section 1), the corresponding  conditions
on the $\vartheta_\bullet$ are equivalent to the following :
$$\begin{array}{lclcl}
\vartheta_{00}\vartheta_{01} \neq \pm \vartheta_{10}\vartheta_{11},&
& \vartheta_{00}\vartheta_{10} \neq \pm
\vartheta_{01}\vartheta_{11},& & \vartheta_{00}\vartheta_{11} \neq
\pm
\vartheta_{01}\vartheta_{10},\\
\vartheta_{00}^2+\vartheta_{01}^2 \neq
\vartheta_{10}^2+\vartheta_{11}^2,& &
\vartheta_{00}^2+\vartheta_{10}^2 \neq
\vartheta_{01}^2+\vartheta_{11}^2,& &
\vartheta_{00}^2+\vartheta_{11}^2 \neq
\vartheta_{01}^2+\vartheta_{10}^2,\\
\vartheta_{00}^2+\vartheta_{01}^2 +
\vartheta_{10}^2+\vartheta_{11}^2 \neq 0\\
\end{array}$$
The fact that $\varphi_{2\Theta}(0)$ is a node means that the four
partial derivatives of the equation (2.1) vanish at
$(\vartheta_{00}:\,\vartheta_{01}:\,\vartheta_{10}:\,\vartheta_{11})$.
It gives a linear system of four equations (depending on the
$\vartheta_\bullet$) that the $k_\bullet$ must satisfy. Solving the
system, we can express the $k_\bullet$ in terms of the
$\vartheta_\bullet$. Namely, one has
$$\begin{array}{lcl}
k_{01}=-\Frac{\vartheta_{00}^4+\vartheta_{01}^4-\vartheta_{10}^4-\vartheta_{11}^4}{\vartheta_{00}^2\vartheta_{01}^2-\vartheta_{10}^2\vartheta_{11}^2},&
&
k_{10}=-\Frac{\vartheta_{00}^4-\vartheta_{01}^4+\vartheta_{10}^4-\vartheta_{11}^4}{\vartheta_{00}^2\vartheta_{10}^2-\vartheta_{01}^2\vartheta_{11}^2},\\
k_{11}=-\Frac{\vartheta_{00}^4-\vartheta_{01}^4-\vartheta_{10}^4+\vartheta_{11}^4}{\vartheta_{00}^2\vartheta_{11}^2-\vartheta_{01}^2\vartheta_{10}^2},
\end{array}$$
and $k_{00}$ can then be computed directly from (2.1). A few more
calculations ([GD], Lemma 2.21) show that one can eliminate the
$\vartheta_\bullet$ from the expressions of the $k_\bullet$ to find
the cubic relationship
\begin{eqnarray*}4+k_{01}k_{10}k_{11}-k_{01}^2-k_{10}^2-k_{11}^2+k_{00}^2=0\end{eqnarray*}
 From this equation, we obtain the following four equations
 \begin{eqnarray*}
(k_{01}+2)(k_{10}+2)(k_{11}+2)&= &(k_{01}+k_{10}+k_{11}+2-k_{00})(k_{01}+k_{10}+k_{11}+2+k_{00})\\
(k_{01}-2)(k_{10}-2)(k_{11}+2)& = &(k_{01}+k_{10}-k_{11}-2-k_{00})(k_{01}+k_{10}-k_{11}-2+k_{00})\\
(k_{01}-2)(k_{10}+2)(k_{11}-2)& =&(k_{01}-k_{10}+k_{11}-2-k_{00})(k_{01}-k_{10}+k_{11}-2+k_{00})\\
(k_{01}+2)(k_{10}-2)(k_{11}-2)& =&
(-k_{01}+k_{10}+k_{11}-2-k_{00})(-k_{01}+k_{10}+k_{11}-2+k_{00})\end{eqnarray*}
Note that neither
$$k_{01}-2=-\left[  \Frac{(\vartheta_{00}^2+\vartheta_{01}^2
-\vartheta_{10}^2-\vartheta_{11}^2)(\vartheta_{00}^2+\vartheta_{01}^2+
\vartheta_{10}^2+\vartheta_{11}^2)}{\vartheta_{00}^2\vartheta_{01}^2-\vartheta_{10}^2\vartheta_{11}^2}\right]$$
nor
$$k_{01}+2=-\left[  \Frac{(\vartheta_{00}^2-\vartheta_{01}^2
+\vartheta_{10}^2-\vartheta_{11}^2)(\vartheta_{00}^2-\vartheta_{01}^2-
\vartheta_{10}^2+\vartheta_{11}^2)}{\vartheta_{00}^2\vartheta_{01}^2-\vartheta_{10}^2\vartheta_{11}^2}\right]$$
can be zero. By symmetry, one has $k_{10}\neq \pm 2$ and $k_{11}\neq
\pm 2$ as well. Together with the four equations deduced from (2.2),
one can conclude. $\square$\\

\subsection{Invariant lines in $M_X$}

We are interested in those linear subspaces of $|2\Theta|$ on which
a non zero $\tau$ in $J[2]$ acts like identity. Let $\tau=(x,\,x^*)$
be an order 2 element in $J[2]$.

A lifting
 $(t,\,x,\,x^*)$ of $\tau$ in $\mathcal{H}$  has order 2 if and only if $t^2=x^*(x)$. Choose once for all a square root $i$
of $-1$ in $k$ and let $\wt$ be $(\mu ,\,x,\,x^*)$ with $\mu=1$
(resp. $\mu=i$) if $x^*(x)=1$ (resp. $x^*(x)=-1$). The corresponding
element $\mu \beta_{x^*}\alpha_x$ in ${\rm GL}(W)$ has order 2 and
determinant $\mu^4=1$. As it cannot be the identity matrix (because
$\tau\neq 0$), $W$ splits in the direct sum
$$W =W^{\wt}\oplus W^{-\wt}$$ of two 2-dimensional spaces of
eigenvectors, associated to the eigenvalues $+1$ and $-1$ of $\wt$
respectively.\\
We construct a basis adapted to this decomposition, that will be
useful in the computations of the Section 4.
\begin{prop}  Let $\tau$ be any non-zero element
of $J[2]$ and let $\wt$ be the order 2 element of $\mathcal{H}$
defined above. There is a basis
$\{\Lambda_0(\tau),\,\Lambda_1(\tau),\,\bar{\Lambda}_0(\tau),\,\bar{\Lambda}_1(\tau)\}$
for $W$ that splits into bases
$\{\Lambda_0(\tau),\,\Lambda_1(\tau)\}$ and
$\{\bar{\Lambda}_0(\tau),\,\bar{\Lambda}_1(\tau)\}$ for $W^{\wt}$
and $W^{-\wt}$ respectively. Furthermore, one can find a theta
structure, i.e., an automorphism of Heisenberg groups $\rho :
\mathcal{H} \to \mathcal{H}$, mapping $(1,00,10)$ to $\wt$, such
that this basis coincides with the theta basis corresponding to
the representation $\mathcal{H} \xrightarrow \rho \mathcal{H}
\xrightarrow U {\rm GL}(W)$.\end{prop} \dem Denote again by $\wt$
the element $U_{\wt}$ of ${\rm GL}(W)$ and by $1$ the identity
matrix. As $\wt^2=1$, it is easily seen that
$$p_{\tau}^+=(\wt+1)/2\ \ {\rm and}\ \ p_{\tau}^-=(1-\wt)/2)$$ is rank 2
projectors of the linear space $W$, and that their images are
$W^{\wt}$ and $W^{-\wt}$ respectively. One can extract a basis of
$W^{\wt}$ (resp. $W^{-\wt}$) from the family
$\{p_{\tau}^+(X_\bullet)\}$ (resp. $\{p_{\tau}^-(X_\bullet)\}$). Let
us distinguish whether $x=00$ or not.

 If $x$ is zero, one has
$$p_{\tau}^+(X_z)=\Frac{1+x^*(z)}{2} X_z =\left|\begin{array}{ll} X_z & \text{ if } x^*(z)=1\\
0 & \text{ if } x^*(z)=-1\end{array}\right.$$ and
$$p_{\tau}^-(X_z)=\Frac{1-x^*(z)}{2} X_z=\left|\begin{array}{ll} 0 & \text{ if } x^*(z)=1\\
X_z& \text{ if } x^*(z)=-1\end{array}\right.$$ Therefore, we let
$\Lambda_\bullet(\tau)$ (resp. $ \bar{\Lambda}_\bullet(\tau)$) be
the non zero elements $p_{\tau}^+(X_\bullet)$ (resp.
$p_{\tau}^-(X_\bullet)$)) with corresponding lexical order. In
other words, we permute the elements of the basis $\{X_\bullet\}$.
Notice that one always has $\Lambda_0(\tau)=X_{00}$ and let
$\{z_1,\,z_2,\,z_3\}$ be the permutation of $\{01,\,10,\,11\}$
such that
$$\Lambda_{1}(\tau)=X_{z_1},\ \bar{\Lambda}_{0}(\tau)=X_{z_2},\ \bar{\Lambda}_{1}(\tau)=X_{z_3}$$
Because $z_1+z_2+z_3=00$, the bijection $\rho_H : H \to H $
well-defined by the conditions $\rho_H(00)=00$, $\rho_H(01)=z_1$ and
$\rho_H(10)=z_2$ is a group automorphism of $H$. Furthermore, there
is  a unique automorphism $\rho_{\hat{H}}$ of $\hat{H}$ such that
the automorphism $$H \times \hat{H} \xrightarrow{\rho_H \times
\rho_{\hat{H}}} H \times \hat{H}$$ is symplectic with respect to
$E$. The automorphism $\rho$ of $\mathcal{H}$ acting like $\rho_H
\times \rho_{\hat{H}}$ on $H \times \hat{H}$ and acting like 1 on
the centers is the one we are looking for. Indeed,
$X_{00}=\Lambda_0(\tau)$ is invariant under the action of $\hat{H}$
and one has
$$\rho_H(01).\Lambda_0(\tau)=\Lambda_1(\tau),\
\rho_H(10).\Lambda_0(\tau)=\bar{\Lambda}_0(\tau),\
\rho_H(11).\Lambda_0(\tau)=\bar{\Lambda}_1(\tau)$$ thus the basis
$\{\Lambda_0(\tau),\,\Lambda_1(\tau),\,\bar{\Lambda}_0(\tau),\,\bar{\Lambda}_1(\tau)\}$
is the theta basis corresponding to the representation
$\mathcal{H} \xrightarrow \rho \mathcal{H} \xrightarrow{U} {\rm
GL}(W)$.\\

If $x$ is non zero, $\wt=(\mu,\,x,\,x^*)$ with $\mu^2=x^*(x)$ and
one has
$$p_{\tau}^+(X_z)=\Frac{1}{2}(X_z+\mu x^*(x+z) X_{x+z}) $$
Thus, the elements of the family $\{p_{\tau}^+(X_\bullet)\}$ are
pairwise colinear. In particular, we find that
$p_{\tau}^+(X_x)=\mu p_{\tau}^+(X_{00})$. We let $\Lambda_0(\tau)$
be $p_{\tau}^+(X_{00})$ and we let $z_0$ be the first (for lexical
order) non-zero element of $H$
 such that
$$\{\Lambda_0(\tau),\,p_{\tau}^+(X_{z_0})\}$$ is a basis for
$W^{\wt}={\rm Im}\, p_{\tau}^+$. This ${z_0}$ is necessarily
different from 00 and $x$ (more precisely, ${z_0}=10$ if $x=01$ and
${z_0}=01$ in the two other cases). We set
$$\Lambda_1(\tau)=p_\tau^+(X_{z_0}),\
\bar{\Lambda}_0(\tau)=p_\tau^-(X_{00}),\
\bar{\Lambda}_1(\tau)=p_\tau^-(X_{z_0})$$ and we obtained the
announced basis of $W$. Taking $\widetilde{\gamma}=(t,\,y,\,y^*)$
in $\mathcal{H}$, one has
$$\widetilde{\gamma}\left(\Lambda_0(\tau)\right)=\Frac{1+x^*(y)y^*(x)\wt}{2}(\widetilde{\gamma}(X_{00}))=
\Frac{ty^*(y)}{2}(X_{y}+\mu (x^*+y^*)(x)X_{x+y})$$ If
$\widetilde{\gamma}\left(\Lambda_0(\tau)\right)=\Lambda_0(\tau)$,
then either $y=0$ or $y=x$. In the former case, one must have $t=1$
and $y^*(x)=1$. There is a unique non-zero element in $\hat{H}$ (say
$z_1^*$) fulfilling this condition and because $x$ and ${z_0}$
generate $H$, one must have $z_1^*({z_0})=-1$. In the latter case,
one must have $t=\mu$ and $y^*(x)=x^*(x)$. The two elements of
$\hat{H}$ fulfilling this condition are $x^*$ and $x^*+z_1^*$.\\
Consider the abelian subgroup
 $\{1,\,(1,00,\,z_1^*),\,\wt,\,(\mu,\,x,\,x^*+z_1^*)\}$ of
 $\mathcal{H}$.
By construction, it fixes $\Lambda_0(\tau)$. Still by
construction, one has
$$\wt (\Lambda_1(\tau))= \Lambda_1(\tau),\ \ \wt (\bar{\Lambda}_0(\tau))= -
\bar{\Lambda}_0(\tau),\ \ \wt (\bar{\Lambda}_1(\tau))=
-\bar{\Lambda}_1(\tau)$$ Finally, because
$E((00,\,z_1^*),\,\tau)=z_1^*(x)=1$ and $z_1^*({z_0})=-1$, one
obtains
$$ (1,00,\,z_1^*)(\Lambda_1(\tau))=
-\Lambda_1(\tau),\ \ (1,00,\,z_1^*) (\bar{\Lambda}_0(\tau))=
\bar{\Lambda}_0(\tau),\ \ (1,00,\,z_1^*)(\bar{\Lambda}_1(\tau))= -
\bar{\Lambda}_1(\tau)$$ Now, choose two elements $z_2^*$ and
$z_3^*$ in $\hat{H}$, not necessarily different, such that
$z_2^*(x)=x^*({z_0})$ and $z_3^*(x)=-1$. It is easily checked
that, on the one hand,
$$(z_2^*({z_0}),\,{z_0},\,z_2^*)(\Lambda_0(\tau))=z_2^*({z_0})^2\Frac{1+x^*({z_0})z_2^*(x)\wt}{2}X_{z_0}=\Lambda_1(\tau)$$
and
$$(1,00,\,z_3^*)(\Lambda_0(\tau))=\Frac{1+z_3^*(x)\wt}{2}X_{00}=\bar{\Lambda}_0(\tau)$$
and that, on the other hand,
$$[(z_2^*({z_0}),\,{z_0},\,z_2^*),\,(1,00,\,z_1^*)]=z_1^*({z_0})=-1,\
[(z_2^*({z_0}),\,{z_0},\,z_2^*),\,\wt]=x^*({z_0})z_2^*(x)=1$$ and
$$[(1,00,\,z_3^*),\,(1,00,\,z_1^*)]=1,\
[(1,00,\,z_3^*),\,\wt]=z_3^*(x)=-1$$ where $[.,.]$ is the
commutator is $\mathcal{H}$. Hence, one can define an automorphism
$\rho : \mathcal{H} \xrightarrow \sim \mathcal{H}$ by setting
$$\begin{array}{lcl} \rho((1,00,01))=(1,00,\,z_1^*), & &
\rho((1,00,10))=\wt, \\
\rho((1,01,00))=(z_2^*({z_0}),\,{z_0},\,z_2^*), & &
\rho((1,10,00))=(1,00,\,z_3^*),\end{array}$$ and by asking that
the map induced on the centers is 1. $\square$\\

\Rq The bases $\{X_{00},\,X_{01},\,X_{10},\,X_{11}\}$ and
$\{\Lambda_0(\tau),\,\Lambda_1(\tau),\,\bar{\Lambda}_0(\tau),\,\bar{\Lambda}_1(\tau)\}$
are the same for $\tau=(0010)$.

\begin{cor}
(1) Given $\tau$ in  $J[2] \setminus \{0\}$, there are two
(disjoint) $\tau$-invariant lines in $|2\Theta|$. The
$\tau$-invariant locus in $|2\Theta|$ is globally invariant under
the action of $J[2]$ and a element $\alpha$ in $J[2]$ permutes the
connected components if and only if $e_2(\alpha,\,\tau)=-1$.\\
(2) If $\Delta(\tau)$ is a line of (1), there are coordinates
$\{\lambda_0,\,\lambda_1\}$ such that $\Delta(\tau)\cap {\rm Kum}_X$
consists in four reduced points with homogeneous coordinates
$$(a:\,b),\ (a:\,-b),\ (b:\,a)\text{  and }(b:\,-a)$$ These scalars don't
depend on the choice of the $\tau$-invariant line $\Delta(\tau)$.
\end{cor}
\dem Because of the proposition, it is enough to prove the
corollary for a particular $\tau$ so let $\tau$ be the element
$(0010)$ of $J[2]$ and let $\wt=(1,\,00,\,10)$ in $\mathcal{H}$.
We will denote by $\Delta^+(\tau)$ (resp. $\Delta^-(\tau)$) the
projective line in $|2\Theta|$ corresponding to $W^{\wt}$ (resp.
$W^{-\wt}$). If $X$ is any eigenvector of $\wt$ and if
$\widetilde{\alpha}$ is any element in $\mathcal{H}$ with class
$\alpha$ in $H\times \hat{H}$, one has
$$\wt(\widetilde{\alpha} (X))=e_2(\alpha,\,\tau)=\widetilde{\alpha}
(\wt(X))$$ which proves (1).

 The ideal $(x_{10},\,x_{11})$ of $\Delta^+(\tau)$ coincides with the
kernel of the dual map $W^* \twoheadrightarrow (W^{\wt})^*$.
Denote by $\lambda_0$ (resp. $\lambda_1$) the image of $x_{00}$
(resp. $x_{01}$) in $(W^{\wt})^ *$. The equation (2.1) of ${\rm
Kum}_X$ maps to
$$ \lambda_0^4+\lambda_1^4+k_{01}\lambda_0^2\lambda_{1}^2=0$$
in ${\rm Sym}^4\,(W^{\wt})^*$ and this quartic generates the ideal
of the scheme-theoretic intersection $\Delta^+(\tau)\cap {\rm
Kum}_X$. Hence, because $k_{01}\neq \pm 2$ (Lemma 2.11.(2)),
$\Delta^+(\tau)\cap {\rm Kum}_X$ consists in four reduced points
with homogeneous coordinates
$$(a:\,b:\,0:\,0),\ \ (b:\,a:\,0:\,0),\ \ (a:\,-b:\,0:\,0),\ \ (b:-a:\,0:\,0).$$
where $a$ and $b$ are pairwise different non zero scalars
satisfying the equality $k_{01}=-\Frac{b^4+a^4}{a^2b^2}$.

Because ${\rm Kum}_X$ is $J[2]$-invariant and because
$\alpha_{10}$ maps $X_{00}$ to $X_{10}$ and $X_{01}$ to $X_{11}$,
it is clear that one would have similar results concerning
$\Delta^-(\tau)$. $\square$\\

We let $\omega(\tau)$ be the unique scalar such that the equation
of the Kummer surface restricts to \begin{eqnarray}
\lambda_0(\tau)^4+\lambda_1(\tau)^4+\omega(\tau)
\lambda_0(\tau)^2\lambda_{1}(\tau)^2=0
\end{eqnarray} on $\Delta(\tau)$. The equations (2.3) exactly tell us that $\omega(\tau)\neq \pm 2 $ for any $\tau$ in $J[2]\setminus
\{0\}$ and we gather their expression in terms of the $k_\bullet$,
computed thanks to the bases constructed in the Proposition 2.12,
in the following chart.

\begin{eqnarray}\begin{array}{|c|c|c|c|c|}

\hline x\setminus x^* & 00 & 01 & 10 & 11 \\
\hline 00 & \star & k_{10}  & k_{01}& k_{11} \\
\hline 01 & \Frac{2(k_{00}+k_{10}+k_{11})}{2+k_{01}} &
\Frac{2(-k_{00}+k_{10}-k_{11})}{2-k_{01}}&
\Frac{2(-k_{00}+k_{10}+k_{11})}{2+k_{01}} &
\Frac{2(k_{00}+k_{10}-k_{11})}{2-k_{01}}\\
\hline 10 & \Frac{2(k_{00}+k_{01}+k_{11})}{2+k_{10}} &
\Frac{2(-k_{00}+k_{01}+k_{11})}{2+k_{10}}&
\Frac{2(-k_{00}+k_{01}-k_{11})}{2-k_{10}} &
\Frac{2(k_{00}+k_{01}-k_{11})}{2-k_{10}}\\
\hline 11 & \Frac{2(k_{00}+k_{01}+k_{10})}{2+k_{11}} &
\Frac{2(k_{00}+k_{01}-k_{10})}{2-k_{11}}&
\Frac{2(-k_{00}+k_{01}-k_{10})}{2-k_{11}} &
\Frac{2(-k_{00}+k_{01}+k_{10})}{2+k_{11}}\\
\hline
\end{array}\hspace{0.3cm} \end{eqnarray} \\

\Rq Let $\Delta(\tau)$ be a $\tau$-invariant line in $|2\Theta|$. We
will prove in the next section that this line identifies with the
quotient $P/\{\pm\}$, where $P$ is the Prym variety associated to
$\tau$, which is an elliptic curve in that case. Furthermore, the
points of the intersection $\Delta(\tau) \cap {\rm Kum}_X$ are the
Weierstrass points of $P$. This gives another proof of the fact that
$\omega(\tau) \neq \pm 2$ for any $\tau$ in
$J[2]\setminus \{0\}$.\\

\section{Prym varieties}

Most of the material in that section has been adapted to our
situation from the very clear exposition found in [Mu3].\\

\subsection{Etale double cover}

Choose a non-zero element $\tau$ of $J[2]$. One can construct an
\'etale double cover $\pi : \widetilde{X} \to X$ as follows : An
isomorphism $\phi : \tau \tens \tau \xrightarrow \sim \O_X$ allows
us to give the direct sum $\O_X \oplus \tau$ a structure of
$\O_X$-algebra with product
$$(a,\,l)(b,\,m)=(ab+\phi(lm),\,am+bl)$$
and we set $\widetilde{X}:=\textbf{Spec}(\O_X\oplus \tau)$, which
is of genus 3 (Hurwitz). Of course, it does not depend (up to
isomorphism) on the choice of the isomorphism $\phi : \tau \tens
\tau \xrightarrow \sim \O_X$. Denote by $\widetilde{J}$ the
Jacobian of $\widetilde{X}$. If $\widetilde{J}\,^2$ stands for the
moduli space of degree 2 line bundles over $\widetilde{X}$, there
is a canonical theta divisor $\widetilde{\Theta} \subseteq
\widetilde{J}\,^2$. As $\pi$ is \'etale, Hurwitz's formula assures
that the canonical divisor $\widetilde{K}$ on $\widetilde{X}$
coincides with the pull-back $\pi^* K$ of the canonical divisor on
$X$. Therefore, $\pi^*\kappa_0$ is a theta characteristic on
$\widetilde{X}$ and pulling-back $\widetilde{\Theta}$ by the
isomorphism $\widetilde{J} \xrightarrow \sim \widetilde{J}^2$
defined by $j \mapsto j\tens \pi^*\kappa_0$), we obtain a
symmetric divisor (still denoted $\widetilde{\Theta}$) that
represents the canonical polarization of
$\widetilde{J}$.\\

We have homomorphisms $\pi^* : J \to \widetilde{J}$ and $\text{Nm}
: \widetilde{J} \to J$, the latter being deduced from the
push-forward of divisors via $\pi$. We easily check that the
composite $\text{Nm}.\pi^*$ is multiplication by 2. Therefore,
$\ker \pi^* \subseteq J[2]$ and it is easily seen to be equal to
$<\tau>$. Furthermore, using divisors, one can show (see, e.g.,
[H], Chapter IV, Ex. 2.6) that, for any $j\in \widetilde{J}$,
\begin{eqnarray}\text{det} (\pi_* j) \cong
\text{det}(\pi_*\O_{\widetilde{X}}) \tens \text{Nm}(j) \cong \tau
\tens \text{Nm}(j)\end{eqnarray} The homomorphisms $\pi^*$ and
$\text{Nm}$ are dual one to each other, i.e., the
following diagrams (equivalent by duality) commute\\

\begin{eqnarray} \end{eqnarray}
\begin{center}
\unitlength=0.6cm
\begin{picture}(17,1)
\put(1,2.5){$\widetilde{J}$} \put(6,2.5){$\widehat{\widetilde{J}}$}
\put(1,0){$J$}\put(6,0){$\widehat{J}$}
\put(1.6,2.7){\vector(1,0){4}} \put(1.6,0.2){\vector(1,0){4}}
\put(1.2,0.7){\vector(0,1){1.5}} \put(6.2,0.7){\vector(0,1){1.5}}
\put(3.4,2.7){$\sim$} \put(3.4,0.2){$\sim$}

\small \put(3.4,3.1){$\lambda_{\widetilde{\Theta}}$}
\put(1.4,1.3){$\pi^*$} \put(3.4,0.5){$\lambda_\Theta$}
\put(6.5,1.3){$\widehat{\text{Nm}}$}

\put(9,0){\put(1,2.5){$\widetilde{J}$}
\put(6,2.5){$\widehat{\widetilde{J}}$}
\put(1,0){$J$}\put(6,0){$\widehat{J}$}
\put(1.6,2.7){\vector(1,0){4}} \put(1.6,0.2){\vector(1,0){4}}
\put(1.2,2.3){\vector(0,-1){1.5}} \put(6.2,2.3){\vector(0,-1){1.5}}
\put(3.4,2.7){$\sim$} \put(3.4,0.2){$\sim$}

\small \put(3.4,3.1){$\lambda_{\widetilde{\Theta}}$}
\put(1.4,1.3){$\text{Nm}$} \put(3.4,0.5){$\lambda_\Theta$}
\put(6.5,1.3){$\widehat{\pi^*}$} }
\end{picture} \end{center}
where $\lambda_\Theta : J \xrightarrow \sim  \widehat{J}$ (resp.
$\lambda_{\widetilde{\Theta}} : \widetilde{J} \xrightarrow \sim
\widehat{\widetilde{J}}$) is the isomorphism of abelian varieties
defined by $j \mapsto \O(T_j^*\Theta - \Theta)$  (resp. $j \mapsto
\O(T_j^*\widetilde{\Theta}-\widetilde{\Theta})$). This implies that
the third diagram\\

\begin{eqnarray} \end{eqnarray}
\begin{center}
\unitlength=0.6cm
\begin{picture}(8,1)
\put(1,2.5){$\widetilde{J}$} \put(6,2.5){$\widehat{\widetilde{J}}$}
\put(1,0){$J$}\put(6,0){$\widehat{J}$}
\put(1.6,2.7){\vector(1,0){4}} \put(1.6,0.2){\vector(1,0){4}}
\put(1.2,0.7){\vector(0,1){1.5}} \put(6.2,2.3){\vector(0,-1){1.5}}
\put(3.4,2.7){$\sim$}

\small \put(3.4,3.1){$\lambda_{\widetilde{\Theta}}$}
\put(1.4,1.3){$\pi^*$} \put(3.4,0.5){$2\lambda_\Theta$}
\put(6.5,1.3){$\widehat{\pi^*}$}
\end{picture} \end{center}
is also commutative. One checks that, from a set-theoretical point
of view,
\begin{eqnarray*}
(\pi^*)^{-1}({\rm Supp}\ \widetilde{\Theta})&=&\{j\in J
|h^0(\widetilde{X},\,\pi^*(j\tens
\kappa_0)\geq 1\}\\
& = & \{j\in J |h^0(X,\,\pi_*\O_{\widetilde{X}} \tens j\tens
\kappa_0)\geq 1\}\\
& =& \{j\in J |h^0(X,\,(j\tens \kappa_0) \oplus (\tau \tens j\tens
\kappa_0))\geq 1\}={\rm Supp}\ \Theta \cup {\rm Supp}\ T_\tau ^*
\Theta
\end{eqnarray*}
Therefore, the divisor $(\pi^*)^{-1}(\widetilde{\Theta})$ is defined
and the commutativity of the latter diagram assures that it is
algebraically equivalent to $2\Theta$.\\

\Rq Note that this result holds for any representative of the
principal polarization $\lambda_{\widetilde{\Theta}}$, i.e., for any
$j$ in $\widetilde{J}$, as soon as the divisor
$(\pi^*)^{-1}(T_j^*\widetilde{\Theta})$ is defined, it is
algebraically equivalent
to $2\Theta$.\\

 Let us introduce the set $$S_\tau=\{z\in J/z^2=\tau\}
\subseteq J[4]$$ which is a principal homogeneous space under
$J[2]$. For any $z$ in $S_\tau$, $\pi^*(z)$ belongs to
$\widetilde{J}[2]$ for $\pi^*(z)^2=\pi^*(z^2)=\pi^*(\tau)=0$. Let
$\widetilde{\Theta}_z$ be the effective divisor
$T_{\pi^*(z)}^*\widetilde{\Theta}$, which is a symmetric
representative for the principal polarization
$\lambda_{\widetilde{\Theta}}$. From a set-theoretical point of
view,
\begin{eqnarray*}
(\pi^*)^{-1}({\rm Supp}\ \widetilde{\Theta}_z)&=&\{j\in J
|h^0(\widetilde{X},\,\pi^*(j\tens z \tens
\kappa_0)\geq 1\}\\
& =& \{j\in J |h^0(X,\,(j\tens z \tens \kappa_0) \oplus (z^{-1}
\tens j\tens \kappa_0))\geq 1\}\\
&  = & {\rm Supp}\  T_z^*\Theta \cup {\rm Supp}\ T_{-z} ^* \Theta
\end{eqnarray*}
The divisor $(\pi^*)^{-1}(\widetilde{\Theta}_z)$ is therefore
defined and algebraically equivalent to $2\Theta$. Note furthermore
that $\pi^*$ being a degree 2 map onto its image, we have, in terms
of divisors, the equality
$(\pi^*)^{-1}(\widetilde{\Theta}_z)=T_z^*\Theta +T_{-z} ^* \Theta$
and the latter is linearly equivalent to $2\Theta$.\\

\subsection{Prym varieties}

Denote by $P$ the abelian variety $\ker (\text{Nm} : \widetilde{J}
\to J)^0$, and choose $z$ in $S_\tau$. Denote by $\sigma$ the
homomorphism
$$J \times P \to \widetilde{J};\ (x,\,q) \mapsto \pi^*(x)+q$$ One
easily sees that $\sigma$ has finite kernel $K_\sigma \subseteq J[2]
\times P[2]$. Thus, it is a separable isogeny and $P$ is an abelian
variety of dimension 1, i.e., an elliptic curve.

\begin{lem} (1) The composite
$$J \times P \hspace{0.3cm}\xrightarrow{\sigma} \hspace{0.3cm}\widetilde{J}\hspace{0.3cm}
\xrightarrow{\lambda_{\widetilde{\Theta}_z}}
\hspace{0.3cm}\widehat{\widetilde{J}}
\hspace{0.3cm}\xrightarrow{\hat{\sigma}} \hspace{0.3cm} \hat{J}
\times \hat{P}$$ defines a polarization for $J\times P$ that splits
into a product $2\lambda_\Theta \times \rho$, where $\rho
: P \to \hat{P}$ is a polarization for $P$.\\
(2) There is a  isomorphism $$\varphi : \tau\,^\perp /<\tau>
\xrightarrow{\sim } P[2]$$ where $\tau\,^\perp=\{\alpha \in J[2]|\
e_2(\alpha,\,\tau)=1\}$, symplectic with respect to $e_{2,\,J}$ and
$e_{2,\,P}$.\\
In particular, $\rho=2\lambda_P$, where $\lambda_P : P \xrightarrow
\sim  \hat{P}$ is a principal polarization.
\end{lem}
\dem The results in that lemma can be found in [Mu3], Sections 2 and 3. We only sketch Mumford's arguments. \\
(1) The polarization defined in the lemma
 may be viewed as a $2\times 2$ matrix
$\left(\begin{array}{cc} \alpha & \beta \\ \gamma & \rho
\end{array} \right)$, where $\alpha$, $\beta$, $\gamma$ and $\rho$ are the
expected homomorphisms of abelian varieties. Because of the diagram
(3.3) and the Remark 3.1 in the previous subsection,
$\alpha=2\lambda_\Theta$. Because any polarization is symmetric,
$\beta=\hat{\gamma}$ and because $\widetilde{J}
\xrightarrow{\lambda_{\widetilde{\Theta}_z}}
\hat{\widetilde{J}}\xrightarrow{\widehat{\pi^*}} \hat{J}$ is zero by
the very definition of $P$ (see the diagram (3.2)), $\beta=0$.
Therefore, the considered polarization splits into the product
$2\lambda_\Theta \times \rho$, where $\rho : P \to \hat{P}$ is a
polarization for
$P$.\\
(2) On the one hand, because $K_\sigma \cap (\{0\} \times
P[2])=\{0\}$, one can find a subgroup $K \subset J[2]$, containing
$\ker \pi^*= <\tau>$, and an injective homomorphism $\varphi :
K/<\tau> \to P[2]$ such that $K_\sigma$ identifies with $K$ by
means of the map $$K \xrightarrow \sim K_\sigma, \ \alpha \mapsto
(\alpha,\,\varphi(\bar{\alpha}))$$where $\bar{\alpha}$ is the
class of $\alpha$ via the quotient map $H \twoheadrightarrow
H/<\tau>$. It is clear that
\begin{enumerate}
\item[(a)] $\varphi(K/<\tau>)
\subset \ker \rho \subset P[2]$.\end{enumerate} On the other hand,
because $\lambda_{\widetilde{\Theta}_z}$ is a principal
polarization, $K_\sigma$ is a maximal isotropic subgroup of $\ker(2
\lambda_\Theta \times \rho)\cong J[2] \times \ker (\rho)$, with
respect to the product bilinear form $e_{2,\,J} \times e_\rho$ (see
Theorem 2.2.(2)). Thus,
\begin{enumerate}
\item[(b)] $\text{card}(K)^2=\text{card}(K_\sigma)^2 =
\text{card}(J[2]) \text{card}(\ker \rho)$,
\item[(c)] for any $\alpha,\,\beta$ in $K$,
$e_{2,\,X}(\alpha,\,\beta)e_\rho(\varphi(\bar{\alpha
}),\,\varphi(\bar{\beta }))=1$.\end{enumerate} Use (c) to obtain the
inclusion  $<\tau> \subseteq K^\perp$ and combine (a)  and (b) to
show that it is an equality. Thus, because $P$ is an elliptic curve,
the inclusions of (a) are equalities and (b) assures that the
corresponding isomorphism is symplectic.\\
It ends the proof of the lemma since the fact that $\ker
(\rho)=P[2]$ implies that $\rho$ is twice the principal polarization
of $P$. $\square$\\

Choose a symmetric divisor $\Xi$ representing the principal
polarization $\lambda_P$. The line bundle $\O(2\Xi)$ is canonical
and because the polarization defined above splits, one has
$$\sigma^*(\O(\widetilde{\Theta}_z))\cong \O(2\Theta) \boxtimes
\O(2\Xi)$$ Consider the Heisenberg group $\mathcal{G}(\O(2\Theta)
\boxtimes \O(2\Xi))$ associated to that ample line bundle over $J
\times P$. It is a central extension
$$ 1 \to \mathbb{G}_m \to \mathcal{G}(\O(2\Theta)
\boxtimes \O(2\Xi)) \to J[2] \times P[2] \to 1$$ and it is
isomorphic to $\mathcal{G}(\O(2\Theta))\times
\mathcal{G}(\O(2\Xi))/((t,\,t^{-1}),\,t\in \mathbb{G}_m)$.
 Because of  the Theorem 2.2.(1), there is a unique level subgroup $\widetilde{K}_\sigma \subset
\mathcal{G}(\O(2\Theta) \boxtimes \O(2\Xi))$, isomorphic to
$K_\sigma=\ker \sigma$, such that $$\sigma_*(\O(2\Theta) \boxtimes
\O(2\Xi))^{\widetilde{K}_\sigma} \cong \O(\widetilde{\Theta}_z)$$
Denote by $\wt$ the image of $\tau$ via the lifting $K_\sigma
\xrightarrow{\sim} \widetilde{K}_\sigma $.

\begin{prop} (1) There is an isomorphism (well defined up to a
multiplicative scalar)
\begin{eqnarray}\chi :
H^0(J,\,\O(2\Theta))^{\wt} \xrightarrow{\sim}
H^0(P,\,\O(2\Xi))^*\end{eqnarray} (2) The application
$\delta_{\tau,\,z} : P \to |2\Theta|$ that maps a point $q$ in $P$
to the well-defined divisor
$(\pi^*)^{-1}(T^*_{q}\widetilde{\Theta}_z)$ is a morphism that
factors as the composite
$$ P \xrightarrow{\varphi_{2\Xi}}\P
H^0(P,\,\O(2\Xi))^* \xrightarrow{\sim} \P H^0(J,\,\O(2\Theta))^{\wt}
\subset |2\Theta|$$ where the (well-defined) isomorphism is deduced
from $\chi$.
\end{prop}
\dem The results of that proposition can be found in [Mu3], Sections
4 and 5, and we sketch Mumford's arguments again.\\
(1) Consider the homomorphism $\pi^* : J \to \widetilde{J}$. Its
kernel is $<\tau>$ and its image identifies with $J/<\tau>$. We let
$i_1$ be the (separable) isogeny $J \to J/<\tau>$ and $\sigma :
J\times P \to \widetilde{J}$ factors as
$$J \times P \xrightarrow{i_1\times {\rm Id}} J/<\tau> \times P \xrightarrow{i_2} \widetilde{J}$$
One has (Theorem 2.2)
$$i_2^*(\O(\widetilde{\Theta}_z)) \cong ((i_1\times {\rm
Id})_*(\O(2\Theta) \boxtimes \O(2\Xi)))^{\wt} \cong ({i_1}_*
\O(2\Theta))^{\wt} \boxtimes \O(2\Xi)$$ and, letting $Z(\wt)$ be the
centralizer of $\wt$ in $\mathcal{G}(\O(2\Theta))$, one has the
isomorphism of Heisenberg groups (of weight 1 in the sense that it
induces the identity on the centers $\mathbb{G}_m$)
 $$\mathcal{G}(({i_1}_*
\O(2\Theta))^{\wt}) \cong Z(\wt)/<\wt>$$ Furthermore, the space
$H^0(J,\,\O(2\Theta))^{\wt}$ is the unique (up to isomorphism)
irreducible representation of weight 1 of the Heisenberg group
$Z(\wt)/<\wt>$. \\
The symplectic isomorphism  $\varphi : \tau\,^\perp /<\tau>
\xrightarrow{\sim } P[2]$ constructed in the Lemma 3.2 allows us to
construct an isomorphism of Heisenberg groups (of weight $-1$ in the
sense hat it induces $\lambda \mapsto \lambda^{-1}$ on the centers)
$$\widetilde{\varphi} : Z(\wt)/<\wt> \to \mathcal{G}(\O(2\Xi))$$
Indeed, since $e_{2,\,X}(\alpha,\,\beta)=e_{2,\,P}(\varphi
(\bar{\alpha}),\,\varphi (\bar{\beta}))^{-1}$ for any
$\alpha,\,\beta$ in $\tau\,^\perp$ (see (c) in the proof of the
Lemma 3.2), one can find an homomorphism $\widetilde{\varphi}$ such
that the following diagram commutes
$$\begin{array}{rcccccl}
1 \to & \mathbb{G}_m & \to & Z(\wt)/<\wt> & \to &
\tau\,^\perp/<\tau> & \to 1\\
 & \hspace{0.3cm} \wr \downarrow \lambda^{-1} & & \wr \downarrow
\widetilde{\varphi} & & \wr \downarrow \varphi & \\
1 \to & \mathbb{G}_m & \to & \mathcal{G}(\O(2\Xi)) & \to & P[2] &
\to 1
\end{array}$$
and since  $\varphi$ is an isomorphism, $\widetilde{\varphi}$ ditto.
Let $H^0(P,\,\O(2\Xi))^*$ be the dual vector space of
$H^0(P,\,\O(2\Xi))$. It is the unique irreducible representation of
weight $-1$ of $\mathcal{G}(\O(2\Xi))$ and the latter isomorphism of
Heisenberg groups, together with the Theorem 2.1, assures that we
have an isomorphism (unique up to scalar as it can be shown using
Schur's lemma) \begin{eqnarray*}\chi : H^0(J,\,\O(2\Theta))^{\wt}
\xrightarrow{\sim} H^0(P,\,\O(2\Xi))^*\end{eqnarray*}

 (2) Let $s_0$ be a non zero section of
$H^0(\widetilde{J},\,\O(\widetilde{\Theta}_z))$. Its pull-back
$\sigma^*(s_0)$ is the unique (up to scalar)
$\widetilde{K}_\sigma$-invariant element of the space of global
sections $$H^0(J\times P,\,\O(2\Theta) \boxtimes \O(2\Xi)) \cong
H^0(J,\,\O(2\Theta) ) \tens H^0(P,\,\O(2\Xi))$$ Let $q$ be a point
of $P$. The zero locus of $\sigma^*(s_0)_{|J \times \{q\}}$
identifies set-theoretically with the inverse image
$(\pi^*)^{-1}({\rm Supp}\,(T^*_{q}\widetilde{\Theta}_z))$. Suppose
that is proper subset of $J$. Then,
$$(\pi^*)^{-1}(T^*_{q}\widetilde{\Theta}_z)$$ is a well-defined
divisor on $J$, that is algebraically equivalent to $2\Theta$ (see
the Remark 3.1). Because one has taken $q$ in $P$, it is in fact
linearly equivalent to $2\Theta$. Indeed, the line bundle associated
to the algebraically trivial divisor
$(\pi^*)^{-1}(T^*_{q}\widetilde{\Theta}_z)-(\pi^*)^{-1}(\widetilde{\Theta}_z)$
corresponds to the point
$\widehat{\pi^*}(\lambda_{\widetilde{\Theta}}(q))$ in $\hat{J}$.
Because of the diagrams (3.2), $P$ identifies
 with $\lambda_{\widetilde{\Theta}}^{-1}(\ker (\widehat{\pi^*})^0)$,
 and we find that
$$(\pi^*)^{-1}(T_q^*\widetilde{\Theta}_z)\sim_{\rm lin} (\pi^*)^{-1}(\widetilde{\Theta}_z)\sim_{\rm lin}2\Theta$$
Thus $\sigma^*(s_0)$ determines a rational map
$$\delta_{\tau,\,z} : P \dashrightarrow |2\Theta|$$ that maps $q$ to
the divisor $\delta_{\tau,\,z}(q)=
(\pi^*)^{-1}(T_q^*\widetilde{\Theta}_z)$ when defined.\\
Choose a basis $\{\Lambda_0,\,\Lambda_1\}$ of
$H^0(J,\,\O(2\Theta))^{\wt}$. The isomorphism  $\chi$ allows us to
construct a basis $\{\Gamma_0,\,\Gamma_1\}$ of $H^0(P,\,\O(2\Xi))$,
uniquely defined up to a multiplicative scalar, and satisfying the
conditions
$$\chi(\Lambda_i)(\Gamma_j)=\delta_{ij} \text{ for any }i,\,j=0 \text{ or }1$$
By construction, $\sum_i \Lambda_i \boxtimes  \Gamma_i$ is
$\widetilde{K}_\sigma$-invariant, hence equal to $\sigma^*(s_0)$
after suitable normalization. Therefore, the restriction
$\sigma^*(s_0)_{|J \times \{q\}}$ coincides with the section
$$\sum_i \Gamma_i(q).\Lambda_i $$ of $H^0(J,\,\O(2\Theta))$, that is
non zero since the linear system $|2\Xi|$ is base-point free. Thus
the rational map $\delta_{\tau_,\,z}$ is actually a morphism $P \to
|2\Theta|$. Furthermore, the corresponding linear map
$H^0(J,\,\O(2\Theta))^* \to H^0(P,\,\O(2\Xi))$ obviously factors as
$$H^0(J,\,\O(2\Theta))^* \twoheadrightarrow  (H^0(J,\,\O(2\Theta)
)^{\wt})^* \xrightarrow \sim  H^0(P,\,\O(2\Xi))$$ where the
surjection is the dual of the inclusion and the isomorphism is
deduced from $\chi$, hence the proposition. $\square$\\

\subsection{Prym varieties and $M_X$}

For any $q$ in $P$, one can construct a semi-stable rank 2 vector
bundle with trivial determinant over $X$, namely $$\pi_*q \tens
z,\ \text{where } z \in S_\tau=\{z\in J|\,z^2=\tau\}$$ Indeed, the
corresponding sheaf is locally free because $\pi$ is flat, it has
determinant $${\rm Nm}(q) \tens \tau \tens z^2 \cong \O_X$$ (see
the isomorphism (3.1)) and if there were an invertible sub-sheaf
$L \subset \pi_*q \tens z$ with non negative degree, the
projection formula would give a non zero map $\pi^* L \to q \tens
\pi^*(z)$, which is contradictory since
${\rm deg } \pi^* L> {\rm deg} (q \tens \pi^*(z))=0$. \\
Taking the universal line bundle $\mathcal{L}$ over $\widetilde{X}
\times P$ and considering $$((\pi\times {\rm Id})_*\mathcal{L})
\boxtimes z$$ as a family of semi-stable rank 2 vector bundles with
trivial determinant over $X$, parameterized by $P$, we obtain, using
the universal property of $M_X$, a morphism $$d_{\tau,\,z} : P
\to M_X$$ depending on $\tau$ and $z$.\\

\begin{lem}
The following diagram is commutative.
\begin{center}
\unitlength=0.6cm
\begin{picture}(8,4)
\put(0,2){$P$} \put(6,3.7){$M_X$} \put(6.3,0.3){$|2\Theta|$}
\put(0.7,2.4){\vector(4,1){5}} \put(0.7,2){\vector(4,-1){5}}
\put(6.9,3.2){\vector(0,-1){2}} {\footnotesize
\put(2.9,3.5){$d_{\tau,\,z}$} \put(2.9,0.6){$\delta_{\tau,\,z}$}
\put(7.2,1.9){$D$}}
\end{picture}
\end{center}
The intersection $\delta_{\tau,\,z}(P) \cap {\rm Kum}_X$ in
$|2\Theta|$ is the image $\delta_{\tau,\,z}(P[2])$ of $P[2]$.
 \end{lem}
\dem For any $q$ in $P$, the set ${\rm
Supp}\,(\delta_{\tau,\,z}(q))$ is by definition the set
$${\rm Supp}\,((\pi^*)^{-1}(T_{q}^*\widetilde{\Theta}_z))=\{j\in
J|\,h^0(\widetilde{X},\,\pi^*(j\tens z \tens \kappa_0)\tens q)\geq
1\}$$ Because of the adjunction formula, it coincides with the set
$$\{j \in J | \, h^0(X,\,(\pi_*q\tens z) \tens j \tens \kappa_0)\geq
1\}$$ which is precisely the support of the divisor
$D([\pi_*q\tens z])$. As the divisors linearly equivalent to
$2\Theta$ are determined by their support (see [NR], Proposition
6.4), the diagram commutes.

Suppose $\delta_{\tau,\,z}(q)$ is in $\delta_{\tau,\,z}(P) \cap {\rm
Kum}_X$. It means that $\pi_*q\tens z$ is a non-stable bundle. If
$L$ is a degree 0 invertible sub-sheaf of $\pi_*q\tens z$, then the
projection formula gives a non zero morphism of invertible sheaves
$$\pi^*L \to q\tens \pi^*z$$ over $\widetilde{X}$ and $q\tens
\pi^*(z\tens L^{-1})$ must be the trivial sheaf. Thus, with the
notations of Section 3.2, $z\tens L^{-1}$ must be an element of
$K=\tau^\perp$ and $q$ must be an element of $P[2]$, corresponding
to the class $\overline{z\tens L^{-1}}$ of $z\tens L^{-1}$ in
$\tau^\perp/<\tau>$ via the symplectic isomorphism $\varphi$ of the
Lemma 3.2. As $\varphi$ is surjective, we are
done. $\square$\\

\Rq One can ask what happens if one takes $z'\neq z$ in $S_\tau$. In
particular, we shall compare $\delta_{\tau,\,z}$ and
$\delta_{\tau,\,-z}$ in the next section. Let $\alpha$ in $ J[2]$ be
the difference $z'-z$ and consider the maps $d_{\tau,\,z'}$ and
$\delta_{\tau,\,z'}$. The map $d_{\tau,\,z'}$ (resp.
$\delta_{\tau,\,z'}$) coincides with the composite of $d_{\tau,\,z}$
(resp. $\delta_{\tau,\,z}$) with the automorphism of $M_X$ (resp.
$|2\Theta|$) induced by the action of $\alpha$. More precisely, the
unique (up to scalar) isomorphism $\O(\widetilde{\Theta}_{z'})
\xrightarrow{\sim} T^*_{\pi^*(\alpha)}\O(\widetilde{\Theta}_{z})$
 provides an isomorphism $\Psi(\alpha)$
defined as the composite

\begin{center} \unitlength=0.6cm
\begin{picture}(25,4)
\put(0.3,3.5){$ \O(2\Theta)\boxtimes \O(2\Xi)$} \put(0,0){$
T_\alpha^*\O(2\Theta)\boxtimes \O(2\Xi)$}
\put(7.5,0){$T_\alpha^*(\sigma^*\O(\widetilde{\Theta}_{z}))$}
\put(13.5,0){$\sigma^*
(T^*_{\pi^*(\alpha)}\O(\widetilde{\Theta}_{z}))$}
\put(21,0){$\sigma^*\O(\widetilde{\Theta}_{z'})$} \put(20,3.5){$
\O(2\Theta)\boxtimes \O(2\Xi)$} \put(6.1,0){$\xrightarrow \sim $}
\put(12.1,0){$\xrightarrow \sim $} \put(19.4,0){$\xrightarrow \sim
$} \put(2.7,3.1){\vector(0,-1){2}} \put(22.3,1){\vector(0,1){2}}
\put(2.9,2){$\wr$} \put(22.5,1.7){$\wr$}
\dashline{0.2}(5.5,3.7)(19,3.7) \put(19,3.7){\vector(1,0){0.2}}
\small \put(0.7,2){$\widetilde{\alpha} \tens 1$}
\put(11.5,3.9){$\Psi(\alpha)$}
\end{picture}
\end{center}
 where
$\widetilde{\alpha}$ is any lifting of $\alpha$ in
$\mathcal{G}(\O(2\Theta))$ and where we use the splitting of
$\sigma^*(\O(\widetilde{\Theta}_z)$ into $\O(2\Theta) \boxtimes
\O(2\Xi)$. Now, the automorphism of
$\mathcal{G}(\O(2\Theta)\boxtimes \O(2\Theta))$ defined by
$$\gamma \mapsto \Psi(\widetilde{\alpha})\circ \gamma \circ
\Psi(\widetilde{\alpha})^{-1}$$ depends only on the choice of $z$
and $z'$ and it gives an isomorphism $\widetilde{K}_\sigma
\xrightarrow{\sim} \widetilde{K}'_\sigma$, where
$\widetilde{K}'_\sigma$ is the unique lifting of $K_\sigma$ in the
theta group scheme $\mathcal{G}(\O(2\Theta)\boxtimes \O(2\Theta))$
such that $\sigma_*(\mathcal{G}(\O(2\Theta)\boxtimes
\O(2\Theta)))^{\widetilde{K}'_\sigma} \cong
\O(\widetilde{\Theta}_{z'})$. Let $\gamma_J$ be the image of
$\gamma$ through the quotient map $\mathcal{G}(\O(2\Theta)\boxtimes
\O(2\Xi))\to J[2]\times P[2] \to J[2]$. We find that
$$\Psi(\widetilde{\alpha})\circ \gamma \circ
\Psi(\widetilde{\alpha})^{-1}=e_{2,X}(\alpha,\,\gamma_J)\,.\gamma$$
If $\wt$ (resp. $\wt \,' $) is the image of $\tau$ via the
isomorphism $K_\sigma \xrightarrow \sim \widetilde{K}_\sigma$ (resp.
$K_\sigma \xrightarrow \sim \widetilde{K}'_\sigma$), one has
$$\wt\,'= \Psi(\widetilde{\alpha})\circ \wt \circ
\Psi(\widetilde{\alpha})^{-1}=e_{2,X}(\alpha,\,\tau)\,.\wt$$ and one
finds that if $\alpha$ is not in $\tau^\perp$, $\delta_{\tau,\,z}$
and $\delta_{\tau,\,z'}$ don't map $P$ onto the same projective line
in $|2\Theta|$.  If $\alpha$ belongs to $\tau^\perp$, then
$\widetilde{K}_\sigma = \widetilde{K}'_\sigma$ and the automorphism
of $\widetilde{K}_\sigma $  induced by $\gamma \mapsto
\Psi(\widetilde{\alpha})\circ \gamma \circ
\Psi(\widetilde{\alpha})^{-1}$ only depends  on the class of
$\alpha$ in $\tau^\perp/<\tau> \cong P[2]$. In other words, the
morphisms $\delta_{\tau,\,z}$ and $\delta_{\tau,\,z'}$ differ from
the involutional translation on $P$ corresponding to the class of
$\alpha$ in $P[2]$. In particular, one has
$\delta_{\tau,\,-z}=\delta_{\tau,\,z}$.\\

\Rq Recall (Corollary 2.14) that we are able to give the homogeneous
coordinates of the points in $\delta_{\tau,\,z}(P) \cap {\rm Kum}_X$
in $|2\Theta|$ in terms of the coefficients of the Kummer surface
(see the equation (2.4) and the chart (2.5)). As $\delta_{\tau,\,z}$
coincides with the canonical map $\varphi_{2\Xi}$, this set is
precisely the set of ramification points of the canonical map
$\varphi_{2\Xi}$. Therefore, we are able
to characterize the elliptic curves arising that way.\\

\section{The generalized Verschiebung $V : M_{X_1} \dashrightarrow M_X$}

\subsection{Review of Theta groups in characteristic $p$}

The curve $X$ is now supposed to have $p$-rank 2, i.e., to be an
ordinary genus 2 curve. Because of the Remark 2.10, one knows that a
general Kummer surface ${\rm Kum}_X$ in $\P^3$ is associated to
such a curve.\\

 Let $X_1$  be the $p$-twist
of $X$. Denote by $i$ the semi-$k$-linear isomorphism $X_1 \to X$
and by $F$ the relative Frobenius $X \to X_1$, which is radicial and
flat. Notice that there is a canonical bijection from the set of
Weierstrass points of $X$ to the set of Weierstrass points of $X_1$.
Thus, the choice of the effective theta characteristic $\kappa_0$
for $X$ determines an effective theta characteristic $i^*(\kappa_0)$
for $X_1$, still denoted by $\kappa_0$.\\
Denote by $J_1$ the $p$-twist of $J$ and let $F : J \to J_1$ (resp.
$i : J_1 \to J$) the relative Frobenius (resp. the semi-$k$-linear
isomorphism) which is flat. The abelian variety $J_1$ coincides with
the Jacobian of $X_1$ and $\kappa_0$ enables us, as before, to give
a symmetric representative $\Theta_1$ for the principal polarization
on $J_1$.
It is easily seen that $\Theta_1=i^*\Theta$.\\

Let $\hat{G}$ denote the kernel of $F : J \to J_1$. It is a local
group scheme and as $J$ is an ordinary abelian variety, $\hat{G}$ is
the local part of the group scheme $J[p\,]$ of $p$-torsion points of
$J$. Denote by $G$ the reduced part of $J[p\,]$. The relative
Frobenius maps $G$ isomorphically onto the kernel of the isogeny
$$J_1\cong J/\hat{G} \to J/J[p\,] \cong J$$ which is separable and of
degree $p^g=p^2$, and is called the Verschiebung $V$. It maps a
degree 0 line bundle $\zeta_1$ over $X_1$ to the degree 0 line
bundle $F^*\zeta_1$ over $X$. Note that both composites $$J
\xrightarrow F J_1 \xrightarrow V J\text{ and }J_1\xrightarrow{V} J
\xrightarrow{F} J_1$$ are multiplication by $p$.\\

The finite group $J_1[p\,]$ is self-dual and consequently, we have
the Göpel system
$$J_1[p\,] \cong G \times \hat{G} \cong (\Z/p\,\Z)^g\times (\mu_p)^g$$
Consider the line bundle $\O(p\,\Theta_1)$ over $J_1$. Its
automorphism group $\mathcal{G}(\O(p\,\Theta_1)$) which can be
obtained as a central extension $$ 1 \to \mathbb{G}_m \to
\mathcal{G}(\O(p\,\Theta_1)) \to J_1[p\,] \to 0$$ and
$\O(p\,\Theta_1)$ is no longer of separable type. Nevertheless,
[Sek] proves that the main results about theta groups (recalled in
Theorems 2.2 and 2.3) extend to line bundles of non-separable type.
We gather some useful results given in [LP1] in the

\begin{lem}
(1) There are the three isomorphisms $\O(p\,\Theta_1)\cong
V^*\O(\Theta)$,
$\O(p\,\Theta) \cong F^*\O(\Theta_1)$ and $\O(\Theta_1) \cong i^*\O(\Theta)$. \\
(2) The restrictions of $\mathcal{G}(\O(p\,\Theta))$ to both $G$ and
$\hat{G}$ are canonically split. Therefore, the decomposition
$J[p\,] \cong G \times \hat{G}$ is symplectic (with respect to $e_p$).\\
(3) There exists a basis $\{X_g\}_{g\in\, G}$ of $H^0(J,
\O(p\,\Theta))$, unique up to a multiplicative scalar, which
satisfies the following relations
$$ a.X_g=X_{a+g} \ \ \ \alpha.X_g=e_p(\alpha,\,g) X_g \ \ \
\forall a,\,g \in G,\, \alpha \in \hat{G}$$ (4) For any $g\in G$,
there is a unique $Y_g$ in $H^0(J_1,\O(p\,\Theta_1))$ such that
$X^p_g=F^*Y_g$. The family $\{Y_g\}_{g\in\,G}$ is a basis of
$H^0(J_1,\O(p\,\Theta_1))$ that corresponds to the basis
$\{X_g\}_{g\in\,G}$ via $i^* : H^0(J,\O(p\,\Theta)) \to
H^0(J_1,\O(p\,\Theta_1))$.
\end{lem}
\underline{Sketch of Proof} (Complete proofs can be found in [LP1])
: One needs to define a splitting $$G \hookrightarrow
\mathcal{G}(\O(p\,\Theta_1))$$ for the central extension above.
Because $G$ is reduced, it is enough to find it at the level of
$k$-point and because $k^*$ is divisible, this can be worked out.
Furthermore, because the skew-symmetric form
$$e_p : J_1[p\,] \times J_1[p\,] \to \mathbb{G}_m$$ (associated to
the commutator in $\mathcal{G}(\O(p\,\Theta_1))$) takes its value in
$\mu_p$, this splitting is unique (another one would differ from the
first one by a morphism $G \to \mu_p$ and $\mu_p(k)=\{1\}$).
Therefore, the analog of Theorem 2.3.(1) in the non separable case
 assures the existence and the uniqueness of a line
 bundle $M$ over $J$ such that $\O(p\,\Theta_1)\cong V^*M$.\\
 One can
 show that $M$ defines a principal polarization and that the
 isomorphism $$V^*F^*\O(\Theta_1) \cong \O(p^2\,\Theta_1)\xrightarrow \sim  V^*M^{\tens p}$$
 obtained in taking $p$-powers commutes with the action of $G$,
 hence descends onto an isomorphism
 $$F^*\O(\Theta_1) \xrightarrow \sim  M^{\tens p}$$
Furthermore, if $D$ is the unique effective divisor on $J$ such that
$L=\O(D)$, on has, set-theoretically, $$V^{-1}({\rm
Supp}(D))=\bigcup_{g \in G} T_g^*({\rm Supp}(\Theta_1))$$ Because
$\Theta_1$ is symmetric, $V^{-1}({\rm Supp}(D))$ is symmetric and
finally $M$ is a symmetric representative for the principal
polarization on $J$. Thus, the difference between $M$ and
$\O(\Theta)$ lies in $J[p\,] \cap J[2]= \{0\}$ and $M\cong
\O(\Theta)$. The third isomorphism announced by (1) is tautological
(see our
definition of $\Theta_1$).\\
(2) is the Lemma 2.3 of [LP1].\\
(3) can be obtained in a very similar way as in the construction of
the basis $\{X_\bullet\}$ in Section 2.2, taking of the fact that the invertible sheaf $\O(p\,\Theta)$ is of non separable type.
\\
(4) Each vector $X_g^p$ is invariant under the action of $\hat{G}$
hence of the form $F^*Y_g$ with $Y_g$ in $H^0(J_1,\O(p\,\Theta_1))$.
Since $k$ is divisible, the family $\{Y_g\}$ is free, hence a basis.
$\square$\\

\Rq These results remain true, {\it mutatis mutandis},
for any ordinary principally polarized abelian variety (see [LP1]).\\

\subsection{Extending Verschiebung to $|2\Theta_1|$}

Let $E_1$ be a semi-stable bundle with trivial determinant over
$X_1$. Then, $F^*E_1$ is a rank 2 vector bundle with trivial
determinant over $X$, which may not be semi-stable since the
pull-back by Frobenius  destabilizes some vector bundles ([R]).
Nevertheless, $F$ induces a rational map $V : M_{X_1}
\dashrightarrow M_X$. If $E_1$ is the non stable bundle $j \oplus
j^{-1}$, its pull-back $F^*j\oplus (F^*j)^{-1}$ is semi-stable but
non stable and we find that the following diagram is commutative
\begin{eqnarray} \end{eqnarray}
\begin{center} \unitlength=0.6cm
\begin{picture}(10,2)
\put(0.3,0){$M_{X_1}$} \put(8.3,0){$M_{X}$} \put(0.8,3){$J_1$}
\put(8.9,3){$J$} \put(2,3.2){\vector(1,0){6}}
\put(1.2,2.7){\vector(0,-1){1.8}} \put(9,2.7){\vector(0,-1){1.8}}
\put(7.6,0.2){\vector(1,0){0.1}} \dashline{0.2}(2.6,0.2)(7.6,0.2)

\small \put(5,3.4){$V$} \put(5,0.6){$V$} \put(0.4,1.5){$b_1$}
\put(9.3,1.5){$b$}
\end{picture}
\end{center}
It is a diagram of $J[2]$-equivariant morphisms in the following
sense : On the one hand, because $p$ is odd, $[p\,]$ induces
identity on $J[2]$. Therefore,  $F : J[2] \to J_1[2]$ and $V :
J_1[2] \to J[2]$ are isomorphisms, inverse one to each other. Thus,
one can define an action of $J[2]$ on both $J_1$ and $M_{X_1}$,
compatible with the maps involved in the diagram.\\
 In particular,
the indeterminacy locus of $V$ does not meet the Kummer surface.
Thus, it is a finite set $\mathcal{I}$ and we let $U$ be the Zariski
open subset
$M_{X_1} \setminus \mathcal{I}$.\\

Let $\O(\Delta_1)$ be the determinant line bundle over $M_{X_1}$. It
has been shown ([LP2]) that
\begin{prop} There is an isomorphism $V_2^*(\O(\Delta))
\cong (\O(\Delta_1)^{\tens p})_{|U}$.\end{prop} Because
$\mathcal{I}$ is a finite set, they obtain that the rational map $V$
is given by degree $p$ polynomials.

Using the Proposition 2.6, one obtains a rational map $\widetilde{V}
: |2\Theta_1|^* \dashrightarrow |2\Theta|^*$ and the vertical arrows
in (4.1) become the canonical maps $\varphi_{2\Theta_1}$ and
$\varphi_{2\Theta}$. Thus, the pull-back by the Verschiebung $V^* :
H^0(J,\,\O(2\Theta)) \to H^0(J_1,\,\O(2p\,\Theta_1))$ factors as the
composite
\begin{eqnarray}H^0(J,\,\O(2\Theta)) \xrightarrow
{\widetilde{V}^*} {\rm Sym}^p\, H^0(J_1,\,\O(2\Theta_1)) \to
H^0(J_1,\,\O(2p\,\Theta_1))\end{eqnarray} where the last arrow is
the canonical evaluation map.\\

Recall (Section 2.2) that we have chosen a theta structure
$\mathcal{H} \xrightarrow \sim \mathcal{G}(\O(2\Theta))$ and that
$W:=H^0(J,\,\O(2\Theta))$ is the unique (up to isomorphism)
irreducible representation of $\mathcal{H}$ (of weight 1). We let
$W_1$ denote the vector space
$$W_1:=H^0(J_1,\,\O(2\Theta_1))$$ It is,
analogously, the unique irreducible representation of the Theta
group
$\mathcal{G}(\O(2\Theta_1))$.\\

\begin{lem}
(1) One can endow ${\rm Sym}^p\,W_1 $ and
$H^0(J_1,\,\O(2p\,\Theta_1))$ with an action (of weight $p$) of
$\mathcal{G}(\O(2\Theta_1))$.\\
(2) The evaluation map ${\rm Sym}^p\, W_1 \to
H^0(J_1,\,\O(2p\,\Theta_1))$ is
$\mathcal{G}(\O(2\Theta_1))$-equivariant for these actions.
\end{lem}
\dem $(1)$ The homomorphism $$\varepsilon_p :
{\mathcal{G}(\O(2\,\Theta_1))} \to  \mathcal{G}(\O(2p\,\Theta_1))$$
that maps an isomorphism $\gamma : \O(2\,\Theta_1) \xrightarrow \sim
T_x^*\O(2\,\Theta_1)$ to the isomorphism
$$\gamma^{\tens p} : \O(2p\,\Theta_1) \xrightarrow \sim
T_x^*\O(2p\,\Theta_1)$$ fits into the commutative diagram
\begin{center} \unitlength=0.6cm
\begin{picture}(16,3)
\put(0.5,0){\put(0,2.5){$0 \to \mathbb{G}_m $} \put(0,0){$0\to
\mathbb{G}_m $} \put(15,0){\put(0,2.5){$0$} \put(0,0){$0$}}
\put(5.3,2.5){$\mathcal{G}(\O(2\Theta_1))$}
\put(5,0){$\mathcal{G}(\O(2p\,\Theta_1))$} \put(11.6,2.5){$J_1[2]$}
\put(11.3,0){$J_1[2p\,]$}

\put(2.8,2.7){\vector(1,0){1.9}} \put(2.8,0.2){\vector(1,0){1.8}}
\put(8.8,2.7){\vector(1,0){2.4}} \put(8.9,0.2){\vector(1,0){2.2}}
\put(13.2,2.7){\vector(1,0){1.3}} \put(13.4,0.2){\vector(1,0){1.1}}
\put(2,2.2){\vector(0,-1){1.5}} \put(7,2.2){\vector(0,-1){1.5}}
\put(12.3,2.2){\vector(0,-1){1.5}}}
 \small \put(0,1.3){$p$-power}
\put(7.8,1.3){$\varepsilon_p$} \put(13.1,1.3){inclusion}
\end{picture}
\end{center}
It gives the action of $\mathcal{G}(\O(2\Theta_1))$ onto
$H^0(J_1,\,\O(2p\,\Theta_1))$.
The other case is straightforward.\\
$(2)$ The evaluation map $W_1 \tens \O_{X_1} \to \O(2\Theta_1)$ is
of course $\mathcal{G}(\O(2\Theta_1))$-equivariant and taking its
$p$-symmetric power, one obtains, at the level of global sections,
the canonical map ${\rm Sym}^p\, W_1 \to
H^0(J_1,\,\O(2p\,\Theta_1))$, which is still
$\mathcal{G}(\O(2\Theta_1))$-equivariant for the induced actions on
both spaces. These are the ones of $(1)$, hence the lemma.
$\square$\\

\begin{lem}
(1) There is an homomorphism of theta groups $\mathcal{H} \to
\mathcal{G}(\O(2\Theta_1))$ (of weight $p$). Therefore,
$\mathcal{H}$ has an action of weight $p^2$ on ${\rm Sym}^p\,
W_1$.\\
(2) The map $\widetilde{V}^*  : W \to {\rm Sym}^p\, W_1$ is
injective and $\mathcal{H}$-equivariant, up to a multiplicative
scalar.
\end{lem}
\dem $(1)$ On the one hand, the isomorphism of sheaves $V^*
\O(2\Theta) \cong \O(2p\,\Theta_1)$ induces (by pull-back) a
homomorphism $V^* : \mathcal{H}\to {\mathcal{G}(\O(2p\,\Theta_1)}$
(of weight 1). On the other hand, for any $\gamma : \O(2p\,\Theta_1)
\xrightarrow \sim T_x^* \O(2p\,\Theta_1)$, there is a unique
isomorphism $\rho : \O(2\Theta_1) \xrightarrow{\sim}
T_{p\,x}^*\O(2\Theta_1)$ such that the following diagram commutes
\begin{center}
\unitlength=0.6cm
\begin{picture}(12,4)
\put(0.8,2.5){$\O(2p^2\,\Theta_1)$}
\put(0.5,0){$[p\,]^*\O(2\Theta_1)$} \put(2.2,0.7){\vector(0,1){1.5}}
\put(2.4,1.3){$\wr$}
\put(2.4,0){\put(6.4,2.5){$T_x^*\O(2p^2\,\Theta_1)$}
\put(6,0){$[p\,]^*T_{p\,x}^*\O(2\Theta_1)$}
\put(1.8,2.7){\vector(1,0){4}} \put(2.2,0.2){\vector(1,0){3.4}}
\put(1,0){\put(7.2,0.7){\vector(0,1){1.5}}\put(7.5,1.3){$\wr$}}
\small \put(3.4,3.1){$\gamma^{\tens p}$}
\put(3.4,0.5){$[p\,]^*\rho$} }
\end{picture}
\end{center}
where the vertical isomorphisms come from the fact that
$\O(2\Theta_1)$ is a symmetric line bundle. We let  $\eta_p :
{\mathcal{G}(\O(2p\,\Theta_1))}\to {\mathcal{G}(\O(2\,\Theta_1))}$
be the homomorphism that maps $\gamma $ to $\rho$. It fits into the
commutative diagram of central extensions
\begin{center} \unitlength=0.6cm
\begin{picture}(16,3)
\put(0.5,0){\put(0,2.5){$0 \to \mathbb{G}_m $} \put(0,0){$0\to
\mathbb{G}_m $} \put(15,0){\put(0,2.5){$0$} \put(0,0){$0$}}

\put(5.3,0){$\mathcal{G}(\O(2\Theta_1))$}
\put(5,2.5){$\mathcal{G}(\O(2p\,\Theta_1))$}

\put(11.6,0){$J_1[2]$} \put(11.3,2.5){$J_1[2p\,]$}

\put(2.8,0.2){\vector(1,0){1.9}} \put(2.8,2.7){\vector(1,0){1.8}}
\put(8.8,0.2){\vector(1,0){2.4}} \put(8.9,2.7){\vector(1,0){2.2}}
\put(13.2,0.2){\vector(1,0){1.3}} \put(13.4,2.7){\vector(1,0){1.1}}
\put(2,2.2){\vector(0,-1){1.5}} \put(7,2.2){\vector(0,-1){1.5}}
\put(12.3,2.2){\vector(0,-1){1.5}}}
 \small \put(0,1.3){$p$-power}
\put(7.8,1.3){$\eta_p$} \put(13.1,1.3){$[p\,]$}
\end{picture}
\end{center}
Now, we consider the composite homomorphism
$$\mathcal{H}\xrightarrow{V^*}{\mathcal{G}(\O(2p\,\Theta_1))}
\xrightarrow{\eta_p} {\mathcal{G}(\O(2\Theta_1))}$$ It has weight
$p$. Using the latter and the natural action (of weight $p$) of
$\mathcal{G}(\O(2\Theta_1))$ on ${\rm Sym}^p\,W_1$, one obtains an
action of weight $p^2$ of $\mathcal{H}$ on ${\rm
Sym}^p\,W_1$.\\
$(2)$ Let
$$\mathcal{G}(\O(2p\,\Theta_1))_2$$  be  the maximal subgroup of
$\mathcal{G}(\O(2p\,\Theta_1))$ lying above $J_1[2]$, viewed as a
sub-group of $J_1[2p\,]$. It is obviously the image of the
homomorphism  $\varepsilon_p$ and since $[p\,]$ acts trivially on
$J_1[2]$, $\eta_p$ restricts to an homomorphism
$$\mathcal{G}(\O(2p\,\Theta_1))_2 \to \mathcal{G}(\O(2\Theta_1))$$
with kernel $\mu_p$. The composite homomorphism $\varepsilon_p \circ
\eta_p$ fits into the commutative diagram of central extensions
\begin{center} \unitlength=0.6cm
\begin{picture}(16,3)
\put(0.5,0){\put(0,2.5){$0 \to \mathbb{G}_m $} \put(0,0){$0\to
\mathbb{G}_m $} \put(15,0){\put(0,2.5){$0$} \put(0,0){$0$}}

\put(5,0){$\mathcal{G}(\O(2p\,\Theta_1))_2$}
\put(5,2.5){$\mathcal{G}(\O(2p\,\Theta_1))_2$}

\put(11.6,0){$J_1[2]$} \put(11.6,2.5){$J_1[2]$}

\put(2.8,0.2){\vector(1,0){1.9}} \put(2.8,2.7){\vector(1,0){1.8}}
\put(9.2,0.2){\vector(1,0){2.2}} \put(9.2,2.7){\vector(1,0){2.2}}
\put(13.2,0.2){\vector(1,0){1.3}} \put(13.4,2.7){\vector(1,0){1.1}}
\put(2,2.2){\vector(0,-1){1.5}} \put(7,2.2){\vector(0,-1){1.5}}
\put(12.3,2.2){\vector(0,-1){1.5}}}
 \small \put(0,1.3){$p^2$-power}
\put(7.8,1.3){$\varepsilon_p \circ \eta_p$} \put(13.1,1.3){Id}
\end{picture}
\end{center}
Now, since $V^*(\mathcal{H}) \subseteq
\mathcal{G}(\O(2p\,\Theta_1))_2$, the map $V^* : W\to
H^0(J_1,\,\O(2p\,\Theta_1))$ is $\mathcal{H}$-equivariant (up to a
multiplicative scalar) when one endows $H^0(J_1,\,\O(2p\,\Theta_1))$
with the action $\mathcal{H}$ induced by the composite
$\varepsilon_p \circ \eta_p \circ V^*$. As this map is non-zero, it
is injective and the map $\widetilde{V}^*$ is injective as well. The
evaluation map must induce an isomorphism between the image of
$\widetilde{V}^*$ and the image of $\widetilde{V}^*$, that is
$\mathcal{H}$-equivariant (Lemma 4.5) and $\widetilde{V}^*$ is
$\mathcal{H}$-equivariant (up to a multiplicative scalar) as well.
$\square$\\

\Rq The reason why $\widetilde{V}^*$ is not
$\mathcal{H}$-equivariant is that the action of $\mathcal{H}$ on the
spaces $W$ and ${\rm Sym}^p\,W_1$ do not have the same weight, in
contradiction with the $k$-linearity of $\widetilde{V}^*$. This
obstruction vanishes when one considers the induced action of the
subgroup $\hat{H}\subseteq \mathcal{H}$ (resp. $H\subseteq
\mathcal{H}$) on both spaces.\\

\Rq One notices that taking $\rho = i^*\gamma : \O(2\Theta_1) \to
T_{\alpha_1}^*\O(2\Theta_1)$ makes the following diagram commutative
\begin{center}
\unitlength=0.6cm
\begin{picture}(10,4)
\put(0.4,2.5){$\O(2p^2\,\Theta_1)$}
\put(0.3,0){$[p\,]^*\O(2\Theta_1)$} \put(2,0.7){\vector(0,1){1.5}}
\put(2.2,1.3){$\wr$}
 \put(2,0){\put(6.4,2.5){$T_{\alpha_1}^*\O(2p^2\Theta_1)$}
\put(6.2,0){$[p\,]^*T_{\alpha_1}^*\O(2\Theta_1)$}
\put(1.8,2.7){\vector(1,0){4}}\put(2.2,0.2){\vector(1,0){3.4}}
\put(8.2,0.7){\vector(0,1){1.5}} \small
\put(3,3.1){$(V^*\,\gamma)^{\tens p}$} \put(3.4,0.5){$[p\,]^*\rho$}
\put(8.5,1.3){$\wr$}}
\end{picture}
\end{center}
Namely, one has, using the fact that $F_{abs}^*$ is the $p$-power,
that
$$[p\,]^*i^*\gamma=V^*(F^*i^*)\gamma=V^*(\gamma^{\tens
p})=(V^*\,\gamma)^{\tens p}$$ Therefore, the homomorphism $\eta_p
\circ V^*$  (of weight $p$) coincides with the homomorphism $i^* :
\mathcal{H} \to \mathcal{G}(\O(2\Theta_1))$ induced by the pull-back
by the quasi-isomorphism $i^*$.\\

\Rq Using the map $\eta_p\circ V^* $, one finds that, for any two
elements $\bar{\alpha}$ and $\bar{\beta}$ in $J[2]$,
$e_2(F(\bar{\alpha}),\,F(\bar{\beta}))=e_2(\bar{\alpha},\,\bar{\beta})^p$.
Because $J[2]$ is reduced and because $e_2$ takes its values in
$\mu_2$,we find that $F$ (hence $V$) is a symplectic isomorphism.
This implies that the Göpel system $J[2] \cong H \times \hat{H}$
determines  a Göpel system $J_1[2] \cong H \times \hat{H}$, that a
theta structure $\mathcal{H} \xrightarrow \sim
\mathcal{G}(2\Theta)$ determines a theta structure $\mathcal{H}_1
\xrightarrow \sim \mathcal{G}(2\Theta_1)$ (where $\mathcal{H}_1: =
\mathcal{H}\tens _{F_{\rm abs},\,k} k$), and that the basis
$\{X_\alpha\}_{\alpha\in \,H}$ determines a basis
$\{Y_{\alpha_1}\}_{\alpha_1 \in \,H_1}$,
compatible in the sense that $Y_{\alpha_1}=i^*X_{V(\alpha_1)}$.\\

\begin{prop} Let $X$ be a smooth and proper ordinary  curve of genus 2 over an algebraically
closed field of characteristic $p\geq 3$. Then, the generalized
Verschiebung $V : M_{X_1} \dashrightarrow M_X$ is completely
determined by an irreducible sub-representation of the
$\mathcal{G}(\O(2\Theta))$ in ${\rm
Sym}^p\,H^0(J_1,\,\O(2\Theta_1))^*$ (isomorphic to
$H^0(J,\,\O(2\Theta))^*$ as a vector space).
\end{prop} \dem As
previously, we identify the space $W$ (resp.  $W_1$) with its dual
$W^*$ (resp. $W_1$) by means of the isomorphism of the Lemma 2.4,
and we let $\{x_\bullet\}$ (resp. $\{y_\bullet\}$) be the basis of
$W^*$ (resp. $W_1^*$) dual to $\{X_\bullet\}$ (resp.
$\{Y_\bullet\}$) constructed in the Section 2.2 (resp. in the Remark
4.8). Denote by $V_i$ ($i=00,\,01,\,10$ and $11$) the degree $p$
homogeneous polynomials (see the Proposition 4.3 above, due to
[LP2]) such that the
  rational map $\widetilde{V} : |2\Theta_1| \dashrightarrow |2\Theta|$ corresponding
  to $V$ is given by \begin{eqnarray}\begin{array}{cccl}
  \widetilde{V} : & |2\Theta_1| & \dashrightarrow & |2\Theta|\\
& (y_i) & \mapsto & (V_i(\underline{y})) \end{array}\end{eqnarray}
Using the Lemma 4.5 and the Remark 4.6, we find that $V_{00}$ is
invariant under the action of the subgroup $\hat{H}$ of
$\mathcal{H}$. It coincides therefore, up to a multiplicative
scalar, with $\widetilde{V}(x_{00})$ and upon normalizing suitably
$V_{00}$, one can suppose that
$$\widetilde{V}^*(x_{00})=V_{00}$$
By construction again, $\widetilde{V}^*$ is $H$-equivariant thus one
obtains $V_i$ (for $i=01,\,10,\,11$) as the transform of $V_{00}$
under the action of the unique element of $H$ that maps $x_{00}$ to
$x_i$. $\square$\\

\subsection{Prym varieties and Frobenius}

Let $\tau$ be a non zero element of $J[2]$ and let $\pi
:\widetilde{X} \to X$ be the corresponding \'etale double cover.
The base change induced by the Frobenius morphism on the base
field gives an \'etale double cover $\pi_1 :\widetilde{X}_1 \to
X_1$ corresponding to a non zero $\tau_1$ of $J[2]$ which is the
image of $\tau$ under the isomorphism $J[2 ] \to J_1[2]$ defined
earlier. The following lemma is well-known (see [SGA1]) :
\begin{lem} If $F$ is the relative Frobenius, the following
diagram
\begin{center}
\unitlength=0.6cm
\begin{picture}(8,4)
\put(1,2.5){$\widetilde{X}$} \put(6,2.5){$\widetilde{X}_1$}
\put(1,0){$X$}\put(6,0){$X_1$} \put(1.6,2.7){\vector(1,0){4}}
\put(1.6,0.2){\vector(1,0){4}} \put(1.2,2.3){\vector(0,-1){1.5}}
\put(6.2,2.3){\vector(0,-1){1.5}} \small \put(3.5,0.4){$F$}
\put(3.5,2.9){$F$} \put(0.5,1.4){$\pi$} \put(6.6,1.4){$\pi_1$}
\end{picture}\end{center}
is cartesian.
\end{lem}
\dem The diagram is certainly commutative, because of Frobenius
functoriality, thus there is a unique $ h : \widetilde{X} \to
\widetilde{X}_1 \times _{X_1} X$ that makes the following diagram
commutative
\begin{center}
\unitlength=0.6cm
\begin{picture}(12,6)
\put(0,4.5){$\widetilde{X}$} \put(0.6,4.6){\vector(4,-1){7}}
\put(0.4,4.2){\vector(2,-3){2.5}} \put(0.5,4.4){\vector(3,-2){1.6}}
\put(1.8,2.5){$\widetilde{X}_1 \times _{X_1} X$} \put(2,0){

\put(6,2.5){$\widetilde{X}_1$} \put(1,0){$X$}\put(6,0){$X_1$}
\put(2.8,2.7){\vector(1,0){2.8}} \put(1.6,0.2){\vector(1,0){4}}
\put(1.2,2.3){\vector(0,-1){1.5}} \put(6.2,2.3){\vector(0,-1){1.5}}
\small \put(3.5,0.4){$F$} \put(4,2.2){$pr_1$} \put(1.4,1.4){$pr_2$}
\put(6.6,1.4){$\pi_1$} \put(3.5,1.2){$\square$} }

\small \put(1.6,3.8){$h$} \put(1,2){$\pi$} \put(4.6,3.8){$F$}
\end{picture}\end{center}
As $\pi$ and $\pi_1$ are \'etale and proper, we find that $pr_2$,
hence $h$, are \'etale and proper as well, hence a finite \'etale
covering. Now, $F$ being radicial, $h$ is radicial as well and
therefore, an
isomorphism. $\square$\\

As above, write $\widetilde{J}_1$ for the Jacobian variety of
$\widetilde{X}_1$. It coincides with the $p$-twist of
$\widetilde{J}$ and we can define the relative Frobenius and the
Verschiebung. Denote again by $\widetilde{\Theta}_1$ the Theta
divisor obtained as the pull-back of the canonical Theta divisor on
$\widetilde{J}_1\,^2$ by means of the Theta characteristic
$\pi_1^*(\kappa_0)$ (where $\kappa_0$ has to be understood here as
the effective Theta characteristic of $X_1$ we have constructed).

\begin{lem}
The morphisms $\pi^* : J \to \widetilde{J} $ and $\emph{Nm} :
\widetilde{J} \to J$ commute with $V$. In other words, the two
following diagrams are commutative

\begin{center}
\unitlength=0.6cm
\begin{picture}(17,3.5)
\put(0.8,2.5){$\widetilde{J}_1$} \put(6,2.5){$\widetilde{J}$}
\put(0.9,0){$J_1$}\put(6,0){$J$} \put(1.6,2.7){\vector(1,0){4}}
\put(1.6,0.2){\vector(1,0){4}} \put(1.2,0.7){\vector(0,1){1.5}}
\put(6.2,0.7){\vector(0,1){1.5}}

{\small \put(3.7,2.9){$V$} \put(0.4,1.3){$\pi_1^*$}
\put(3.7,0.4){$V$} \put(6.5,1.3){$\pi^*$}}

\put(9,0){

\put(0.8,2.5){$\widetilde{J}_1$} \put(6,2.5){$\widetilde{J}$}
\put(0.9,0){$J_1$}\put(6,0){$J$} \put(1.6,2.7){\vector(1,0){4}}
\put(1.6,0.2){\vector(1,0){4}} \put(1.2,2.3){\vector(0,-1){1.5}}
\put(6.2,2.3){\vector(0,-1){1.5}}

\small \put(3.7,2.9){$V$} \put(0,1.3){$\emph{Nm}$}
\put(3.7,0.4){$V$} \put(6.5,1.3){$\emph{Nm}$} }
\end{picture} \end{center}
\end{lem}
\dem The commutation of the left-hand diagram is a straightforward
consequence of the commutation of the diagram in the Lemma 4.10. For
the right-hand one, one takes an element $j\in \widetilde{J}_1$ and
sees, using the Lemma 4.10 again, that $F^*(\pi_* j) \cong \pi_*
(F^*j)$. Thus, using (3.1), one can write that
\begin{eqnarray*}
\text{Nm}(V(j))=\text{det}(F^*j)\tens \tau & \cong &
F^*(\text{det}(\pi_*\,j)) \tens \tau\\
& \cong & F^*(\text{det}(\pi_*\,j) \tens \tau_1)= V(\text{Nm}(j)) \
\ \ \ \square \end{eqnarray*}

\begin{prop}
The following diagram
\begin{center}
\unitlength=0.6cm
\begin{picture}(10,3.5)
\put(0.3,0){\put(0.8,0){\put(0.2,2.5){$J\times P$}
 \put(0,0){$J_1\times
P_1$}\put(1.1,0.7){\vector(0,1){1.5}}

\put(0.8,0){\put(6,0){$\widetilde{J}_1$}
\put(6,2.5){$\widetilde{J}$} \put(1.6,2.7){\vector(1,0){4}}
\put(1.6,0.2){\vector(1,0){4}} \put(6.2,0.7){\vector(0,1){1.5}}}}

\small \put(4.7,2.9){$\sigma$} \put(4.7,0.4){$\sigma_1$}
\put(8,1.1){$V$}} \small \put(0.3,1.1){$V \times V$}
\end{picture} \end{center} is commutative.\\
Furthermore, $\sigma$ induces an isomorphism $J[p\,]\times P[p\,]
\xrightarrow{\sim } \widetilde{J}[p\,]$. In particular, if $J$ is
an ordinary abelian variety, then $\widetilde{J}$ is ordinary if
and only if $P$ is ordinary.
\end{prop}
\dem Because of the right-hand diagram in  the previous lemma, the
Prym variety $$P_1:=\text{ker}(\text{Nm})^0 \subseteq
\widetilde{J}_1$$ coincides with the $p$-twist of
$P:=\text{ker}(\text{Nm})^0 \subseteq \widetilde{J}$ and it is
mapped by $V$ onto $P$. Furthermore, the restriction $V_{|P_1} : P_1
\to P$ being the pull-back of particular line bundles over
$\widetilde{X}_1$ by the relative Frobenius, it coincides with the
Verschiebung $V : P_1 \to P$ for $P$. Therefore, the commutation of
the diagram in the proposition follows from the left-hand
commutative diagram in the lemma and from the fact that $V$ is a homomorphism.\\
One knows that $\text{ker}\, \sigma \subseteq J[2] \times P[2]$ and,
as $A[p\,]\cap A[2]=\{0\}$ for any abelian variety, we have the
isomorphism announced. Recalling that $X$ was supposed to be
ordinary, the last assertion follows from the induced isomorphism
$$J[p\,]_{\rm red}\times P[p\,]_{\rm red}
\xrightarrow{\sim } \widetilde{J}[p\,]_{\rm red}$$ on the reduced
parts of these group schemes. $\square$\\

The following result enables us to apply the rsults gathered in
the Lemma 4.1 to both $\widetilde{J}$ and $P$ for any sufficiently
general curve $X$.

\begin{prop}[B. Zhang] Let $X$ be a general, proper and smooth connected curve
over an algebraically closed field of characteristic $p$ and let $f
: Y \to X$ be an étale cover with abelian Galois group G. Then $Y$
is ordinary.\end{prop} \dem [Zh]\\

 Let us
investigate a bit further in the relationship between $P$ and $P_1$.
Choose an element $z$ in $S_{\tau}=\{z\in J|\,z^2=\tau \} \subset
J[4]$. In particular, it determines the image of the map
$\delta_{\tau,\,z} : P \to |2\Theta|$, namely, one of the two
$\tau$-invariant projective line in $|2\Theta|$.

 As $[p\,] : J \to
J$ induces $[-1]^{\frac{p-1}{2}}$ on $J[4]$, $F$ and $V$ induce
isomorphisms between $J[4]$ and $J_1[4]$, and we let $z_1$ be the
isomorphic image of $z$ via $F$. Thus,
$F^*z_1=V(z_1)=(-1)^{\frac{p-1}{2}}\, z$ and as
$$F^*((\pi_1)_*(q_1) \tens
z_1)\cong \pi_*(F^*(q_1)) \tens F^*z_1$$ for any $q_1$ of $P_1$
(Lemma 4.10), one sees (Lemma 3.4 and Remark 3.5) that the following
diagram commutes

\begin{eqnarray} \end{eqnarray}
\begin{center}
\unitlength=0.6cm
\begin{picture}(10,1)
\put(0.3,0){\put(0.8,0){\put(1,2.5){$P$}
 \put(0.9,0){$
P_1$}\put(1.1,0.7){\vector(0,1){1.5}}

\put(0.8,0){\put(5.6,0){$M_{X_1}$} \put(5.8,2.5){$M_{X}$}
\put(1.2,2.7){\vector(1,0){4}} \put(1.2,0.2){\vector(1,0){4}}
\dashline{0.2}(6.4,0.7)(6.4,2.1)\put(6.4,2.1){\vector(0,1){0.2}}}}

\small \put(4.2,3){$d_{\tau,\,z}$} \put(4.,0.5){$d_{\tau_1,\,z_1}$}
\put(8.4,1.1){$V$}} \small \put(1.5,1.1){$V$}
\end{picture} \end{center}

Let $\Xi$ (resp. $\Xi_1$) be a symmetric representative for the
principal polarization of $P$ (resp. $P_1$). Upon choosing $\Xi$ and
$\Xi_1$ suitably, one can ask that $\O(p\,\Xi_1) \cong V^*\O(\Xi)$
(Lemma 4.1). Let $\varphi_{2\Xi}$ (resp. $\varphi_{2\Xi_1}$) be the
canonical map $$P \to \P H^0(P,\,\O(2\Xi))^*\ (\text{resp. }P_1 \to
\P H^0(P_1,\,\O(2\Xi_1))^*)$$ Because $V$ commutes with $[-1]$, it
induces a map $\widetilde{V} $ such that the following diagram
commutes
\begin{center} \unitlength=0.6cm
\begin{picture}(12,4)
\put(0,0){$\P H^0(P_1,\,\O(2\Xi_1))$} \put(8.4,0){$\P
H^0(P,\,\O(2\Xi))$} \put(1.3,1.7){$\varphi_{2\Xi_1}$}
\put(1.5,0){\put(0.8,3){$P_1$} \put(8.9,3){$P$}
\put(2,3.2){\vector(1,0){6}} \put(1.2,2.7){\vector(0,-1){1.8}}
\put(9,2.7){\vector(0,-1){1.8}} \put(4,0.2){\vector(1,0){2.6}}
\small \put(5,3.4){$V$} \put(5,0.6){$\widetilde{V}$}
 \put(9.4,1.7){$\varphi_{2\Xi}$}}
\end{picture}
\end{center}
In other words, $V^*  : H^0(P,\,\O(2\Xi)) \to
H^0(P_1,\,\O(2p\,\Xi_1))$ factors as the composite
\begin{eqnarray}H^0(P,\,\O(2\Xi)) \xrightarrow{\widetilde{V}^*} {\rm
Sym}^p\,H^0(P_1,\,\O(2\Xi_1)) \to
H^0(P_1,\,\O(2p\,\Xi_1))\end{eqnarray}

Let $\wt$ (resp $\wt_1$) be the lifting of $\tau$ (resp. $\tau_1$)
in $\mathcal{H}$ (resp. $\mathcal{H}_1$) corresponding to our choice
of $z$ in $S_\tau$ (resp. to $z_1=F(z)$ in $S_{\tau_1}$).  Consider
the basis $\{\Lambda_0,\,\Lambda_1\}$ of $W^{\wt}$ constructed in
the Proposition 2.12 and construct the basis
 $\{\Gamma_0,\,\Gamma_1\}$ of
$H^0(P,\,\O(2\Xi))$  by means of the isomorphism $\chi$, as it is
done in the proof of the Proposition 3.3. Upon normalizing $\chi_1$
suitably, the two bases
$$\{\Lambda^{(p)}_0,\, \Lambda^{(p)}_1\}\  (\text{of }
W_1^{\wt_1})\text{ and }\{\Gamma^{(p)}_0,\,\Gamma^{(p)}_1\}\
(\text{of }H^0(P_1,\,\O(2\Xi_1)))$$ obtained via $i^*$ correspond
one to the other via $\chi_1$.\\
Now, let $Q_0,\, Q_1$
 be the degree $p$ homogeneous polynomials such that
\begin{eqnarray}\widetilde{V}^*(\Gamma_i)=Q_i(\underline{\Gamma}^{(p)})\end{eqnarray}
These polynomials depend only on the elliptic curve $P$, hence on
$\tau$, and can be explicitly computed (see Section 5). The diagram
4.4, together with the Proposition 3.3.(2), gives the following
commutative diagram

\begin{eqnarray} \end{eqnarray}
\begin{center}
\unitlength=0.6cm
\begin{picture}(16,2)
\put(0,0.5){\put(0,0){ \put(0.4,2.5){$\P H^0(P,\,\O(2\Xi))$}
 \put(0.2,0){$\P  H^0(P_1,\,\O(2\Xi_1))$}
 \put(2.5,0.7){\vector(0,1){1.5}} \small \put(1.5,1.1){$\widetilde{V}$} }
\put(11.6,0){ \put(0.1,2.5){$\subset \ |2\Theta|$}
 \put(0,0){$\subset \ |2\Theta_1|$}
 \dashline{0.2}(1.8,0.7)(1.8,2.1) \put(1.8,2.1){\vector(0,1){0.1}} }
\put(4.2,0){\put(5.4,0){$\P W_1^{\wt_1}$} \put(5.5,2.5){$\P
W^{\wt}$} \put(1.2,2.7){\vector(1,0){3.9}} \put(3,2.8){$\sim$}
\put(1.4,0.2){\vector(1,0){3.7}} \put(3,0.3){$\sim$} } \small
\put(13.7,1.2){$V$}}
\end{picture} \end{center}
Therefore, there corresponds to $\widetilde{V}$ a morphism $\P
W_1^{\wt_1} \to \P W^{\wt}$ still denoted by $\widetilde{V}$ and
letting $\{\lambda_\bullet\}$ (resp. $\{\lambda^{(p)}_\bullet\}$) be
the basis of $(W^{\wt})^*$ (resp. $(W_1^{\wt_1})^*$), dual to the
basis $\{\Lambda_\bullet\}$ (resp. $\{\Lambda_\bullet^{(p)}\}$,
obtained thanks to the Lemma 2.4, one has
$\widetilde{V}(\lambda_i)=Q_i(\underline{\lambda}^{(p)})$ for
$i=0,\,1$.\\

\subsection{Determining the equations of the Verschiebung}

\begin{prop} Let $\lambda_0$ be the image of $x_{00}$ via $W^*
\twoheadrightarrow (W^{\wt})^*$. Via the canonical surjection
$${\rm Sym}\,^p\,W_1^* \to {\rm Sym}\,^p\,(W_1^{\wt_1})^*$$ the
element $V_{00}$ of ${\rm Sym}\,^p\,W_1^*$ (defined at (4.3)) maps
to $$\widetilde{V}^*(\lambda_0)=Q_0(\lambda_\bullet^{(p)})$$
(where $Q_0$ is the polynomial defined at (4.6)).
\end{prop}
\dem Let $K_{V\times V}:=J_1[p\,]_{\rm red} \times P_1[p\,]_{\rm
red}$ denote the kernel of $V\times V$ and let $f : J_1\times P_1
\to \widetilde{J}$ be the diagonal map $\sigma \circ (V\times V) =V
\circ \sigma_1$ in the diagram of the Proposition 4.12. Its kernel
$K_f$ is isomorphic to $K_{\sigma_1} \times K_{V\times V}$ and we
let
$$\vartheta : K_f\xrightarrow{\sim} \widetilde{K}_f\subseteq
\mathcal{G}(\O(2p\,\Theta_1) \boxtimes \O(2p\,\Xi_1))$$ be the
unique level structure such that $$f_*(\O(2p\,\Theta_1) \boxtimes
\O(2p\,\Xi_1))^{\widetilde{K}_f}\cong \O(\widetilde{\Theta}_z)$$ Let
$s_0$ (resp. $s_1$) be the unique (up to scalar) non zero section of
$\O(\widetilde{\Theta}_z)$ (resp. $\O(\widetilde{\Theta}_{z_1})$.
Upon normalizing suitably $s_1$, one has
$$V^*(s_0)^p=[p\,]^*(s_1)$$ or, equivalently, $s_1=i^*(s_0)$ (see
Remark 4.7 for an analogous fact). The morphism $\delta_{\tau,\,z}$
(resp. $\delta_{\tau_1,\,z_1}$) is defined by $\sigma^*(s_0)$ (resp.
$\sigma^*(s_1)$) (see Proposition 3.3).\\
 The pull-back $f^*(s_0)$
has to be the unique non zero section (up to scalar) of
$$H^0(J_1,\,\O(2p\,\Theta_1)) \tens H^0(P_1,\,\O(2p\,\Xi_1))$$
invariant under the action of $\widetilde{K}_f$ and it induces an
arrow $$ H^0(J_1,\,\O(2p\,\Theta_1))^* \to
H^0(P_1,\,\O(2p\,\Xi_1))$$ Now, insert this map in  the following
diagram
\begin{center}
\unitlength=0.7cm
\begin{picture}(25,7)
\put(0,0.8){\put(0.4,6){$H^0(J,\,\O(2\Theta))^*$}

\put(0,3){$\left(H^0(J,\,\O(2\Theta))^{\wt}\right)^*$}

\put(0.6,0){$H^0(P,\,\O(2\Xi))$}

\put(8.8,0){${\rm Sym}^p H^0(P_1,\,\O(2\Xi_1))$}

\put(8,3){${\rm
Sym}\,^p\,\left(H^0(J_1,\,\O(2\Theta_1))^{\wt_1}\right)^*$}

\put(8.3,6){${\rm
Sym}\,^p\,\left(H^0(J_1,\,\O(2\Theta_1))\right)^*$}

\put(18,6){$H^0(J_1,\,\O(2p\,\Theta_1))^*$}

\put(18.2,0){$H^0(P_1,\,\O(2p\,\Xi_1))$}

\put(5,6.2){\vector(1,0){3}} \put(5,0.2){\vector(1,0){3}}
\put(15,6.2){\vector(1,0){2.5}} \put(15,0.2){\vector(1,0){2.5}}

 \put(2.2,5.7){\vector(0,-1){2}}  \put(2.2,5.7){\vector(0,-1){1.8}}
\put(2.2,2.7){\vector(0,-1){2}} \put(11.3,5.7){\vector(0,-1){2}}
\put(11.3,5.7){\vector(0,-1){1.8}} \put(11.3,2.7){\vector(0,-1){2}}
\put(20.2,5.7){\vector(0,-1){5}} \put(2.5,1.6){$\wr$}
\put(11.6,1.6){$\wr$}}
\end{picture}\end{center}
where
\begin{itemize}
 \item the top line is the factorization (4.2) of $V^* : H^0(J,\,\O(2\Theta))^* \to
H^0(J_1,\,\O(2p\,\Theta_1))^* $
\item the bottom line
 is the factorization (4.5) of  $V^* : H^0(J,\,\O(2\Xi)) \to
H^0(J_1,\,\O(2p\,\Xi_1))$) \item the isomorphisms are deduced from
$\chi$ and $\chi_1$ (see Proposition 3.3)\item the surjections are
deduced from the canonical restriction maps.
\end{itemize}
It is enough to show that the left-hand square is commutative. Note
that the composite map in the left-hand column is the one induced by
$\sigma^* (s_0)$. Therefore, the big square is commutative for
$f^*(s_0)=(V\times V)^*(\sigma^* (s_0))$. Note as well that the
composite map in the middle column is the one induced by
$\sigma_1^*(s_1)$ by taking symmetric $p$-powers. Because
$$\sigma_1^*((s_1)^p)=\sigma_1^*(V^*(s_0))= (V\times V)^*(\sigma^*
(s_0))=f^*(s_0)$$ the right-hand square commutes as well. The
following lemma ends the proof. $\square$\\

\begin{lem} The evaluation map ${\rm Sym}^p\,H^0(P_1,\,\O(2\Xi_1)) \to
H^0(P_1,\,\O(2p\,\Xi_1))$ is injective.\end{lem} \dem On the one
hand, as the evaluation map $H^0(P_1,\,\O(2\Xi_1)) \tens \O_{P_1}
\to \O(2\Xi_1)$ corresponds to the canonical map $\varphi_{2\Xi_1}$,
the map
$${\rm Sym}^p\,H^0(P_1,\,\O(2\Xi_1)) \tens \O_{P_1} \to
\O(2p\,\Xi_1)$$ corresponds to the composite $$P_1
\xrightarrow{\varphi_{2\Xi_1}} \P^1 \xrightarrow {\rho_p} \P^p$$
where $\rho_p$ is the $p$-uple embedding. On the other hand, because
$\O(2p\,\Xi_1)$ is very ample for any $p>2$, the map
$$H^0(P_1,\,\O(2p\,\Xi_1)) \tens \O_{P_1} \to \O(2p\,\Xi_1)$$
corresponds to the embedding $P_1 \hookrightarrow \P
H^0(P_1,\,\O(2p\,\Xi_1))^*$. Now, the evaluation map $${\rm
Sym}^p\,H^0(P_1,\,\O(2\Xi_1)) \to H^0(P_1,\,\O(2p\,\Xi_1))$$ induces
a map $$\P H^0(P_1,\,\O(2p\,\Xi_1))^* \to \P ^p$$ and the former is
injective if and only if the latter is non-degenerate. But all these
maps fit into the commutative diagram
\begin{center}
\unitlength=0.6cm
\begin{picture}(12,4)
\put(0.1,0){$\P^1$} \put(0,3){$P_1$} \put(8,0){$\P^p$}
\put(5.4,3){$\P H^0(P_1,\,\O(2p\,\Xi_1))^* $}
\put(1,0.2){\vector(1,0){6}} \put(1,3.2){\vector(1,0){4}}
\put(0.3,2.7){\vector(0,-1){2}} \put(8.2,2.7){\vector(0,-1){2}}
\end{picture}
\end{center}
and the required injectivity comes from the fact that the image of
the $p$-uple embedding is non-degenerate. $\square$\\

This enables us to state the main result of this paper :

\begin{theo} Let $X$ be a smooth, proper, curve of genus 2, sufficiently general, over an
algebraically closed field of characteristic $p=3,\,5$ or 7. The
generalized Verschiebung $V : M_{X_1} \dashrightarrow M_X$ is
completely determined by its restriction to the projective lines
that are invariant under the action of a non zero element of $J[2]$.
\end{theo}

\begin{cor} For $p=3,\,5,\,7$, one can compute the homogeneous
degree $p$ polynomials $V_i$ ($i=00,\,01,\,10$ and $11$) such that
$$\begin{array}{cccl}
  \widetilde{V} : & |2\Theta_1| & \dashrightarrow & |2\Theta|\\
& y_i & \mapsto & V_i(\underline{y}) \end{array}$$
\end{cor}

\subsubsection{Proof of the Theorem 4.16 : A digression in combinatorial algebra.}

Using Remark 3.5 and Lemma 4.5, we find that the direct sum ${\rm
Sym}^p\,(W_1^{\wt_1})^*\oplus {\rm Sym}^p\,(W_1^{-\wt_1})^*$ (which
depends only on the choice of $\tau$) is endowed with an action of
$\mathcal{H}$ (of weight $p^2$) and the canonical map ${\rm
Sym}^p\,W_1^*\to {\rm Sym}^p\,(W_1^{\wt_1})^*\oplus {\rm
Sym}^p\,(W_1^{-\wt_1})^*$ is equivariant for the action of
$\mathcal{H}$ on both spaces. Taking all order 2 elements of $J[2]$
together, we find a morphism of $\mathcal{H}$-representations
$$ R_P : {\rm
Sym}\,^p\,W_1^* \to BG:=\bigoplus_{\tau \in\, J[2]\setminus
\{0\}}{\rm Sym}^p\,(W_1^{\wt_1})^*\oplus {\rm
Sym}^p\,(W_1^{-\wt_1})^*$$ Because of the Proposition 4.14, one
knows the image of the irreducible sub-representation $W^* \subset
{\rm Sym}\,^p\,W_1^*$ that determines $V$ (see Proposition 4.9) in
$BG$ and one can ask whether or not these data allow us to determine
completely this sub-representation.

A necessary condition is that the above map $R_P$ is injective. It
cannot be the case for large $p$ since
$$ {\rm
dim}\ \left({\rm Sym}\,^p\,W_1^*\right)=
\left(\begin{array}{c}p+3\\3\end{array} \right)\sim \Frac{p^3}{6}
\text{ for large }p$$ and $$ {\rm dim}\, BG =30(p+1)$$ More
precisely, it
cannot be injective for any prime $p>7$.\\

As the map $R_P$ is $\mathcal{H}$-equivariant, it is injective if
and only if its restriction to the subspace $\left({\rm
Sym}^p\,W_1^*\right)^{\hat{H}}$ is. One has the

\begin{lem} The reunion of
the two families
$$A(p)=\left\{A_{\underline{f}}=y_{00}\left[\prod_{i\in \,H}
y_{i}^{f_{i}}\right]^2,\text{ with
}|\underline{f}|=\Frac{p-1}{2}\right\}$$ and
$$B(p)=\left\{B_{\underline{f}}=y_{01}y_{10}y_{11}\left[\prod_{i\in \,H} y_{i}^{f_{i}}\right]^2,\text{
with }|\underline{f}|=\Frac{p-3}{2}\right\}$$
 (where $\underline{f}$ is, in
both cases, a multi-index $(f_{00},\,f_{01},\,f_{10},\,f_{11})$ with
$|\underline{f}|= \sum f_i$) is a basis for the space
$\left({\rm Sym}^p\,W_1^*\right)^{\hat{H}}$.\\
In particular, there are scalars $a_{\underline{f}}$ and
$b_{\underline{f}}$ such that
\begin{eqnarray}V_{00}=\widetilde{V}^*(x_{00})=\sum_{|\underline{f}|=\frac{p-1}{2}}
a_{\underline{f}} A_{\underline{f}}
+\sum_{|\underline{f}|=\frac{p-3}{2}} b_{\underline{f}}
B_{\underline{f}}\end{eqnarray}
\end{lem} \dem The
subspace $\left({\rm Sym}^p\,W_1^*\right)^{\hat{H}}$ is generated by
the free family of monomials $\prod_{i\in \,H} y_{i}^{e_{i}}$ with
$e_{01}+e_{10}\equiv e_{01}+e_{11}\equiv e_{10}+e_{11}\equiv 0$
(mod. 2) and we can divide it into the two families of the lemma. $\square$\\

Let $\tau=(x,\,x^*)$ be an non zero element of $J[2]$ and
$\wt=(\mu,\,x,\,x^*)$ with $\mu=1$ if $x^*(x)=1$ and $\mu=i$ if
$x^*(x)=-1$. We will use the basis
$\{\Lambda_0(\tau),\,\Lambda_1(\tau)\}$ (resp.
$\{\bar{\Lambda}_0(\tau),\,\bar{\Lambda}_1(\tau)\}$) of $W^{\wt}$
(resp. $W^{-\wt}$) constructed in the Proposition 2.12. Taking
account of the fact that we are working over the $p$-twist $J_1$ of
$J$, one can give the images of the $y_\bullet$ via the restriction
maps $$W_1^* \to (W_1^{\wt_1})^*\  \text{(resp. }W_1^* \to
(W_1^{-\wt_1})^*)$$ for every $\tau $, in terms of the
$\lambda_\bullet^{(p)}(\tau_1)$ (resp. the
$\bar{\lambda}_\bullet^{(p)}(\tau_1)$) and one can then deduce the
images of the $A_{\underline{f}}$ and the $B_{\underline{f}}$ in
$${\rm Sym}^p\,(W_1^{\wt_1})^*\ \text{(resp. }{\rm
Sym}^p\,(W_1^{-\wt_1})^*)$$ For sake of readability, we shall write
$\lambda_\bullet$ (resp. $\bar{\lambda}_\bullet$) instead of
$\lambda_\bullet^{(p)}(\tau_1)$ (resp.
$\bar{\lambda}_\bullet^{(p)}(\tau_1)$).\\

If $x=00$, the $A_{\underline{f}}$ map to 0 in ${\rm
Sym}^p\,(W_1^{-\wt_1})^*$ for $y_{00}$ maps to 0 in
$(W_1^{-\wt_1})^*$,  and the $B_{\underline{f}}$ map to 0  in both
${\rm Sym}^p\,(W_1^{\wt_1})^*$ and ${\rm Sym}^p\,(W_1^{-\wt_1})^*$
for at least one of the three $y_{01},\,y_{10}$ and $y_{11}$ maps
to 0 in the corresponding space $(W_1^{\wt_1})^*$ or
$(W_1^{-\wt_1})^*$. In the following chart, we have gathered the
images of the $A_{\underline{f}}$ in ${\rm
Sym}^p\,(W_1^{\wt_1})^*$.

\begin{eqnarray}\begin{array}{|c|c|} \hline

\tau & A_{\underline{f}}\\
& \\
 \hline

\hline 0001 & \lambda_0^{2f_{00}+1}\lambda_1^{f_{10}}\text{ if
}f_{01}=f_{11}=0,\ 0 \text{ else.}
\\
\hline 0010 & \lambda_0^{2f_{00}+1}\lambda_1^{f_{01}}\text{ if
}f_{10}=f_{11}=0,\ 0 \text{ else.}
\\
\hline 0011 & \lambda_0^{2f_{00}+1}\lambda_1^{f_{11}}\text{ if
}f_{01}=f_{10}=0,\ 0 \text{ else.}
\\
\hline \end{array}\end{eqnarray}

If $x\neq 00$, upon identifying $\lambda_\bullet$ and
$\bar{\lambda}_\bullet$, the $A_{\underline{f}}$ (resp. the
$B_{\underline{f}}$) have the same images in both ${\rm
Sym}^p\,(W_1^{\wt_1})^*$ and ${\rm Sym}^p\,(W_1^{-\wt_1})^*$.
Straightforward calculations give the  results gathered in the
following chart (4.10) on the next page.\\

Now, write an element of $\left({\rm
Sym}^p\,W_1^*\right)^{\widetilde{\hat{H}}}$ under the form
$$\sum_{|\underline{f}|=\frac{p-1}{2}} a_{\underline{f}} A_{\underline{f}} +\sum_{|\underline{f}|=\frac{p-3}{2}}
b_{\underline{f}} B_{\underline{f}}$$ and suppose that it is in the
kernel of $R_P$.\\
Because of the computations summed up in (4.9), one has
$a_{\underline{f}}=0$ as soon as two of the
$f_{01},\,f_{10},\,f_{11}$ are zero.

\begin{eqnarray}\begin{array}{|c|c|c|}\hline

\tau & A_{\underline{f}} & B_{\underline{f}}\\
& & \\

\hline

 \hline 0100 &
\lambda_0^{1+2(f_{00}+f_{01})}\lambda_1^{2(f_{10}+f_{11})} &
\lambda_0^{1+2(f_{00}+f_{01})}\lambda_1^{2(f_{10}+f_{11}+1)}
\\
\hline 0101 &
(-1)^{f_{01}+f_{11}}\lambda_0^{1+2(f_{00}+f_{01})}\lambda_1^{2(f_{10}+f_{11})}
& (-1)^{1+f_{01}+f_{11}}
\lambda_0^{1+2(f_{00}+f_{01})}\lambda_1^{2(f_{10}+f_{11}+1)}
\\
\hline 0110 &
\lambda_0^{1+2(f_{00}+f_{01})}\lambda_1^{2(f_{10}+f_{11})} &
-\lambda_0^{1+2(f_{00}+f_{01})}\lambda_1^{2(f_{10}+f_{11}+1)}
\\
\hline 0111 &
(-1)^{f_{01}+f_{11}}\lambda_0^{1+2(f_{00}+f_{01})}\lambda_1^{2(f_{10}+f_{11})}
&
(-1)^{f_{01}+f_{11}}\lambda_0^{1+2(f_{00}+f_{01})}\lambda_1^{2(f_{10}+f_{11}+1)}
\\
\hline

\hline 1000 &
\lambda_0^{1+2(f_{00}+f_{10})}\lambda_1^{2(f_{01}+f_{11})} &
\lambda_0^{1+2(f_{00}+f_{10})}\lambda_1^{2(f_{01}+f_{11}+1)}
\\
\hline 1001 &
\lambda_0^{1+2(f_{00}+f_{10})}\lambda_1^{2(f_{01}+f_{11})} &
-
\lambda_0^{1+2(f_{00}+f_{10})}\lambda_1^{2(f_{01}+f_{11}+1)}
\\
\hline 1010 &(-1)^{f_{10}+f_{11}}
\lambda_0^{1+2(f_{00}+f_{10})}\lambda_1^{2(f_{01}+f_{11})}&
(-1)^{1+f_{10}+f_{11}}\lambda_0^{1+2(f_{00}+f_{10})}\lambda_1^{2(f_{01}+f_{11}+1)}
\\
\hline 1011 &
(-1)^{f_{10}+f_{11}}\lambda_0^{1+2(f_{00}+f_{10})}\lambda_1^{2(f_{01}+f_{11})}
&
(-1)^{f_{10}+f_{11}}\lambda_0^{1+2(f_{00}+f_{10})}\lambda_1^{2(f_{01}+f_{11}+1)}
\\
\hline

\hline 1100 &
\lambda_0^{1+2(f_{00}+f_{11})}\lambda_1^{2(f_{01}+f_{10})} &
\lambda_0^{1+2(f_{00}+f_{11})}\lambda_1^{2(f_{01}+f_{10}+1)}
\\
\hline 1101 &
(-1)^{f_{10}+f_{11}}\lambda_0^{1+2(f_{00}+f_{11})}\lambda_1^{2(f_{01}+f_{10})}
&
(-1)^{f_{10}+f_{11}}\lambda_0^{1+2(f_{00}+f_{11})}\lambda_1^{2(f_{01}+f_{10}+1)}
\\
\hline 1110 &(-1)^{f_{10}+f_{11}}
\lambda_0^{1+2(f_{00}+f_{11})}\lambda_1^{2(f_{01}+f_{10})}&
(-1)^{1+f_{10}+f_{11}}\lambda_0^{1+2(f_{00}+f_{11})}\lambda_1^{2(f_{01}+f_{10}+1)}
\\
\hline 1111 &
\lambda_0^{1+2(f_{00}+f_{11})}\lambda_1^{2(f_{01}+f_{10})} &
-\lambda_0^{1+2(f_{00}+f_{11})}\lambda_1^{2(f_{01}+f_{10}+1)}
\\
\hline
\end{array}\end{eqnarray}

\vspace{1cm} \noindent Because of the computations summed up in the
chart (4.10) above, one finds that, for any $0\leq k\leq
\Frac{p-1}{2}$,
$$\sum_{f_{00}+f_{01}=k\text{ and }f_{01}+f_{11} \text{ even
}}a_{\underline{f}}=0;\hspace{1cm} \sum_{f_{00}+f_{01}=k\text{ and
}f_{01}+f_{11} \text{ odd }}a_{\underline{f}}=0$$
$$\sum_{f_{00}+f_{10}=k\text{ and }f_{10}+f_{11} \text{ even
}}a_{\underline{f}}=0;\hspace{1cm} \sum_{f_{00}+f_{10}=k\text{ and
}f_{10}+f_{11} \text{ odd }}a_{\underline{f}}=0$$
$$\sum_{f_{00}+f_{11}=k\text{ and }f_{10}+f_{11} \text{ even
}}a_{\underline{f}}=0;\hspace{1cm} \sum_{f_{00}+f_{11}=k\text{ and
}f_{10}+f_{11} \text{ odd }}a_{\underline{f}}=0$$ and we have
analogous results for the $b_{\underline{f}}$ (with $0\leq k \leq
\Frac{p-3}{2}$).\\

Thus, we have reduced our problem to the following combinatorial
situation : Being given a set of scalars $a_{\underline{f}}$, where
the $\underline{f}$ are four-letters multi-index $\underline{f}$
with $|\underline{f}|=r$, satisfying the $6(r+1)$ relations stated
above, are these scalars meant to be zero ? If they are for
$r=\Frac{p-3}{2}$ but are not for $r=\Frac{p-1}{2}$, does the
indetermination only concern the $a_{\underline{f}}$ for which two
of the $f_{01},\,f_{10},\,f_{11}$ are zero (the indetermination
would therefore vanish because of the additional data deduced from
(4.9)) ?\\

 In fact, the first question has a positive answer for
$r=0,\,1,\,2,\,3$ (we leave the proof of this assertion to the
reader) and
the theorem follows. \\

\section{Computing the equations of $V_2$ for small $p$}

\subsection{Multiplication by $p$ on an elliptic curve}

Let $k$ be an algebraically closed field of characteristic $p\geq 3$
and let $(E,\,q_0)$ be an elliptic curve.

Let us recall briefly how the group law of $E$ can be recovered from
the geometry of the curve (see, e.g., [Sil] for further references
on that question). The sheaf $\O(q_0)$ gives a principal
polarization thus $\O(3q_0)$ is very ample and determines an
embedding $E \hookrightarrow \P^2$. Three points $q_1,\,q_2$ and
$q_3$ on the curve lie on the same projective line in $\P^2$ if and
only if $\O(q_1+q_2 +q_3) \cong \O(3q_0)$. On the other hand, any
projective line in $\P^2$ intersects with $E$ in three points
(counted with multiplicities). It is easily seen that $E$ is
isomorphic to its Jacobian variety by means of $q \mapsto \O(q-q_0)$
and the group law on the latter gives the group law on the former.
Namely, one sets $q_1+q_2=-q_3$, where $q_3$ is the unique point in
$E$ such that $\O(q_1+q_2 +q_3) \cong \O(3q_0)$.

The projection $\P^2 \to \P^1$ from the point $q_0$ induces the
canonical map $E \to \P^1$, which is a ramified double cover, and
the choice of a suitable rational coordinate $x$ on the projective
line allows us to give a birational model
$$y^2=x(x-1)(x-\mu)$$
of $E$, with $\mu$ different from $0$ and $1$.\\

 One can determine explicitly
the group law over $E$ in intersecting this plane curve with lines.
Namely, one has the following duplication and addition formulae
found in [Sil] (III, §2). For convenience, we let $P$ be the
polynomial $x(x-1)(x-\mu)$ and we let $P'$ be its derivative.

\underline{Duplication formula} : Let $q_1$ be a $k$-point on the
curve with coordinates $(x_1,\,y_1)$ and let $q_2$, with coordinates
$(x_2,\,y_2)$, be $[2]\,( q_1)$. The opposite of the latter is the
unique point of $E$ (different from $q_1$) lying on the tangent line
to $E$ at the point $q_1$. This tangent line has equation $y=\alpha
x + \beta$ with
$$\alpha=\Frac{P'(x_1)}{2y_1}=\Frac{3x_1^2-2(\mu+1)x_1+\mu}{2y_1} \ \
{\rm and } \ \ \beta=\Frac{2P(x_1)-x_1 P'(X_1)}{2y_1}=\Frac{x_1(\mu
-x_1^2)}{2y_1}$$ Thus, one has
\begin{eqnarray}\left\{ \begin{array}{rcl}
x_2& =& \alpha^2+(\mu+1)-2 x_1 = \Frac{(P'(x_1)) ^2-4 (2 x_1 -(\mu + 1)) P(x_1)}{4 P(x_1)}\\
& = & \Frac{1}{4}\Frac{(x_1^2-\mu)^2}{x_1(x_1-1)(x_1-\mu)} \\
\\
y_2 & =& -(\alpha x_2 + \beta)\\
\end{array}\right. \end{eqnarray}

\underline{Addition formula} : Let $q_1$ and $q_2$ be two $k$-points
on the curve with coordinates $(x_i,\,y_i)$ $(i=1,\,2)$ and let
$q_3$, with coordinates ($x_3,\,y_3)$, be the sum  $q_1+q_2$.
Suppose that $q_1\neq \pm q_2$, i.e., that $x_1\neq x_2$. The unique
line passing through $q_1$ and $q_2$ has equation
$$y=\alpha x + \beta \ \ {\rm with } \ \alpha=\Frac{y_2-y_1}{x_2-x_1} \ \ {\rm
and } \ \ \beta=\Frac{y_1x_2-y_2x_1}{x_2-x_1}$$ Thus, the third
intersection point of that line and the plane curve being
$-(q_1+q_2)$, $q_3$ has coordinates
\begin{eqnarray}\left\{\begin{array}{rcl}
x_3& =& \alpha^2+(\mu+1)-(x_1+x_2)\\
& = &  \left(\Frac{y_2-y_1}{x_2-x_1}\right)^2+(\mu+1)-(x_1+x_2)\\
y_3& =& -(\alpha x_3 + \beta)\\
\end{array}\right.\end{eqnarray}

Combining these two formulae, we are theoretically able to give the
coordinates $(x_n,\,y_n)$ of the point $q_n=[n](q_1)$, in terms of
$x_1$ and $y_1$, at least for a general point $q_1$. Note that two
opposite points of $E$ collapse in $\P^1$, i.e., that the canonical
$E\to \P^1$ is a quotient under the action of $\{\pm\}$ and that the
branched points $0,\,1,\,\infty$ and $\mu$ of that map are precisely
the order 2 points of $E$. As the action of $\{\pm\}$ commutes with
multiplication by $p$, the latter induces a map $\P^1 \to\P^1$. It
has total degree $p^2$ and separable degree $p$. Hence, if we let
$z$ be $1/x$, so that $\{x,\,z\}$ is a basis for
$H^0(E,\,\O(2q_0))\cong H^0(\P^1,\,\O(1))$, one can find two
homogeneous polynomials of degree $p$ (say $D$ and $N$) such that
the map induced by $[p\,]$ on $\P^1$ is given by
\begin{eqnarray*}
\P^1 & \to & \P^1\\
(x:\,z) & \mapsto & (N(x^p,\,z^p):\,D(x^p,\,z^p))\end{eqnarray*} If
$E_1$ is the $p$-twist of $E$, the map $\P^1 \to \P^1$, induced by
the separable part $V : E_1 \to E$ of multiplication by $p$, is
therefore given by
\begin{eqnarray}
(x^{(p)}:\,z^{(p)}) & \mapsto &
(N(x^{(p)},\,z^{(p)}):\,D(x^{(p)},\,z^{(p)}))\end{eqnarray} where
$x^{(p)}$ and $z^{(p)}$ are the $p$-twisted coordinates of $\P^1$
corresponding to $x$ and $z$ respectively.\\

\underline{Division Polynomials} : In the case $p\geq 5$, [Sil]
gives, as an exercise (Ex. 3.7.), the following formulae, that are
more convenient to implement when trying to determine the
polynomials $N$ and $D$ using a computer. Take an elliptic curve
\begin{eqnarray} y'^2=x'^3+Ax'+B \end{eqnarray} Define $\psi_m$ in $\Z[A,\,B,\,x',\,y']$
inductively by :
\begin{eqnarray} \left\{\begin{array}{l}
\psi_1=1,\ \\
\psi_2=2y', \\
\psi_3=3x'^4+6Ax'^2+12Bx'-A^2,\\
\psi_4=4y'(x'^6+5Ax'^4+20Bx'^3-5A^2x'^2-4ABx'-8B^2-A^3),\\
\psi_{2m+1}=\psi_{m+2}\psi_{m}^3-\psi_{m-1}\psi_{m+1}^3
\ \ (m\geq 2)\\
2y'\psi_{2m}=\psi_m (
\psi_{m+2}\psi_{m-1}^2-\psi_{m-2}\psi_{m+1}^2) \ \ (m\geq
2).\end{array}\right. \end{eqnarray} Define furthermore
\begin{eqnarray}\phi_m=x'\psi_m^2-\psi_{m+1}\psi_{m-1}\end{eqnarray} Then, using the
equation (5.4), one checks that, for $m$ odd, the polynomials
$\phi_m$ and $\psi_m$ (of degree $m^2$ and $m^2-1$ respectively)
lie in fact in $\Z[A,\,B,\,x']$, and that the map $\P^1 \to \P^1$
induced by $[p\,]$ is given by
\begin{eqnarray*}
 \P^1 & \to & \P^1\\
 (x' :\, z')& \mapsto &
\left(
z'^{p^2}\,\phi_p(x'/z'):\,z'^{p^2}\,\psi_p(x'/z')\right)\end{eqnarray*}
where $z'$ is the rational coordinate of $\P^1$ defined by
$z'=1/x'$. We let $N'$ (resp. $D'$) be the homogeneous degree $p$
polynomial such that $N'(x'^p,\,z'^p)=z'^{p^2}\,\phi_p(x'/z')$
(resp. $D'(x'^p,\,z'^p)=z'^{p^2}\,\phi_p(x'/z')$). As above, the
map $\P^1 \to \P^1$ induced by $V : E_1 \to E$ is therefore given
by
\begin{eqnarray*}
\P^1 & \to & \P^1\\
(x'^{(p)}:\,z'^{(p)}) & \mapsto &
(N'(x'^{(p)},\,z'^{(p)}):\,D'(x'^{(p)},\,z'^{(p)}))\end{eqnarray*}
where $x'^{(p)}$ and $z'^{(p)}$ are the $p$-twisted coordinates of
$\P^1$
corresponding to $x'$ and $z'$ respectively.\\

One finds a coordinates change transforming the elliptic curve
$y^2=x(x-1)(x-\mu)$ in the (isomorphic) elliptic curve
$y'^2=x'^3+Ax+B$, e.g.,
$$x'=x-\Frac{\mu +1}{3} z,\ \ z'=z,$$ and the scalars $A$ and $B$ are
therefore $$A=\Frac{\mu^2-\mu+1}{3},\ \
B=\Frac{(\mu+1)^3-3(\mu^3+1)}{27}.$$ Finally, one obtains the
polynomials $N$ and $D$ defined in (5.3) as follows :
\begin{eqnarray}
\left\{ \begin{array}{l} N(x^{(p)},\,z^{(p)}) =
N'\left(x^{(p)}-\Frac{\mu^p +1}{3}z^{(p)},\,z^{(p)}\right) +
\Frac{\mu
+1}{3}\ D'\left(x^{(p)}-\Frac{\mu^p +1}{3}z^{(p)},\,z^{(p)}\right)\\

\\
D(x^{(p)},\,z^{(p)}) = D'\left(x^{(p)}-\Frac{\mu^p
+1}{3^p}z^{(p)},\,z^{(p)}\right)\end{array}\right. \end{eqnarray}\\

Using the Proposition 2.12, the Corollary 2.13 and the Lemma 3.4, we
find that there is a unique linear automorphism of $\P^1$, i.e., an
element of $PGL(k,\,2)$, which maps $(a:\,b)$ to 0, $(a:\,-b)$ to
$1$ and and $(b:\,a)$ to $\infty$. It maps $(b:\,-a)$ to ($\mu:\,1)$
with
$$\mu=\left(\Frac{b^2+a^2}{2ab}\right)^2=\Frac{2-\omega(\tau)}{4}$$
If $x$ and $z$ are the corresponding rational coordinates of $\P^1$,
one has
$$x=-\Frac{1}{a} \lambda_0+\Frac{1}{b}\lambda_1,\ \ z= \Frac{-2a}{a^2+b^2}\lambda_0+\Frac{2b}{a^2+b^2}
\lambda_1,$$ and the elliptic curve $P$ has equation
$y^2=x(x-1)(x-\mu)$.\\

If we let $\lambda_0^{(p)}$ and $\lambda_1^{(p)}$ be the $p$-twisted
coordinates  of $\P^1$ corresponding to $\lambda_0$ and $\lambda_1$,
and if we denote by $Q_0$ and $Q_1$ the homogeneous polynomials of
degree $p$ such that the map $\P^1 \to \P^1$ induced by $V : E_1 \to
E$ is given by
$$(\lambda_0^{(p)}:\, \lambda_1^{(p)}) \mapsto (Q_0(\lambda_0^{(p)},\,
\lambda_1^{(p)}):\,Q_1(\lambda_0^{(p)},\, \lambda_1^{(p)}))$$ then
one has, writing $\lambda_\bullet$ instead of
$\lambda_\bullet^{(p)}$ for sake of readability
\begin{eqnarray}
\begin{array}{rcl}
Q_0(\lambda_0,\, \lambda_1)& =& \Frac{2ba}{a^2+b^2}
N\left(-\Frac{1}{a^p} \lambda_0+\Frac{1}{b^p}\lambda_1,
\Frac{-2a^p}{(a^{2}+b^{2})^p}\lambda_0+\Frac{2b^p}{(a^2+b^2)^p}
\lambda_1\right)\\
& & \ \ - \Frac{1}{b} D\left(-\Frac{1}{a^p}
\lambda_0+\Frac{1}{b^p}\lambda_1,
\Frac{-2a^p}{(a^{2}+b^{2})^p}\lambda_0+\Frac{2b^p}{(a^2+b^2)^p}
\lambda_1\right)\end{array}\end{eqnarray}

\Rq Note that we are only interested in $Q_0$, onto which $V_{00}$
restricts. Furthermore, $Q_1$ can be obtained from $Q_0$ under
the action of a suitable element of $P[2]$.\\

\Rq The final result $Q_0$ should not depend on $a$ and $b$ but only
on the constant $\omega(\tau)=-\Frac{a^4+b^4}{a^2b^2}$.\\

\begin{lem}
With the notations given above, one has\\
$\bullet \ p=3$. \begin{eqnarray} Q_0(\lambda_0,\,
\lambda_1)=\lambda_0^3-\omega\lambda_0\lambda_1^2\end{eqnarray}
\noindent $\bullet \ p=5$.
\begin{eqnarray}Q_0=\lambda_0^5+\omega(\omega^2+2)\lambda_0^3\lambda_1^2+(\omega^2+2)
\lambda_0\lambda_1^4\end{eqnarray} \noindent $\bullet \ p=7$.
\begin{eqnarray}Q_0=\lambda_0^7-2\omega(\omega^4-1)\lambda_0^5\lambda_1^2
+\omega^2(\omega^2-1)(\omega^2-2)
\lambda_0^3\lambda_1^4-\omega(\omega^2-1)\lambda_0\lambda_1^6
\end{eqnarray}
\end{lem} \dem In the case $p=3$, one computes $x_3$ directly using formulae (5.1) and (5.2). One finds
$$x_3=\Frac{x_1^9+2\mu (\mu+1)x_1^6+\mu^2 (\mu+1)^2x_1^3}{((\mu+1)x_1^3+\mu^2)^2}
$$
Thus, one has  $$N(x,\,z)=x(x+\mu(\mu+1)z)^2 \text{   and  }
D(x,\,z)=z((\mu+1)x+\mu^2z)^2$$ This result is consistent with the
fact that, in characteristic 3, the unique supersingular elliptic
curve has parameter $\mu$ equal to $-1$ (see [H], Chapter IV,
Example 4.23.1). A straightforward application of the formula (5.8)
gives
the formula (5.9).\\
In the cases $p=5$ and 7, the computations cannot be worked out by
hand. We use Maple 9 to compute the division polynomials defined in
(5.5) and (5.6), then we apply the formula (5.7) to find
$$N(x,z)= x\left[x^2-\mu(\mu+1)(\mu^2-\mu+1)xz+\mu^4(\mu^2-\mu+1)z^2\right]^2$$ and
$$D(x,z)=z\left[(\mu^2-\mu+1)\left[x^2-\mu^2(\mu+1)xz\right]+\mu^6z^2\right]^2$$
when $p=5$ and \begin{eqnarray*} N(x,z)& =&
x\left[x^3+2\mu(\mu+1)(\mu-2)(\mu-4)(\mu^2+3\mu+1)x^2z\right.\\
& &
\hspace{0.5cm}\left.+\mu^4(\mu+1)^2\mu-2)(\mu-4)(\mu^2+1)xz^2+\mu^9(\mu+1)(\mu-2)(\mu-4)z^3\right]^2\end{eqnarray*}
and
$$D(x,z)=z\left[(\mu+1)(\mu-2)(\mu-4)\left[x^3+\mu^2(\mu+1)(\mu^2+1)x^2z+\mu^6(\mu^2+3\mu+1)xz^2\right]+\mu^{12}z^3\right]^2$$
when $p=7$. These results are consistent with the fact that, in
characteristic $5$ (resp. $7$), the only supersingular elliptic
curves  have parameter $\mu$ equal to $j$ or $-j$ with $j^3=1$
(resp. $-1,\,2$ or $4$) (see [H], Chapter IV, Example 4.23.2 (resp.
4.23.3)). Applying the formula (5.8), we obtain the formulae (5.10)
and (5.11). $\square$\\

\subsection{Equations of $V$ for $p=3$}

Let $V_{00}$ be the $\widetilde{\hat{H}}$-invariant element of the
sub-representation determining $V$. It can be written under the
form (4.8) (see Lemma 4.18) which is, for $p=3$,
$$a_{00}y_{00}^3+a_{01}y_{00}y_{01}^2+a_{10}y_{00}y_{10}^2+a_{11}y_{00}y_{11}^2+b y_{01}y_{10}y_{11}$$
We do not need to determine the Prym varieties for every non zero
$\tau$ in $J[2]$. Doing it in the cases
$\tau=0001,\,\tau=0010,\,\tau=0011$ and $\tau=0100$ is enough.\\
 We
fix $a_{00}=1$ and we obtain, using the formula (5.9) and the
expression (given in the chart (2.5)) of the needed $\omega(\tau)$
in terms of the coefficients $k_\bullet$ of the Kummer surface ${\rm
Kum}_X$, the following :
$$V_{00}(\underline{y})=y_{00}^3+2k_{01}y_{00}y_{01}^2+2k_{10}y_{00}y_{10}^2+2k_{11}y_{00}y_{11}^2+2k_{00}
y_{01}y_{10}y_{11}$$ Then, one can deduce the $V_i$
($i=01,\,10,\,11$) by permuting suitably the coordinate functions
$y_\bullet$ in $V_{00}$ (see Proposition 4.9).\\

Notice that $V_{i}$ is the partial (with respect to $y_{i}$) of a
quartic surface
$$ S+ 2k_{00} P+k_{10}Q_{01}+k_{01}Q_{10}+k_{11}Q_{11}$$
(with
$$\begin{array}{c}S= y_{00}^4+y_{01}^4+y_{10}^4+y_{11}^4;\hspace{1cm}
P=y_{00}y_{01}y_{10}y_{11};\\
Q_{01}=y_{00}^2y_{01}^2+y_{10}^2y_{11}^2;\hspace{1cm}Q_{10}=y_{00}^2y_{10}^2+y_{01}^2y_{11}^2;\hspace{1cm}
Q_{11}=y_{00}^2y_{11}^2+y_{01}^2y_{10}^2.)\end{array}$$ isomorphic
to ${\rm Kum}_X$. Thus, one recovers the second assertion of the :

\begin{theo}[Laszlo-Pauly] Let $X$ be a smooth and projective curve of genus 2
over an algebraically closed filed of characteristic $3$.\\
(1) There is an embedding $\alpha : {\rm Kum}_{X} \hookrightarrow
|2\Theta_1|$ such that the equality of divisors in $|2\Theta_1|$
$$\widetilde{V}^{-1}({\rm Kum}_{X})= {\rm Kum}_{X_1}+2\alpha({\rm
Kum}_{X})$$ holds \emph{scheme-theoretically}.\\
(2) The cubic equations of $\widetilde{V}$ are given by the 4
partial derivatives of the quartic equation of the Kummer surface
$\alpha({\rm Kum}_X) \subseteq |2\Theta_1|$. In other words,
$\widetilde{V}$ is the polar map of the surface $\alpha({\rm
Kum}_X)$.\end{theo} \dem
[LP2], Theorem 6.1.\\

 The inverse image
$\widetilde{V}^{-1}({\rm Kum}_{X})$ can be computed explicitly in
our situation as it is defined by the ideal generated by the
pull-back $\widetilde{V}^*(K)$ of the equation (2.1) of ${\rm
Kum}_{X}$, more precisely by its image via the $k$-linear
homogeneous ring map of degree $p$
$$\widetilde{V}^* : {\rm Sym} W^* \to {\rm Sym} W_1^*$$
In other words, a few more computations enable us to recover the
first assertion of the Theorem. Namely, one knows (see the diagram
4.1) that the equation $K_1$ of ${\rm Kum}_{X_1}$ divides
$\widetilde{V}^*(K)$. Let $Q$ be the exact quotient
$\widetilde{V}^*(K)/K_1$. Using Magma, one computes the square root
of $Q$ (e.g., as the greatest common divisor of the partial
derivative $\Frac{\partial Q}{\partial y_{11}}$ and $Q$). This
homogeneous polynomial furthermore coincides with $K_X$.\\

Note by the way that the base locus $\mathcal{I}$ of the rational
map $\widetilde{V} : |2\Theta_1| \dashrightarrow |2\Theta|$ is
tautologically contained in the zero locus of $\widetilde{V}^*(K)$.
As $\widetilde{V}$ restricts to a morphism on ${\rm Kum}_{X_1}$,
$\mathcal{I}$ is contained in the zero locus of
$Q=\widetilde{V}^*(K)/K_1$ and one checks that it is actually
contained in the zero locus of $A$. In other words, $\mathcal{I}$ is
a reduced zero dimensional sub-scheme of $\alpha({\rm Kum}_X)$ which
coincides furthermore with its singular locus.\\

 \Rq Notice that this theorem is
true for \emph{any} curve $X$ (in particular, with no particular
assumptions concerning its $p$-rank)
whereas our calculations only give the result for a sufficiently general ordinary curve.\\

\subsection{Equations of $V$ for $p=5$}

\begin{prop} Let $X$ be a general proper and smooth curve of genus 2
over an algebraically closed field of characteristic 5. There are
coordinates $\{x_\bullet\}$ and $\{y_\bullet\}$ for $|2\Theta|$ and
$|2\Theta_1|$ respectively such that the Kummer surface ${\rm
Kum}_X$ in $|2\Theta|$ has an equation of the form (2.1) and such
that, if the polynomials ($V_\bullet)$ define $\widetilde{V} :
|2\Theta_1| \dashrightarrow |2\Theta|,\ (y_i) \mapsto
(V_i(\underline{y}))$, then
\begin{eqnarray*} V_{00}& = &
y_{00}^5+a_{1100}y_{00}^3y_{01}^2+a_{1010}y_{00}^3y_{10}^2+a_{1001}y_{00}^3y_{11}^2
+a_{0200}y_{00}y_{01}^4+a_{0110}y_{00}y_{01}^2y_{10}^2\\
& & +a_{0101}y_{00}y_{01}^2y_{11}^2+a_{0020}y_{00}y_{01}^4
+a_{0011}y_{00}y_{10}^2y_{11}^2+a_{0002}y_{00}y_{11}^4\\
& &
+b_{00}y_{00}^2y_{01}y_{10}y_{11}+b_{01}y_{01}^3y_{10}y_{11}+b_{10}y_{01}y_{10}^3y_{11}+b_{11}y_{01}y_{10}y_{11}^3
\end{eqnarray*}
with $$\begin{array}{lclcl}
a_{1100}=k_{01}(k_{01}^2+2), & & a_{1010}=k_{10}(k_{10}^2+2), & & a_{1001}=k_{11}(k_{11}^2+2),\\
a_{0200}=(k_{01}^2+2), & & a_{0020}=(k_{10}^2+2),& &
a_{0002}=(k_{11}^2+2),\\ \end{array}$$
$$
\begin{array}{rcl}
a_{0110}& = &3k_{11}(k_{00}^2+k_{11}^2)+k_{01}k_{10}(1-k_{11}^2),\\
a_{0101}& =& 3k_{10}(k_{00}^2+k_{10}^2)+k_{01}k_{11}(1-k_{10}^2),\\
a_{0011}& = & 3k_{01}(k_{00}^2+k_{01}^2)+k_{10}k_{11}(1-k_{01}^2),\\
b_{00}& = & 3k_{00}(k_{00}^2+1)+k_{00}k_{01}k_{10}k_{11},\\
\end{array}$$
$$b_{01}=k_{00}(k_{01}+3k_{10}k_{11}),\ b_{10}=k_{00}(k_{10}+3k_{01}k_{11}),\
b_{11}=k_{00}(k_{11}+3k_{01}k_{10})$$where the $k_\bullet$ are the
coefficients of the equation (2.1) of ${\rm Kum}_X$. The $V_i$
($i=01,\,10,\,11$) can be deduced from $V_{00}$ by a suitable
permutation of the coordinate functions $y_\bullet$, namely the
unique pairwise permutation that exchanges $y_{00}$ and $y_i$.
\end{prop}
\dem Define
$$\begin{array}{l}
\alpha_\tau=\omega(\tau)(\omega(\tau)^2+2)\\
\beta_\tau=\omega(\tau)^2+2\\
\end{array}$$ so that the
formula (5.10) can be written
$$\lambda_0^5+\alpha_\tau\lambda_0^3\lambda_1^2+\beta_\tau
\lambda_0\lambda_1^4$$ Using the equation (4.8) of the Lemma 4.10,
normalized by the condition $a_{2000}=1$, one can look for $V_{00}$
under the form given in the proposition. \\
Using the data (4.9) for $\tau=0001$, one finds that the two
equations
$$\lambda_0^5+a_{1010}\lambda_0^3\lambda_1^2+a_{0020}\lambda_0\lambda_1^4$$
and
$$\lambda_0^5+\alpha_{0001}\lambda_0^3\lambda_1^2+\beta_{0001}\lambda_0\lambda_1^4$$
coincide up to a multiplicative scalar. Therefore, one obtains
$$\begin{array}{lcl}
a_{1010}=\alpha_{0001}, & & a_{0020}=\beta_{0001},\\
\end{array}$$
Similarly, using the data (4.9) for $\tau=0010$ and $0011$
respectively, one finds
$$\begin{array}{lcl}
a_{1010}=\alpha_{0001}, & & a_{0020}=\beta_{0001},\\
a_{1100}=\alpha_{0010}, & & a_{0200}=\beta_{0010},\\
a_{1001}=\alpha_{0011}, & & a_{0002}=\beta_{0011}.\\
\end{array}$$
Now, using the data (4.10) for $\tau=0100$, one finds that the the
two equations
$$(1+a_{1100}+a_{0200})\lambda_0^5+(a_{1010}+a_{1001}+a_{0110}+a_{0101}+b_{00}+b_{01})\lambda_0^3\lambda_1^2+(a_{0020}+a_{0011}+a_{0002}+b_{10}+b_{11})\lambda_0\lambda_1^4$$
and
$$\lambda_0^5+\alpha_{0100}\lambda_0^3\lambda_1^2+\beta_{0100}\lambda_0\lambda_1^4$$
coincide up to a multiplicative scalar. Therefore, one obtains
$$\left\{\begin{array}{l}
a_{1010}+a_{1001}+a_{0110}+a_{0101}+b_{00}+b_{01}=  (1+a_{1100}+a_{0200})\alpha_{0100}\\
a_{0020}+a_{0011}+a_{0002}+b_{10}+b_{11}=
(1+a_{1100}+a_{0200})\beta_{0100}\end{array}\right.$$ Similarly,
using the data (4.10) for $\tau=0101,\,0110$ and $0111$
respectively, one finds
$$\begin{array}{l}\left\{\begin{array}{l}
a_{1010}-a_{1001}-a_{0110}+a_{0101}-b_{00}+b_{01}=  (1-a_{1100}+a_{0200})\alpha_{0101}\\
a_{0020}-a_{0011}+a_{0002}-b_{10}+b_{11}=
(1-a_{1100}+a_{0200})\beta_{0101}\end{array}\right.\\
\left\{\begin{array}{l}
a_{1010}+a_{1001}+a_{0110}+a_{0101}-b_{00}-b_{01}=  (1+a_{1100}+a_{0200})\alpha_{0110}\\
a_{0020}+a_{0011}+a_{0002}-b_{10}-b_{11}=
(1+a_{1100}+a_{0200})\beta_{0110}\end{array}\right.\\
\left\{\begin{array}{l}
a_{1010}-a_{1001}-a_{0110}+a_{0101}+b_{00}-b_{01}=  (1-a_{1100}+a_{0200})\alpha_{0111}\\
a_{0020}-a_{0011}+a_{0002}+b_{10}-b_{11}=
(1-a_{1100}+a_{0200})\beta_{0111}\end{array}\right. \end{array}$$
Combining these results, one can express the $a_\bullet$ as well
as the $b_\bullet$ in terms of the $\alpha_\tau$ and the
$\beta_\tau$. Finally, we use the data (2.5) to express the
$\alpha_\bullet$ and the $\beta_\bullet$ in terms of the
coefficients of ${\rm Kum}_X$ and Maple 9 gives expressions that,
up to a multiple of the equation
(2.2) between the $k_\bullet$, are those stated in the Proposition. $\square$\\

Using Magma to exploit these formulae, one can show the following
corollary :

\begin{cor} There is a degree $2p-2=8$ hypersurface $S$ is $|2\Theta_1|$ such that
the equality of divisors in $|2\Theta_1|$
$$\widetilde{V}^{-1}({\rm Kum}_{X})= {\rm Kum}_{X_1}+2 S $$ holds
\emph{ scheme-theoretically}.\end{cor} \dem Define the field $L$ as
the extension
$$\mathbb{F}_p(k_{01},\,k_{10},\,k_{11})[k_{00}]/(k_{00}^2-k_{01}^2-k_{10}^2-k_{11}^2+k_{01}k_{10}k_{11}+4)$$
of the prime field $\mathbb{F}_p$, and define $R$ as the $L$ vector
space generated by $y_{00},\,y_{01},\,y_{10}$ and $y_{11}$. The
homogeneous polynomials $V_{00},\,V_{01},\,V_{10}$ and $V_{11}$
define a $L$-linear ring homomorphism $\widetilde{V}^* : R \to R$
(defined by $\widetilde{V}^*(y_i)=V_i$). Letting $K$ (resp. $K_1$)
be the equation (2.1) of the Kummer surface ${\rm Kum}_X$ in
$|2\Theta|$ (resp. ${\rm Kum}_{X_1}$ in $|2\Theta_1|$), Magma checks
that $K_1$ divides $\widetilde{V}^*(K)$. Letting $Q$ be the exact
quotient $\widetilde{V}^*(K)/K_1$, Magma checks that it is a square.
$\square$

\subsection{Equations of $V$ for $p=7$}

\begin{prop} Let $X$ be a general proper and smooth curve of genus 2
over an algebraically closed field of characteristic 7. There are
coordinates $\{x_\bullet\}$ and $\{y_\bullet\}$ for $|2\Theta|$ and
$|2\Theta_1|$ respectively such that the Kummer surface ${\rm
Kum}_X$ in $|2\Theta|$ has an equation of the form (2.1) and such
that, if the polynomials ($V_\bullet)$ define $\widetilde{V} :
|2\Theta_1| \dashrightarrow |2\Theta|,\ (y_i) \mapsto
(V_i(\underline{y}))$, then
\begin{eqnarray*} V_{00}& = &
y_{00}^7+a_{2100}y_{00}^5y_{01}^2+a_{2010}y_{00}^5y_{10}^2+a_{2001}y_{00}^5y_{11}^2+a_{1200}y_{00}^3y_{01}^4
+a_{1110}y_{00}^3y_{01}^2y_{10}^2\\
& & +a_{1101}y_{00}^3y_{01}^2y_{11}^2+a_{1020}y_{00}^3y_{10}^4
+a_{1011}y_{00}^3y_{10}^2y_{11}^2+a_{1002}y_{00}^3y_{11}^4\\
& &
+a_{0300}y_{00}y_{01}^6+a_{0210}y_{00}y_{01}^4y_{10}^2+a_{0201}y_{00}y_{01}^4y_{11}^2+a_{0120}y_{00}y_{01}^2y_{10}^4
+a_{0111}y_{00}y_{01}^2y_{10}^2y_{11}^2\\
& & +a_{0102}y_{00}y_{01}^2y_{11}^4+a_{0030}y_{00}y_{10}^6+a_{0021}y_{00}y_{10}^4y_{11}^2+a_{0012}y_{00}y_{10}^2y_{11}^4+a_{0003}y_{00}y_{11}^6\\
& &
+b_{2000}y_{00}^4y_{01}y_{10}y_{11}+b_{1100}y_{00}^2y_{01}^3y_{10}y_{11}+b_{1010}y_{00}^2y_{01}y_{10}^3y_{11}+b_{1001}y_{00}^2y_{01}y_{10}y_{11}^3\\
& & +b_{0200}y_{01}^5y_{10}y_{11}+b_{0110}y_{01}^3y_{10}^3y_{11}+b_{0101}y_{01}^3y_{10}y_{11}^3\\
& &
+b_{0020}y_{01}y_{10}^5y_{11}+b_{0011}y_{01}y_{10}^3y_{11}^3+b_{0002}y_{01}y_{10}y_{11}^5
\end{eqnarray*}
with $$\begin{array}{lclcl}
a_{2100}=-2k_{01}(k_{01}^4-1), & & a_{2010}=-2k_{10}(k_{10}^4-1), & & a_{2001}=-2k_{11}(k_{11}^4-1),\\
a_{1200}=k_{01}^2(k_{01}^2-1)(k_{01}^2-2), & &
a_{1020}=k_{10}^2(k_{10}^2-1)(k_{10}^2-2),& & a_{1002}=k_{11}^2(k_{11}^2-1)(k_{11}^2-2),\\
a_{0300}=-k_{01}(k_{01}^2-1),& &a_{0030}=-k_{10}(k_{10}^2-1),&
&a_{0003}=-k_{11}(k_{11}^2-1),\\
\end{array}$$

$$
\begin{array}{rcl}
a_{0111}& = &
3k_{00}^6+2(k_{01}^6+k_{10}^6+k_{11}^6)+2(k_{00}^4+k_{01}^4+k_{10}^4+k_{11}^4)+k_{00}^2(k_{00}^2+4)k_{01}k_{10}k_{11}\\
& &
+(4k_{00}^2+1)k_{01}^2k_{10}^2k_{11}^2-2k_{00}^2(k_{01}^4+k_{10}^4+k_{11}^4+4)-2k_{01}k_{10}k_{11}(k_{01}^2k_{10}^2k_{11}^2-1)+1,\\
\end{array}$$

$$\begin{array}{rcl}a_{0210}& =& 3k_{01}^2k_{10}k_{11}^4+3k_{01}^3k_{11}^3-2k_{01}k_{10}^2k_{11}^3
-k_{01}k_{11}^5+4k_{10}^3k_{11}^2+2k_{10}k_{11}^4+k_{01}^3k_{11}+k_{01}k_{10}^2k_{11}\\
& &
+k_{01}^2k_{10}+4k_{10}^3-2k_{10}k_{11}^2+3k_{01}k_{11}-2k_{10},\\
a_{0201}& =&
3k_{01}^2k_{11}k_{10}^4+3k_{01}^3k_{10}^3-2k_{01}k_{11}^2k_{10}^3
-k_{01}k_{10}^5+4k_{11}^3k_{10}^2+2k_{11}k_{10}^4+k_{01}^3k_{10}+k_{01}k_{11}^2k_{10}\\
& &
+k_{01}^2k_{11}+4k_{11}^3-2k_{11}k_{10}^2+3k_{01}k_{10}-2k_{11},\\
a_{0120}& =&
3k_{10}^2k_{01}k_{11}^4+3k_{10}^3k_{11}^3-2k_{10}k_{01}^2k_{11}^3
-k_{10}k_{11}^5+4k_{01}^3k_{11}^2+2k_{01}k_{11}^4+k_{10}^3k_{11}+k_{10}k_{01}^2k_{11}\\
& &
+k_{10}^2k_{01}+4k_{01}^3-2k_{01}k_{11}^2+3k_{10}k_{11}-2k_{01},\\
a_{0102}& =&
3k_{11}^2k_{01}k_{10}^4+3k_{11}^3k_{10}^3-2k_{11}k_{01}^2k_{10}^3
-k_{11}k_{10}^5+4k_{01}^3k_{10}^2+2k_{01}k_{10}^4+k_{11}^3k_{10}+k_{11}k_{01}^2k_{10}\\
& &
+k_{11}^2k_{01}+4k_{01}^3-2k_{01}k_{10}^2+3k_{11}k_{10}-2k_{01},\\
a_{0021}& =&
3k_{10}^2k_{11}k_{01}^4+3k_{10}^3k_{01}^3-2k_{10}k_{11}^2k_{01}^3
-k_{10}k_{01}^5+4k_{11}^3k_{01}^2+2k_{11}k_{01}^4+k_{10}^3k_{01}+k_{10}k_{11}^2k_{01}\\
& &
+k_{10}^2k_{11}+4k_{11}^3-2k_{11}k_{01}^2+3k_{10}k_{01}-2k_{11},\\
a_{0012}& =&
3k_{11}^2k_{10}k_{01}^4+3k_{11}^3k_{01}^3-2k_{11}k_{10}^2k_{01}^3
-k_{11}k_{01}^5+4k_{10}^3k_{01}^2+2k_{10}k_{01}^4+k_{11}^3k_{01}+k_{11}k_{10}^2k_{01}\\
& &
+k_{11}^2k_{10}+4k_{10}^3-2k_{10}k_{01}^2+3k_{11}k_{01}-2k_{10},\\
\end{array}$$

$$\begin{array}{rcl}
a_{1110} &= &
k_{11}(2(k_{00}^4-k_{00}^2k_{11}^2-k_{11}^4)+2k_{00}^2+k_{11}^2)+k_{01}k_{10}(3k_{11}^4-k_{00}^4+2k_{00}^2-k_{11}+2)\\
& &+k_{01}^2k_{10}^2k_{11}(4k_{11}^2+k_{00}^2+2)+4k_{01}^3k_{10}^3k_{11}^2,\\
a_{1101} &= &
k_{10}(2(k_{00}^4-k_{00}^2k_{10}^2-k_{10}^4)+2k_{00}^2+k_{10}^2)+k_{01}k_{11}(3k_{10}^4-k_{00}^4+2k_{00}^2-k_{10}+2)\\
& &+k_{01}^2k_{11}^2k_{10}(4k_{10}^2+k_{00}^2+2)+4k_{01}^3k_{11}^3k_{10}^2,\\
a_{1011} &= &
k_{01}(2(k_{00}^4-k_{00}^2k_{01}^2-k_{01}^4)+2k_{00}^2+k_{01}^2)+k_{11}k_{10}(3k_{01}^4-k_{00}^4+2k_{00}^2-k_{01}+2)\\
& &+k_{11}^2k_{10}^2k_{01}(4k_{01}^2+k_{00}^2+2)+4k_{11}^3k_{10}^3k_{01}^2,\\
\end{array}$$

$$\begin{array}{rcl}
b_{2000} &= &
k_{00}(2k_{00}^4+k_{01}^4+k_{10}^4+k_{11}^4+2k_{00}^2(k_{01}k_{10}k_{11}-2)+4k_{01}^2k_{10}^2k_{11}^2), \hspace{2.5cm}\\
b_{0200} & =& k_{00}(4k_{00}^2+2k_{01}^2+1+4k_{10}^2k_{11}^2),\\
b_{0020} & =& k_{00}(4k_{00}^2+2k_{10}^2+1+4k_{01}^2k_{11}^2),\\
b_{0002} & =& k_{00}(4k_{00}^2+2k_{11}^2+1+4k_{01}^2k_{10}^2),\\
\end{array}$$

$$\begin{array}{rcl}
b_{1100} &= &
-k_{00}(k_{01}(k_{01}^4+(k_{00}^2+3)(k_{01}^2+2k_{00}^2+4))-3(k_{01}^2+4)(k_{01}^2-2k_{00}^2-2)k_{10}k_{11}\\
& &\ -2k_{01}(k_{01}^2-2)k_{10}^2k_{11}^2,\\
b_{1010} &= &
-k_{00}(k_{10}(k_{10}^4+(k_{00}^2+3)(k_{10}^2+2k_{00}^2+4))-3(k_{10}^2+4)(k_{10}^2-2k_{00}^2-2)k_{01}k_{11}\\
& &\ -2k_{10}(k_{10}^2-2)k_{01}^2k_{11}^2,\\
b_{1001} &= &
-k_{00}(k_{11}(k_{11}^4+(k_{00}^2+3)(k_{11}^2+2k_{00}^2+4))-3(k_{11}^2+4)(k_{11}^2-2k_{00}^2-2)k_{01}k_{10}\\
& &\ -2k_{11}(k_{11}^2-2)k_{01}^2k_{10}^2,\\
\end{array}$$

$$\begin{array}{rcl}
b_{0110} &= &
k_{00}(k_{11}(5k_{11}^4+(k_{00}^2+2)(3k_{11}^2+1))-3(k_{11}^2+2)(k_{11}^2+1)k_{01}k_{10}, \hspace{2.3cm}\\
b_{0101} &= &
k_{00}(k_{10}(5k_{10}^4+(k_{00}^2+2)(3k_{10}^2+1))-3(k_{10}^2+2)(k_{10}^2+1)k_{01}k_{11},\\
b_{0011} &= &
k_{00}(k_{01}(5k_{01}^4+(k_{00}^2+2)(3k_{01}^2+1))-3(k_{01}^2+2)(k_{01}^2+1)k_{10}k_{11},\\
\end{array}$$
where the $k_\bullet$ are the coefficients of the equation (2.1) of
${\rm Kum}_X$. The $V_i$ ($i=01,\,10,\,11$) can be deduced from
$V_{00}$ by a suitable permutation of the coordinate functions
$y_\bullet$, namely the unique pairwise permutation that exchanges
$y_{00}$ and $y_i$. .
\end{prop}
\dem We make the same kind of calculations as in the case of the
proof of the Proposition 5.6, except that we use the formula 5.11
instead of the formula 5.10 and that we must use more data to
determine all the coefficients involved. $\square$\\

The computer the author used to carry out these computations was not
powerful enough to obtain the same result as in the case $p=5$.
Actually, one  could only check that the equation of the twisted
Kummer surface ${\rm Kum}_{X_1}$ did divide the image of the the
equation on the Kummer surface ${\rm Kum}_{X}$ via the ring
homomorphism defined by the above equations. \\

\section{Further questions}

{\bf Question 1 : } Notice (Remark 5.5) that the Theorem 6.1 of
[LP2] holds for any proper and smooth curve of genus 2 over an
algebraically closed field of characteristic 3. Therefore, one can
ask if the formulae given in Propositions 5.6 and 5.8 remain valid
for any smooth and proper curve in characteristic 5 and 7.

One way to answer this question is, being given any curve $X$ over
$k$,  to consider a family of genus 2 curve $\mathcal{X}$ over ${\rm
Spec}\,k[[t]]$ with sufficiently general generic fiber and special
fiber isomorphic to $X$ and to study how do the equations specialize
(see [LP2] or [Du] for examples of that method in characteristic 2).
For that purpose, we would need a more refined description of the
moduli space $\mathfrak{M}_K(k)$ of Kummer surfaces in $\P^3_k$, or
equivalently (see Remark 2.10), of the moduli space of smooth and
proper curve over $k$.

Note that, combining the chart 2.5 giving the expressions of the
$\omega(\tau)$ in terms of the $k_\bullet$ and the classification of
supersingular elliptic curve in small characteristic (see [H],
Chapter IV, Examples 4.23.1 and followings), one can determine the
open subset of $\mathfrak{M}_K(k)$ the closed points of which
correspond to curves with Prym varieties associated to \'{e}tale
double cover that all are ordinary. \\

\noindent {\bf Question 2 : } One would like to use the formulae
given is Proposition 5.6 and 5.8 to say if, as in the characteristic
3 case, the generalized Verschiebung base locus is reduced in
characteristics 5 and 7 (see [LnP] for a discussion on that topic).
We would like to study the singular locus of the surface $S$ in
Corollary 5.7, with the aim of showing that it is zero dimensional
and comparing it with the base locus of the Verschiebung. We plan to
use the computers of the MEDICIS at the Ecole Polytechnique to
carry out these computations.\\

\noindent {\bf Question 3 : } The most interesting question arisen
by this work is how can we recover the results obtained by
computational means in a more geometric way. We plan to return to
this latter question in a future work.

\newpage

\vspace{1cm}

\noindent \Large \textbf{References}\\
\normalsize
\begin{enumerate}

\item[{[B]}] \textsc{A. Beauville}: \emph{Fibrés de rang 2 sur une courbe,
fibrés déterminant et fonctions thêta}, Bull. Soc. Math. France
{\bf116}, 1988, 431--448.

\item[{[Du]}] \textsc{L. Ducrohet}: \emph{The action of the Frobenius map on rank 2 vector bundles over a supersingular curve in characteristic 2},
arXiv : math.AG/0504500.

\item[{[vG]}] \textsc{B. van Geemen}: \emph{Schottky-Jung
relations and vector bundles on hyperelliptic curves},
 Math. Ann. {\bf281} (1988), 431--449.

\item[{[GD]}] \textsc{M. Gonz\'alez-Dorrego}: \emph{(16-6)-configurations and Geometry of Kummer surfaces in $\P^3$},
Memoirs of the American Math. Society, Vol {\bf107}, 1994.

\item[{[H]}] \textsc{R. Hartshorne}: \emph{Algebraic geometry},
Graduate Texts in Mathematics {\bf52}, Springer, New-York, 1977.

\item[{[Hud]}] \textsc{R. Hudson}: \emph{Kummer's quartic surface},
Cambridge Univ. Press, 1905.

\item[{[JX]}] \textsc{K. Joshi, E.Z. Xia}: \emph{Moduli of vector bundles on curves in positive characteristic},
Compositio Math. {\bf122}, N°3 (2000), 315--321.

\item[{[LP1]}] \textsc{Y. Laszlo}, \textsc{C. Pauly}: \emph{The
action of the Frobenius map on rank 2 vector bunbles in
characteristic 2}, J. of Alg. Geom. {\bf11} (2002), 219--243.

\item[{[LP2]}] \textsc{Y. Laszlo}, \textsc{C. Pauly}: \emph{The
Frobenius map, rank 2 vector bunbles and Kummer's quartic surface in
characteristic 2 and 3}, Advances in Mathematics {\bf185} (2004),
246--269.

\item[{[LnP]}] \textsc{H. Lange}, \textsc{C. Pauly}: \emph{On
Frobenius-destabilized rank-2 vector bundles over curves}, arXiv :
math.AG/0309456 (2003)

\item[{[LS]}] \textsc{H. Lange}, \textsc{U. Stuhler}:
\emph{Vektorb\"{u}ndel auf Kurven und Darstellungen der
algebraischen Fundamentalgruppe}, Math. Zeit. {\bf156} (1977),
73--83.

\item[{[Mu1]}] \textsc{D. Mumford}:  \emph{Abelian varieties},
Tata Institute of Fundamental Research Studies in Mathematics {\bf
5}, Bombay, 1970

\item[{[Mu2]}] \textsc{D. Mumford}:  \emph{On equations defining
abelian varieties. I.}, Invent. Math. {\bf1} (1966), 287--354.

\item[{[Mu3]}] \textsc{D. Mumford}:  \emph{Prym varieties. I.},
Contributions to analysis, 325--350, London, New York Academic
Press, 1974.

\item[{[NR]}] \textsc{M.S. Narasimhan, S. Ramanan}: \emph{Moduli
of vector bundles on a compact Riemann surface}, Ann. of Math.
{\bf89}
 (1969), 14--51.

\item[{[R]}] \textsc{M. Raynaud}: \emph{Sections des fibr\'es
vectoriels sur une courbe}, Bull. Soc. Math. France {\bf110} (1982),
103--125.

\item[{[Sek]}] \textsc{T. Sekiguchi}: \emph{On projective
normality of abelian varieties. II.}, J. Math. Soc. Japan {\bf29}
(1977), 709--727.

\item[{[Sil]}] \textsc{J. H. Silverman}: \emph{The  arithmetic of elliptic curves} , Graduate Texts in
Math. {\bf106}, Springer-Verlag, New-York, 1986.

\item[{[Zh]}] \textsc{B. Zhang}: \emph{Revêtements étales abeliens de
courbes génériques et ordinarité}, Ann. Fac. Sci. Toulouse Sér. 6,
1992, 133--138.

\end{enumerate}

\vspace{2cm}

\noindent Laurent Ducrohet\\
Universit\'e Pierre et Marie Curie\\
Analyse alg\'ebrique, UMR 7586\\
4, place Jussieu\\
75252 Paris Cedex 05 France\\
e-mail: ducrohet@math.jussieu.fr

\end{document}